\documentclass[10pt]{amsart}

\newif\iffinalrun

\usepackage{hyperref}
\usepackage{amssymb,eucal}
\usepackage[all,cmtip,poly]{xy}
\usepackage{enumitem}
\usepackage{color}
\usepackage{appendix}

\newcommand{\WD}{{\operatorname{WD}}}
\newcommand{\rec}{{\operatorname{rec}}}
\newcommand{\Fpbar}{\overline{\F}_p}
\newcommand{\Fp}{\F_p}

\newcommand{\mf}{\mathfrak}
\newcommand{\m}{\mathfrak{m}}
\renewcommand{\a}{\mathfrak{a}}
\newcommand{\p}{\mathfrak{p}}

\newcommand{\rbar}{\overline{r}}
\newcommand{\taubar}{\overline{\tau}}

\newcommand{\GQp}{{G_{\Qp}}}

\newcommand{\T}{\mathbb{T}}

\newcommand{\cA}{\mathcal{A}}
\newcommand{\cB}{\mathcal{B}}

\newcommand{\cD}{\mathcal{D}}

\newcommand{\cF}{\mathcal{F}}
\newcommand{\cG}{\mathcal{G}}

\newcommand{\cM}{\mathcal{M}}
\newcommand{\cN}{\mathcal{N}}
\newcommand{\cO}{\mathcal{O}}
\newcommand{\cP}{\mathcal{P}}

\newcommand{\cS}{\mathcal{S}}
\newcommand{\cT}{\mathcal{T}}

\newcommand{\cX}{\mathcal{X}}

\newcommand{\fM}{\mathfrak{M}}

\newcommand{\A}{\mathbb A}
\newcommand{\C}{\mathbb C}
\newcommand{\F}{\mathbb F}
\newcommand{\fp}{\F_p}

\newcommand{\fpb}{\overline \F_p} 
\newcommand{\fpn}[1]{\F_{p^{#1}}}

\renewcommand{\O}{\mathcal O}

\renewcommand{\oe}{\O_E}
\newcommand{\Q}{\mathbb Q}
\newcommand{\qb}{\overline\Q}
\newcommand{\qp}{\Q_p}

\newcommand{\qpb}{\overline\Q_p} 
\newcommand{\R}{\mathbb R}

\newcommand{\Z}{\mathbb Z}
 
\newcommand{\zp}{\Z_p}

\newcommand{\Zp}{\Z_p}

\newcommand{\into}{\hookrightarrow}

\newcommand{\onto}{\twoheadrightarrow}
\newcommand{\congto}{\xrightarrow{\,\sim\,}}
\newcommand{\isoto}{\congto}

\newcommand{\s}{^\times}

\newcommand{\ab}{^{\mathrm{ab}}}

\newcommand{\cont}{^{\mathrm{cont}}}

\newcommand{\rhobar}{\overline{\rho}}

\newcommand{\norm}[1]{\left\vert#1\right\vert} 

\newcommand{\GL}{\operatorname{GL}}

\newcommand{\Fil}{\operatorname{Fil}}

\newcommand{\ad}{\operatorname{ad}}

\newcommand{\Qp}{\Q_p}

\newcommand{\Qpbar}{\overline{\Q}_p}
\newcommand{\Zpbar}{\overline{\Z}_p}

\newcommand{\M}{\mathcal{M}}

\newcommand{\gr}{\operatorname{gr}} 
\newcommand{\FBrModdd}[1][r]{\text{$\F$-$\operatorname{BrMod}_{\mathrm{dd}}^{#1}$}}
\newcommand{\FBrModddN}[1][r]{\text{$\F$-$\operatorname{BrMod}_{\mathrm{dd},0}^{#1}$}}
\newcommand{\FBrModN}[1][r]{\text{$\F$-$\operatorname{BrMod}_{0}^{#1}$}}
\newcommand{\FBrMod}[1][r]{\text{$\F$-$\operatorname{BrMod}^{#1}$}}

\newcommand{\BrMod}[1][r]{\text{$\operatorname{BrMod}^{#1}$}}
\newcommand{\BrModN}[1][r]{\text{$\operatorname{BrMod}_0^{#1}$}}
\newcommand{\BrModdd}[1][r]{\text{$\operatorname{BrMod}_{\mathrm{dd}}^{#1}$}}
\newcommand{\BrModddN}[1][r]{\text{$\operatorname{BrMod}_{\mathrm{dd},0}^{#1}$}}
\newcommand{\OEModdd}[1][r]{\text{$\cO_E$-$\Mod_{\mathrm {dd}}^{#1}$}}

\DeclareMathOperator{\Mod}{Mod}
\DeclareMathOperator{\Art}{Art}

\newcommand{\ord}{\mathrm{ord}_p}

\DeclareMathOperator{\End}{End}
\DeclareMathOperator{\Hom}{Hom}
\DeclareMathOperator{\Spec}{Spec}

\DeclareMathOperator{\rank}{rank}

\DeclareMathOperator{\Gal}{Gal}

\DeclareMathOperator{\tr}{tr}

\DeclareMathOperator{\res}{res}

\DeclareMathOperator{\Ind}{Ind}

\DeclareMathOperator{\red}{{red}}
\DeclareMathOperator{\rad}{rad}

\DeclareMathOperator{\im}{im}
\DeclareMathOperator{\Frob}{Frob}
\DeclareMathOperator{\Sym}{Sym}
\DeclareMathOperator{\supp}{supp}

\newcommand{\plim}{\varprojlim}  
\newcommand{\ilim}{\varinjlim}   

\renewcommand{\o}[1]{\overline{#1}}
\newcommand{\wt}[1]{\widetilde{#1}}
\newcommand{\wh}[1]{\widehat{#1}}

\newcommand{\N}{\mathbb{N}}

\newcommand{\st}{\mathrm{st}}

\newcommand{\textD}{\mathrm{D}}

\newcommand{\textT}{\mathrm{T}}

\newcommand{\Dst}{\textD_{\st}}

\newcommand{\Tst}{\textT_{\st}}
\newcommand{\Tcris}{\textT_{\mathrm{cris}}}

\newcommand{\adj}{\mathrm{adj}}
\newcommand{\un}[1]{\mathrm{nr}_{#1}}
\newcommand{\und}[1]{\underline{#1}}

\newcommand{\maq}[9]{\left(\begin{array}{ccc} #1&#2&#3\\#4&#5&#6\\#7&#8&#9 \end{array}\right)}
\newcommand{\maqdue}[4]{\left(\begin{array}{cc} #1&#2\\#3&#4 \end{array}\right)}

\newcommand{\defeq}{\stackrel{\mathrm{\text{\tiny{\upshape def}}}}{=}} \newcommand{\barS}{\overline{{S}}}
\newcommand{\teich}[1]{\widetilde{#1}}
\newcommand{\columnvct}[1]{\left(\begin{array}{c} #1\end{array}\right)}

\newcommand{\bA}{\mathbb{A}}
\newcommand{\bT}{\mathbb{T}}

\newcommand{\mat}{\mathrm{Mat}}
\newcommand{\glob}{_\mathrm{glob}}
\newcommand{\FL}{\mathrm{FL}}

\iffinalrun
  \newcommand{\need}[1]{}
  \newcommand{\mar}[1]{}
\else
  \newcommand{\need}[1]{{\tiny *** #1}}
  \newcommand{\mar}[1]{\marginpar{\raggedright\tiny #1}}
\fi

\renewcommand{\(}{\textup{(}}
\renewcommand{\)}{\textup{)}}

\usepackage{hyperref}
    \usepackage{enumitem}

\theoremstyle{plain} 

\newtheorem{lm}[equation]{Lemma}

\newtheorem{prop}[equation]{Proposition}
\newtheorem{thm}[equation]{Theorem}
\newtheorem{corollary}[equation]{Corollary}

\newtheorem{question}[equation]{Question}
\newtheorem{claim}[equation]{Claim}
\newtheorem{ithm}{Theorem}

\newtheorem{iprop}[ithm]{Proposition}

\theoremstyle{definition}
\newtheorem{df}[equation]{Definition}
\theoremstyle{remark}
\newtheorem{rk}[equation]{Remark}

\theoremstyle{remark} 
\newtheorem{remark}[equation]{Remark}

\numberwithin{equation}{subsection}
\numberwithin{figure}{subsection}

\begin{document}

\title[On mod $p$ local-global compatibility for $\GL_3$]{On mod $p$ local-global compatibility for $\GL_3$ in the ordinary case }
\author{Florian Herzig}
\address{Department of Mathematics,
University of Toronto,
40 St. George Street,
Toronto, ON M5S 2E4, Canada}
\email{herzig@math.toronto.edu}
\thanks{The first-named author was partially supported by a Sloan Fellowship and an NSERC grant}
\author{Daniel Le}
\address{Department of Mathematics,
University of Toronto,
40 St. George Street,
Toronto, ON M5S 2E4, Canada}
\email{le@math.toronto.edu}
\author{Stefano Morra}
\address{Institut Montpelli\'erain Alexander Grothendieck,
Universit\'e de Montpellier,
Case courrier 051, Place Eug\`ene Bataillon,
34095 Montpellier Cedex, France}
\email{stefano.morra@umontpellier.fr}
\thanks{}

\maketitle

\begin{abstract}
  Suppose that $F/F^+$ is a CM extension of number fields in which the prime $p$ splits completely and every other prime is unramified.
  Fix a place $w| p$ of $F$.
  Suppose that $\rbar : \Gal(\o F/F) \to \GL_3(\fpb)$ is a continuous irreducible Galois representation such that $\rbar|_{\Gal(\o F_w/F_w)}$ is upper-triangular, maximally non-split, and generic.
  If $\rbar$ is automorphic, and some suitable technical conditions hold, we show that $\rbar|_{\Gal(\o F_w/F_w)}$ can be recovered from the $\GL_3(F_w)$-action on a space of mod $p$ automorphic forms on a compact unitary group.
  On the way we prove results about weights in Serre's conjecture for $\rbar$, show the existence of an ordinary lifting of $\rbar$, and prove the freeness of certain Taylor--Wiles patched modules in this context.
  We also show the existence of many Galois representations $\rbar$ to which our main theorem applies.
\end{abstract}

\tableofcontents

\section{Introduction}

\subsection{Motivation and statement of main results}
\label{sec:motiv-stat-main}

Suppose that $p$ is a prime and that $\rhobar : \Gal(\qpb/\qp) \to \GL_n(\fpb)$ is a continuous Galois representation.  One
would hope that there is a mod $p$ Langlands correspondence that associates to $\rhobar$ in a natural way a smooth
representation $\Pi(\rhobar)$ of $\GL_n(\qp)$ over $\fpb$ (or maybe a packet of such representations), and similarly for
$p$-adic representations. Unfortunately, at this point, this is only known for $n \le 2$ \cite{breuilI}, \cite{Colmez}.  But
suppose now that $F/F^+$ is a CM extension of number fields in which $p$ splits completely, and fix a place $w| p$.  Even
in the absence of a local mod $p$ Langlands correspondence for $n > 2$, given a \emph{global} automorphic Galois representation
$\rbar : \Gal(\o F/F) \to \GL_n(\fpb)$ we can define smooth representations $\Pi\glob(\rbar)$ of $\GL_n(\qp)$ over $\fpb$ on
spaces of mod $p$ automorphic forms on a definite unitary group, that serve as candidates for $\Pi(\rbar|_{\Gal(\o
  F_w/F_w)})$ (in the spirit of Emerton's local-global compatibility \cite{emerton-local-global}; see also \cite{CEGGPS} in the
$p$-adic setting).  It is not known whether
$\Pi\glob(\rbar)$ depends only on $\rbar|_{\Gal(\o F_w/F_w)}$.  The motivating question of this paper is opposite to this:
do the candidate representations $\Pi\glob(\rbar)$ contain at least as much information
as $\rbar|_{\Gal(\o F_w/F_w)}$? We answer this question in the affirmative in many cases when $n = 3$.

We fix a finite extension $E/\qp$ with residue field $\F$, and consider absolutely irreducible Galois representations $\rbar :
\Gal(\o F/F) \to \GL_3(\F)$. We assume moreover that $\rbar|_{\Gal(\o F_w/F_w)}$ is upper-triangular, maximally non-split, and generic. This means that
\begin{equation}\label{eq:8}
  \rbar|_{\Gal(\o F_w/F_w)}\sim \maq{\omega^{a+1}\un{\mu_2}}{\ast_1}{\ast}{}{\omega^{b+1}\un{\mu_1}}{\ast_2}{}{}{\omega^{c+1}\un{\mu_0}},
\end{equation}
the extensions $\ast_1$, $\ast_2$ are non-split, and $a-b > 2$, $b-c > 2$, $a-c < p-3$.
(Here, $\omega$ is the mod $p$ cyclotomic character and $\un{\mu}$ denotes the unramified character taking value $\mu \in \F\s$ on
geometric Frobenius elements.)
It is not hard to see that once the diagonal characters are fixed, the isomorphism class of $\rbar|_{\Gal(\o F_w/F_w)}$ is
determined by an invariant $\FL(\rbar|_{\Gal(\o F_w/F_w)})$ that can take any value in $\mathbb{P}^1(\F) \setminus \{\mu_1\}$.
(We normalize this invariant using Fontaine--Laffaille theory, see Definition~\ref{definvFL}.)

To explain our results, we briefly describe our global setup, referring to Section~\ref{sec:local-glob-comp} for details.
Fix a unitary group $G_{/F^+}$ such that $G \times F \cong \GL_3$ and $G(F^{+}_{v})\cong U_3(\R)$ for all $v|\infty$.
Choose a model $\cG_{/\cO_{F^+}}$ of $G$ such that $\cG \times \cO_{F^+_{v'}}$ is reductive for all places $v'$ of $F^+$ that
split in $F$. Let $v = w|_{F^+}$. Choose a compact open subgroup $U = U^v \times \cG(\cO_{F^+_v}) \le G(\A_{F^+}^\infty)$ that is sufficiently small,
and unramified at all places dividing $p$. Let $V'$ denote an irreducible smooth representation of $\prod_{v'| p, v'\ne v} \cG(\cO_{F^+_{v'}})$ over
$\F$ determined by a highest weight in the lowest alcove, and let $\wt V'$ denote a Weyl module of $\prod_{v'| p, v'\ne v} \cG(\cO_{F^+_{v'}})$ over $\cO_E$ 
of the same highest weight, so $\wt V'\otimes_{\cO_E} \F \cong V'$. 
(Except for Theorem~\ref{intro:thm-all-invariants-arise} below, the reader may assume for simplicity that $V' = \F$ and $\wt V' = \cO_E$.)
We can then define in the usual way spaces of mod $p$ automorphic forms $S(U^v, V') = \ilim_{U_v \le G(F_v^+)} S(U^v U_v, V')$ and similarly $S(U^v, \wt V')$
that are smooth representations of $G(F_v^+) \cong \GL_3(\qp)$ (where this isomorphism depends on our chosen place $w | v$).

We fix a cofinite subset $\cP$ of places $w' \nmid p$ of $F$ that split over $F^+$, such that $U$ is unramified at $w'|_{F^+}$, and such that $\rbar$ is unramified at $w'$.
Then the abstract Hecke algebra $\T^{\cP}$ generated over $\cO_E$ by Hecke operators at all places in $\cP$ acts on $S(U^v, V')$ and $S(U^v, \wt V')$, commuting with the $\GL_3(\qp)$-actions.
Moreover, $\rbar$ determines a maximal ideal $\m_{\rbar}$ of $\T^{\cP}$. We assume that $\rbar$ is automorphic in this setup, which
means that $S(U^v, V')[\m_{\rbar}] \ne 0$ (or equivalently, $S(U^v, V')_{\m_{\rbar}} \ne 0$).
These $\GL_3(\qp)$-representations, $S(U^v, V')[\m_{\rbar}]$, are the natural candidates $\Pi\glob(\rbar)$ that we mentioned above, at least if the level $U^v$ is chosen optimally. 
It is a consequence of earlier work of the first-named author and C.\ Breuil \cite{BH} that the $\GL_3(\qp)$-representation $S(U^v, V')[\m_{\rbar}]$ determines the ordered triple of diagonal characters $(\omega^{a+1}\un{\mu_2}, \omega^{b+1}\un{\mu_1}, \omega^{c+1}\un{\mu_0})$ of $\rbar|_{\Gal(\o F_w/F_w)}$. 
In fact, the triple $(a,b,c)$ can be recovered from the (ordinary part of the) $\GL_3(\zp)$-socle -- i.e., the Serre weights of $\rbar$ -- by \cite{gee-geraghty} and the $\mu_i \in \F\s$ are determined by the Hecke action at $p$ on the $\GL_3(\zp)$-socle.
It therefore remains for us to show that $S(U^v, V')[\m_{\rbar}]$ determines the invariant $\FL(\rbar|_{\Gal(\o F_w/F_w)})$.
(We note that the representation $\Pi(\rbar|_{\Gal(\o F_w/F_w)})^{\mathrm{ord}}$ of \cite{BH} does not contain this information.)

Let $I$ denote the Iwahori subgroup of $\GL_3(\zp)$, which is the preimage of the upper-triangular matrices $B(\fp)$ in $\GL_3(\fp)$.
If $V$ is a representation of $\GL_3(\zp)$ over $\cO_E$ and $a_i \in \Z$ we write 
\[ V^{I,(a_2,a_1,a_0)} \defeq \Hom_I(\cO_E(\teich{\omega}^{a_2}\otimes\teich{\omega}^{a_1}\otimes\teich{\omega}^{a_0}),
V), \] where the character in the domain denotes the inflation to $I$ of the homomorphism $B(\fp) \to \cO_E\s$, $\text{$\tiny \maq
xyz{}uv{}{}w$} \mapsto \wt x^{a_2}\wt u^{a_1}\wt w^{a_0}$. If $V$ is even a representation of $\GL_3(\qp)$, then 
$V^{I,(a_2,a_1,a_0)}$ affords an action of $U_p$-operators $U_1$, $U_2$ (see \eqref{eq:6}).
Define also $\Pi \defeq \text{$\tiny\maq{}{1}{}{}{}{1}{p}{}{}$}$, which sends $V^{I,(a_2,a_1,a_0)}$ to $V^{I,(a_1,a_0,a_2)}$.

Finally, and crucially, we define explicit group algebra operators $S$, $S' \in \F[\GL_3(\fp)]$, see \eqref{eq:7}. We can now state our first main theorem.

\begin{ithm}[Thm.\ \ref{main local global}]
\label{intro:main local global}
We make the following additional assumptions:
\begin{enumerate}
\item $\FL(\rbar|_{\Gal(\o F_w/F_w)})\not\in\{0,\infty\}$.
\item The $\cO_{E}$-dual of $S(U^v,\wt V')_{\m_{\rbar}}^{I,(-b,-c,-a)}$ is free over $\T$, where 
  $\T$ denotes the $\cO_{E}$-subalgebra of $\End\big(S(U^v,\wt V')_{\m_{\rbar}}^{I,(-b,-c,-a)}\big)$ generated by $\bT^{\cP}$, $U_1$, and $U_2$.
\end{enumerate}
Then we have the equality
\begin{equation*}
S'\circ \Pi =(-1)^{a-b}\cdot\frac{b-c}{a-b}\cdot\FL(\rbar|_{\Gal(\o F_w/F_w)})\cdot {S}
\end{equation*}
of maps
\[ S(U^v,V')[\m_{\rbar}]^{I,(-b,-c,-a)}[U_1,U_2] \to S(U^v,V')[\m_{\rbar}]^{I,(-c-1,-b,-a+1)}.  \]
Moreover, these maps are injective with non-zero domain.
In particular, $\FL(\rbar|_{\Gal(\o F_w/F_w)})$ is determined by the smooth $\GL_3(\Qp)$-representation $S(U^v,V')[\m_{\rbar}]$.
\end{ithm}

The first assumption is related to the surprising fact that the $\GL_3(\zp)$-socle of $S(U^v,V')[\m_{\rbar}]$ changes for the exceptional
two invariants, see Theorem~\ref{intro:theorem weight elimination} below. 
(Incidentally, this means that even in the exceptional cases $\FL(\rbar|_{\Gal(\o F_w/F_w)})$ is determined by $S(U^v,V')[\m_{\rbar}]$.)

We show that the second assumption is often a consequence of the first assumption.
The second assumption is an analogue of Mazur's ``mod $p$ multiplicity one'' result, and thus our result may be of independent interest.
We have not tried to optimize our hypotheses, the most stringent of which is that $U$ may be taken to be unramified at all finite places.

\begin{ithm}[Thm.\ \ref{freeT}]
\label{intro:freeT}
Assume hypotheses \ref{field unram}--\ref{automorphy} in Section \ref{sec:setup}.
If we have $\FL(\rbar|_{\Gal(\o F_w/F_w)}) \ne \infty$, then assumption \(ii\) in Theorem~\ref{intro:main local global} holds.
\end{ithm}

In fact we show that for any value of $\FL(\rbar|_{\Gal(\o F_w/F_w)})$ either the second assumption or its dual holds (see Remark \ref{rk:dual-lattice}).

We also show, using results of \cite{EG}, that for any given local Galois representation as in \eqref{eq:8} we can construct a globalization to which Theorem~\ref{intro:freeT} applies.

\begin{ithm}[Thm.\ \ref{thm:all-invariants-arise}]
  \label{intro:thm-all-invariants-arise}
  Suppose that $\rhobar : \Gal(\qpb/\qp) \to \GL_3(\F)$ is upper-triangular, maximally non-split, and generic. Then, after possibly replacing $\F$ by a
  finite extension, there exist a CM field $F$, a Galois representation $\rbar : \Gal(\o F/F) \to \GL_3(\F)$, a place $w|p$ of $F$,
  groups $G_{/F^+}$ and $\cG_{/\cO_{F^+}}$, and a compact open subgroup $U^v$ \(where $v = w|_{F^+}$\) satisfying all hypotheses
  of the setup in Section \ref{sec:setup} such that $\rbar|_{\Gal(\o F_w/F_w)} \cong \rhobar$.
  In particular, if $\FL(\rhobar) \not \in \{0,\infty\}$, Theorem~\ref{intro:main local global} applies to $\rbar$.
\end{ithm}

As a by-product of our methods we almost completely determine the set of Serre weights of $\rbar$.
Here, the set $W_w(\rbar)$ is defined to be the set of irreducible $\GL_3(\zp)$-representations whose duals occur in the $\GL_3(\zp)$-socle
of $S(U^v, V')_{\m_{\rbar}}$ (for some $U^v$ and $\cP$ as above). See Section \ref{sec:serre weights} for our notation for Serre weights.

\begin{ithm}[Thm.\ \ref{theorem weight elimination}]
\label{intro:theorem weight elimination}
Keep the assumptions on $\rbar$ that precede Theorem~\ref{intro:main local global} above.
\begin{enumerate}
\item If $\FL(\rbar\vert_{\Gal(\o F_w/F_w)})\notin\{0,\infty\}$ 
we have
\begin{align*}
&\{F(a-1,b,c+1)\}\subseteq W_{w}(\rbar)\subseteq\\
&\qquad\subseteq \{F(a-1,b,c+1),\,F(c+p-1,b,a-p+1)\}.
\end{align*}
\item If $\FL(\overline{r}\vert_{\Gal(\o F_w/F_w)})=\infty$ we have
\begin{align*}
&\{F(a-1,b,c+1), F(a,c,b-p+1)\}\subseteq W_{w}(\rbar)\subseteq\\
&\qquad\subseteq \{F(a-1,b,c+1),\,F(c+p-1,b,a-p+1), F(a,c,b-p+1)\}.
\end{align*}
\item  If $\FL(\overline{r}\vert_{\Gal(\o F_w/F_w)})=0$ we have
\begin{align*}
&\{F(a-1,b,c+1), F(b+p-1,a,c)\}\subseteq W_{w}(\rbar)\subseteq\\
&\qquad\subseteq \{F(a-1,b,c+1),\,F(c+p-1,b,a-p+1), F(b+p-1,a,c)\}.
\end{align*}
\end{enumerate}
\end{ithm}

This is one of the first Serre weight results in dimension 3. It was completed in early 2014, before the recent progress of \cite{LLLM}
(on Serre weights in dimension 3 in the generic semisimple case, using different methods).
As we said above, the dependence on $\FL(\rbar\vert_{\Gal(\o F_w/F_w)})$ in this theorem was unexpected, as there were no explicit Serre weight conjectures in the literature that apply to non-semisimple $\rbar|_{\Gal(\o F_w/F_w)}$.

As a consequence of this theorem we also show the existence of an automorphic lift $r$ of $\rbar$ such that $r\vert_{\Gal(\o F_w/F_w)}$ is upper-triangular.
It is in the same spirit as the main results of \cite{BLGG12} (which concerned two-dimensional representations).

\begin{ithm}[Cor.\ \ref{corollary weight elimination}]
\label{intro:corollary weight elimination}
In the setting of Theorem \ref{intro:theorem weight elimination}, $\rbar$ has an automorphic lift $r:\Gal(\o F/F)\to \GL_3(\cO_E)$ \(after possibly enlarging $E$\) such that $r|_{\Gal(\o F_w/F_w)}$ is crystalline and ordinary of Hodge--Tate weights $\{-a-1,-b-1,-c-1\}$.
\end{ithm}

\subsection{Methods used}
\label{sec:methods-used}

Theorems~\ref{intro:main local global} and \ref{intro:freeT} generalize earlier work of Breuil--Diamond \cite{BD} which treated two-dimensional Galois representations of $\Gal(\o F/F)$, where $F$ is totally real and $p$ is unramified in $F$. 
We follow the same general strategy: we lift the Hecke eigenvalues of $\rbar$ to a well-chosen type in characteristic zero, use classical local-global compatibility at $p$, and then study carefully how both the Galois-side and the $\GL_3$-side reduce modulo $p$.
However, it is significantly more difficult to carry out this strategy in the $\GL_3$-setting.

We first prove the upper bound in Theorem~\ref{intro:theorem weight elimination} by lifting to various types in characteristic zero and using integral $p$-adic Hodge theory to reduce modulo $p$.
This is more involved in dimension 3, since we are no longer in the potentially Barsotti--Tate setting.
We crucially use results of Caruso to filter our Breuil module (corresponding to $\rbar\vert_{\Gal(\o F_w/F_w)}$) according to the socle filtration on $\rbar\vert_{\Gal(\o F_w/F_w)}$, see Proposition~\ref{Proposition7Note Florian}.

The following theorem is our key local result on the Galois-side. 
Our chosen type is a tame principal series that contains in its reduction mod $p$ all elements of $W_w(\rbar)$ (unlike in \cite{BD}, where the intersection always consisted of one element).
In contrast to \cite{BD} we get away with a rough classification of the strongly divisible module corresponding to $\rho$. (We do not need to determine Frobenius and  monodromy operators.)
We also note that the relevant information on the Galois-side is independent of the Hodge filtration, so that we can transfer this information to the $\GL_3$-side using classical local-global compatibility.

\begin{ithm}[Thm.\ \ref{main theorem Galois}]
\label{intro:main theorem Galois}
Let $\rho: \Gal(\qpb/\qp)\rightarrow \mathrm{GL}_3(\oe)$ be a potentially semistable $p$-adic Galois representation of Hodge--Tate weights $\{-2,-1,0\}$ and inertial type 
$\teich{\omega}^{a}\oplus\teich{\omega}^{b}\oplus\teich{\omega}^{c}$.
Assume that the residual representation $\rhobar : \Gal(\qpb/\qp)\rightarrow \mathrm{GL}_3(\F)$ is upper-triangular, maximally non-split, and generic as in \eqref{eq:8}.
Let $\lambda\in \oe$ be the Frobenius eigenvalue on $\Dst^{\Qp,2}(\rho)^{I_{\Qp} = \teich{\omega}^{b}}$.
Then the Fontaine--Laffaille invariant of $\rhobar$ is given by:
$$
\FL(\rhobar)=\red\left(p\lambda^{-1}\right),
$$
where $\red$ denotes the specialization map $\mathbb{P}^1(\oe)\rightarrow \mathbb{P}^1(\F)$.
\end{ithm}

On the $\GL_3$-side our main innovation consists of the explicit group algebra operators $S$, $S'$. The analogues of these operators for $\GL_2$ show up in various contexts 
(see, for example, \cite{paskunas-restriction}, Lemma 4.1, \cite{BP}, Lemma 2.7, \cite{breuil-IL}, \S4, and \cite{BD}, Proposition 2.6.1). Our proof is specific to $\GL_3$. It would be interesting to find a more conceptual explanation for them. See also Question~\ref{qu:group-algebra-oper} for a further discussion of such operators.

\begin{iprop}[Prop.\ \ref{main mod p GL3}]
\label{intro:main mod p GL3}\ 
\begin{enumerate}
\item There is a unique non-split extension of irreducible $\GL_3(\fp)$-representations
$$
0\rightarrow F(-c-1,-b,-a+1)\rightarrow V\rightarrow F(-b+p-1,-c,-a)\rightarrow0
$$
and $S$ induces an isomorphism $S:\,V^{I,(-b,-c,-a)}\stackrel{\sim}{\longrightarrow}V^{I,(-c-1,-b,-a+1)}$
of one-dimensional vector spaces.
\item There is a unique non-split extension of irreducible $\GL_3(\fp)$-representations
$$
0\rightarrow F(-c-1,-b,-a+1)\rightarrow V\rightarrow F(-c,-a,-b-p+1)\rightarrow0
$$
and $S'$ induces an isomorphism $S':\,V^{I,(-c,-a,-b)}\stackrel{\sim}{\longrightarrow}V^{I,(-c-1,-b,-a+1)}$ of one-dimensional vector spaces.
\end{enumerate}
\end{iprop}

The reduction mod $p$ result on the $\GL_3$-side is comparatively easier, see Proposition \ref{main char 0 GL3}.

By combining the above results we deduce Theorem~\ref{intro:main local global}.
We note that assumption (ii) is needed for lifting elements of $S(U^v,V')[\m_{\rbar}]^{I,(-b,-c,-a)}[U_1,U_2]$ to suitable Iwahori eigenvectors in characteristic zero.
The $U_i$-operators allow us to deal with the possible presence of the shadow weight $F(c+p-1,b,a-p+1)$ in Theorem~\ref{intro:theorem weight elimination}.
(The term ``shadow weight'' is defined in \cite{EGH}, \S6 and more generally in \cite{GHS}, \S1.5 and \S7.2.)
Namely, if $v \in S(U^v,V')[\m_{\rbar}]^{I,(-b,-c,-a)}[U_1,U_2]$ is non-zero, we show that it generates the representation of Proposition~\ref{intro:main mod p GL3}(i) under the $\GL_3(\zp)$-action.
Similar comments apply to $\Pi v$.
Proposition~\ref{intro:main mod p GL3} then allows us to deduce that the maps in Theorem~\ref{intro:main local global} are well-defined and injective.
(We refer to Remark \ref{rk:main-local-global1} for variations on assumption (ii).
The stronger assumptions appearing in Remark \ref{rk:main-local-global2} are analogous to the multiplicity one conditions appearing in \cite{BD}.) 

Interestingly, the argument proving Theorem~\ref{intro:main local global} also lets us deduce the hardest part of Theorem~\ref{intro:theorem weight elimination}, namely the existence of the shadow weights $F(a,c,b-p+1)$, $F(b+p-1,a,c)$ in the two exceptional cases.
After \cite{EGH}, this is only the second result in the literature proving the existence of shadow weights. (Again this precedes \cite{LLLM}.)

Finally, we establish Theorem~\ref{intro:freeT}. As in \cite{BD} our method relies on the Taylor--Wiles method. However, as we do not know whether
our local deformation ring at $p$ is formally smooth (which in any case should be false if our chosen type intersects $W_w(\rbar)$ in more than one element)
we cannot directly apply Diamond's method \cite{diamond}. Instead we use the patched modules of \cite{CEGGPS} that live over the universal deformation
space at $p$ and use ideas of \cite{EGS} and \cite{le}. See Theorem~\ref{freeR} for our freeness result at infinite level from which we deduce Theorem~\ref{intro:freeT}. 
Similarly to above, we add $U_p$-operators in order to deal with the possible presence of the shadow weight $F(c+p-1,b,a-p+1)$.

\subsection{Acknowledgements}
\label{sec:acknowledgements}

The debt that this paper owes to the work of Breuil--Diamond \cite{BD} should be obvious to the reader.
We thank Christophe Breuil for a number of very helpful conversations as well as for his encouragement.
We thank Guy Henniart for providing a reference for the proof of Lemma~\ref{lm:inertial-llc-principal-series},
Xavier Caruso for helpful correspondence regarding Appendix~\ref{sec:appendix}, and John
Enns for helpful remarks on an earlier version of this paper.
Finally we are thankful to the referee for useful comments that in particular helped us to improve the exposition.

\subsection{Notation}
Let $\overline{\Q}$ be an algebraic closure of $\Q$. All number fields $F/\Q$ will be considered as subfields of $\overline{\Q}$ and we write $G_F\defeq \Gal(\overline{\Q}/F)$ to denote the absolute Galois group of $F$. For any rational prime $\ell\in \Q$, we fix algebraic closures $\overline{\Q}_{\ell}$ of $\Q_{\ell}$ and embeddings $\overline{\Q}\into\overline{\Q}_{\ell}$ (hence inclusions $G_{\Q_{\ell}}\into G_{\Q}$). 
The residue field of $\overline{\Q}_{\ell}$, which is an algebraic closure of $\F_{\ell}$, is denoted by $\overline{\F}_{\ell}$.
As above, all algebraic extensions of $\Q_{\ell}$, $\F_{\ell}$ will be considered as subfields of the fixed algebraic closures $\overline{\Q}_{\ell}$, $\overline{\F}_{\ell}$.

Let $k/\F_p$ be a finite extension of degree $f\geq 1$, and let $K_0\defeq W(k)[\frac{1}{p}]$ be the unramified extension of $\Qp$ of degree $f$. 
Suppose that $e \ge 1$ is any divisor of $p^f-1$. (Starting in Section~\ref{sec:filtrations} we will assume $e = p^f-1$. However, in the appendix
it will be convenient to allow more general $e$, in particular $e = 1$.)
We consider the Eisenstein polynomial $E(u)\defeq u^e+p\in \Qp[u]$ and fix a root $\varpi=\sqrt[e]{-p}\in \Qpbar$.
Let $K\defeq K_0(\varpi)$, a tamely totally ramified cyclic extension of $K_0$ of degree $e$ with uniformizer $\varpi$. 

Let $E$ be a finite extension of $\Qp$. We write $\cO_E$ for its ring of integers, $\F$ for its residue field and $\varpi_E\in\cO_E$ to denote a uniformizer. 
We always assume that $K\subseteq E$.

The choice of $\varpi\in K$ provides us with a homomorphism
\begin{align*}
\teich{\omega}_{\varpi}:\Gal(K/\Qp)&\longrightarrow W(k)\s\\
g&\longmapsto\frac{g(\varpi)}{\varpi}
\end{align*}
whose reduction mod $p$ will be denoted by $\omega_{\varpi}$. Note that the inclusion $k\subseteq \F$ induced by $K\subseteq E$ provides us with a niveau $f$ fundamental character $\omega_f:\Gal(K/K_0)\rightarrow \F\s$, namely $\omega_f= \omega_\varpi|_{\Gal(K/K_0)}$.

We denote by $\omega: G_{\Qp}\rightarrow \Fp\s$ the mod $p$ cyclotomic character, so $\omega=\omega_1$.

Write $\varphi$ for the $p$-power Frobenius on $k$. 
We recall the standard idempotent elements $\epsilon_{\sigma}\in k\otimes_{\Fp}\F$ defined for 
$\sigma\in\Hom(k,\F)$, which verify $(\varphi\otimes 1)(\epsilon_{\sigma})=\epsilon_{\sigma\circ\varphi^{-1}}$ and $(\lambda\otimes1)\epsilon_{\sigma}=(1\otimes\sigma(\lambda))\epsilon_{\sigma}$.
We write $\widehat{\epsilon}_{\sigma}\in W(k)\otimes_{\zp}\oe$ for the standard idempotent elements; they reduce to $\epsilon_{\sigma}$ modulo $\varpi_E$.

Our convention on Hodge--Tate weights is that the cyclotomic character $\varepsilon : G_{\Qp} \to \Qp\s$ has Hodge--Tate weight $-1$.

Given a potentially semistable $p$-adic representation $\rho: G_{\Qp}\rightarrow \GL_n(E)$, we write $\WD(\rho)$ to denote the associated Weil--Deligne representation as defined in  \cite{CDT}, Appendix B.1.
Therefore, $\rho\mapsto\WD(\rho)$ is a \emph{covariant} functor.
We refer to $\WD(\rho)|_{I_{\Qp}}$ as the \emph{inertial type} associated to $\rho$.

\section{The local Galois side}
\label{sec:local-galois-side}

In this section we analyze the local Galois side. In particular, we establish Theorem~\ref{intro:main theorem Galois} of the introduction.

\subsection{Fontaine--Laffaille invariant}
\label{FLp}

Let $\rhobar: \GQp\rightarrow \GL_3(\F)$ be a continuous Galois representation. We assume that $\rhobar$ is \emph{maximally non-split} meaning that $\rhobar$ is uniserial and the graded pieces in its socle filtration are at most one dimensional over $\F$.
(Recall that a finite length module is uniserial if it has a unique composition series.) In other words,
\begin{equation}
\label{Galoisrep}
\rhobar\sim \maq{\omega^{a_2+1}\un{\mu_2}}{\ast_1}{\ast}{}{\omega^{a_1+1}\un{\mu_1}}{\ast_2}{}{}{\omega^{a_0+1}\un{\mu_0}}
\end{equation}
for some $a_i\in\Z$, $\mu_i\in\F\s$ and where $\ast_1,\,\ast_2$ are non-split.
Here, $\un{\mu}: G_{\Qp}\rightarrow \F\s$ denotes the unramified character taking the value $\mu$ on a geometric Frobenius element of $G_{\Qp}$.

\subsubsection{\textbf{Preliminaries on Fontaine--Laffaille theory}}
\label{sec:prel:FL}

We briefly recall the theory of Fontaine--Laffaille modules with $\F$-coefficients and its relation with mod-$p$ Galois representations.

A \emph{Fontaine--Laffaille module} $(M,\Fil^{\bullet}M,\phi_{\bullet})$ over $k\otimes_{\Fp}\F$ is the datum of 
\begin{enumerate}
	\item a finite free $k\otimes_{\Fp}\F$-module $M$;
	\item a decreasing filtration $\{\Fil^{j}M\}_{j\in\Z}$ on $M$ by 
$k\otimes_{\Fp}\F$-submodules such that $\Fil^0M=M$ and $\Fil^{p-1}M=0$;
	\item a $\varphi$-semilinear isomorphism $\phi_{\bullet}: \mathrm{gr}^{\bullet}M\rightarrow M$.
\end{enumerate}

Defining the morphisms in the obvious way, we obtain the abelian category $\F\text{-}\mathcal{FL}^{[0,p-2]}$ of Fontaine--Laffaille modules over $k\otimes_{\Fp}\F$. 
Given a Fontaine--Laffaille module $M$ and $\sigma\in\Hom(k,\F)$, we define the Hodge--Tate weights of $M$ with respect to $\sigma$:
\begin{equation*}
\mathrm{HT}_{\sigma}(M)\defeq\left\{i\in\N : \dim_{\F}\left(\frac{\epsilon_{\sigma}\Fil^iM}{\epsilon_{\sigma}\Fil^{i+1}M}\right)\neq 0\right\}.
\end{equation*}

In the remainder of this paper we focus on Fontaine--Laffaille modules with \emph{parallel} Hodge--Tate weights, i.e.\ we assume that for all $i\in\N$, the submodules $\Fil^iM$ are free over $k\otimes_{\fp}\F$. 

\begin{df}
Let $M$ be a Fontaine--Laffaille module with parallel Hodge--Tate weights. A $k\otimes_{\fp}\F$-basis $\underline{e}=(e_0,\dots,e_{n-1})$ on $M$ is \emph{compatible with the Hodge filtration} if for all $i\in\N$ there exists $j_i\in\N$ such that
$\Fil^iM=\sum_{j=j_i}^{n}(k\otimes_{\fp}\F)\cdot e_j$. 
\end{df}
Note that if the graded pieces of the Hodge filtration have rank at most one, then any two compatible bases on $M$ are related by a lower triangular matrix in $\GL_n(k\otimes_{\fp}\F)$.

Given a Fontaine--Laffaille module and a compatible basis $\underline{e}$, it is convenient to describe the Frobenius action via a matrix $\mathrm{Mat}_{\underline{e}}(\phi_{\bullet})\in\GL_n(k\otimes_{\fp}\F)$, defined in the obvious way using the principal symbols $(\mathrm{gr}(e_0),\dots,\mathrm{gr}(e_{n-1}))$ as a basis of $\mathrm{gr}^{\bullet}M$.

\begin{thm}
\label{MainFL}
There is an exact, fully faithful contravariant functor
$$
\mathrm{T}^*_{\mathrm{cris}}:\,\,\F\text{-}\mathcal{FL}^{[0,p-2]}\rightarrow\mathrm{Rep}_{\F}(G_{K_0})
$$
which is moreover compatible with base change: if $K_0'/K_0$ is finite unramified, with residue field $k'/k$, then 
$$
\mathrm{T}^*_{\mathrm{cris}}(M\otimes_{k} k')\cong \mathrm{T}^*_{\mathrm{cris}}(M)\vert_{G_{K_0'}}.
$$
Also, the essential image of $\mathrm{T}^*_{\mathrm{cris}}$ is closed under subquotients.
\end{thm}
\begin{proof}
The statement with $\fp$-coefficients is in \cite{Fontaine-Laffaille}, Th\'eor\`eme 6.1; its analogue with $\F$-coefficient is a formal argument which is left to the reader (cf.\ also \cite{Gao-Liu}, Theorem 2.2.1).
\end{proof}

\begin{lm}
\label{rhobarFL}
Let $\rhobar:G_{\qp}\rightarrow \GL_3(\F)$ be as in \(\ref{Galoisrep}\).
If the integers $a_i$ verify $a_1-a_0>1,\,\, a_2-a_1>1$ and $p-2> a_2-a_0$ then $\rhobar\otimes\omega^{-a_0}$ is Fontaine--Laffaille, i.e.\ it is in the essential image of $\Tcris^*$. 
\end{lm}
\begin{proof}
This follows, for example, from \cite{gee-geraghty}, Lemma 3.1.5.
\end{proof}

In order to obtain the main results on Serre weights (\S \ref{sec:Weight elimination}) and local-global compatibility (\S \ref{sec:local/global compatibility}), we must assume a stronger genericity condition on the integers $a_i$.

\begin{df}
\label{genericitycondition}
We say that a maximally non-split Galois representation $\rhobar:G_{\qp}\rightarrow \GL_3(\F)$ as in (\ref{Galoisrep}) is \emph{generic} if the triple $(a_2,a_1,a_0)$ satisfies the condition
\begin{equation}
\label{genericityG}
 a_1-a_0>2,\,\, a_2-a_1>2,\,\, p-3> a_2-a_0.
\end{equation}
\end{df}

\subsubsection{\textbf{The Fontaine--Laffaille invariant}}

Let $\rhobar: G_{\Qp}\rightarrow \GL_3(\F)$ be as in Definition \ref{genericitycondition}. 
By Lemma \ref{rhobarFL} there is a Fontaine--Laffaille module $M$ with Hodge--Tate weights $\{1,a_1-a_0+1,a_2-a_0+1\}$ such that $\Tcris^*(M)\cong\rhobar\otimes\omega^{-a_0}$ and which is moreover endowed with a filtration by Fontaine--Laffaille submodules $0\subsetneq M_0\subsetneq M_1\subsetneq M_2=M$ induced via $\Tcris^*$ from the cosocle filtration on $\rhobar$ (cf.\ Theorem \ref{MainFL}).

\begin{lm}
\label{leminvFL}
Let 
$M\in \F\text{-}\mathcal{FL}^{[0,p-2]}$ be such that $\Tcris^*(M)\cong\rhobar\otimes{\omega^{-a_0}}$.
There exists a basis $\underline{e}=(e_0,e_1,e_2)$ of $M$ such that 
$$
M_i\cap \Fil^{a_i-a_0+1}M= \F\cdot e_i
$$
for all $i\in\{0,1,2\}$.
\end{lm}
\begin{proof}
This follows by noting that $M_i\cap \Fil^{a_i-a_0+1}M=\Fil^{a_i-a_0+1}M_i$ and that $M_i$ has Hodge--Tate weights ${a_j-a_0+1}$ with $0\leq j\leq i$.
\end{proof}

Note that a basis $\und e$ as in Lemma \ref{leminvFL} is compatible with the Hodge filtration.
From Lemma \ref{leminvFL} we deduce a useful observation on the Frobenius action on $M$.
\begin{corollary}
\label{corinvFL}
Let $\rhobar$, $M$, $\underline{e}$ be as in Lemma \ref{leminvFL}. 
Then
$$
\mathrm{Mat}_{\underline{e}}(\phi_{\bullet})=\maq{\mu_0}{\alpha_{01}}{\alpha_{02}}
{}{\mu_1}{\alpha_{12}}{}{}{\mu_2}\in\GL_3(\F),
$$
where moreover $\alpha_{01},\,\alpha_{12}\in \F\s$.
\end{corollary}
\begin{proof}
This follows from the lemma, by recalling that the Fontaine--Laffaille module associated to a $G_{\Qp}$-character $\omega^r\un{\mu}$ has Hodge--Tate weight $r$ and $\phi=\mu$.
Note that $\alpha_{01},\,\alpha_{12}\neq 0$ as $\rhobar$ is maximally non-split.
\end{proof}

Conversely, we also note that any such matrix defines a Fontaine--Laffaille module whose associated Galois representation is maximally non-split as in \ref{Galoisrep}.

The Fontaine--Laffaille invariant $\FL(\rhobar)$ associated to $\rhobar$ is defined in terms of $\mathrm{Mat}_{\underline{e}}(\phi_{\bullet})$. 
\begin{lm}
\label{leminvFL1}
Keep the hypotheses of Corollary \ref{corinvFL}.
The element $\frac{\alpha_{02}}{\alpha_{01}\alpha_{12}}$ deduced from $\mathrm{Mat}_{\underline{e}}(\phi_{\bullet})$ does not depend on the choice of the basis $\underline{e}$.
\end{lm}
\begin{proof}
By Lemma \ref{leminvFL} any other such basis is of the form $\beta_i e_i$ for $\beta_i\in\F\s$. The lemma follows.
\end{proof}

\begin{df}
\label{definvFL}
Let $\rhobar:\GQp\rightarrow \GL_3(\F)$ be maximally non-split and generic as in Definition \ref{genericitycondition}.  
Let $M$ be the Fontaine--Laffaille module associated to $\overline{\rho}\otimes{\omega^{-a_0}}$, $\underline{e}=(e_0,e_1,e_2)$ a basis of $M$ as in Lemma \ref{leminvFL} and let
$$
\mathrm{Mat}_{\underline{e}}(\phi_{\bullet})=\maq{\mu_0}{\alpha_{01}}{\alpha_{02}}{}{\mu_1}{\alpha_{12}}{}{}{\mu_2}
$$
be the matrix of the Frobenius action on $M$.

The Fontaine--Laffaille invariant associated to $\rhobar$ is defined as
$$
\FL(\rhobar)\defeq\frac{\det\begin{pmatrix}\alpha_{01}&\alpha_{02}\\\mu_1&\alpha_{12}\end{pmatrix}}{-\alpha_{02}}
\in \mathbb{P}^1(\F)=\F\cup\{\infty\}.
$$
By Lemma \ref{leminvFL1} it is well-defined.
\end{df}

\begin{rk}
Let $\rhobar$ be maximally non-split as in (\ref{Galoisrep}). The isomorphism class of $\rhobar$ is determined by the diagonal characters $\omega^{a_i+1}\un{\mu_i}$ and the Fontaine--Laffaille invariant $\FL(\rhobar)$.
Note that $\FL(\rhobar)$ can take any value in $\mathbb{P}^1(\F)$ except for $\mu_1$.
(Similarly, a maximally non-split Galois representation $\rhobar:\GQp\rightarrow \GL_n(\F)$ is determined by the diagonal characters and $\binom{n-1}2$ invariants.)
\end{rk}

\begin{rk}
Note that in the situation of Definition \ref{definvFL}, if $\F'/\F$ is a finite field extension, then $\FL(\rhobar\otimes_\F \F')=\FL(\rhobar)$.
\end{rk}

\begin{rk}
\label{rmk FL duality}
We leave it to the reader to show that the Fontaine--Laffaille module associated to $\rhobar^{\vee}\otimes\omega^{a_2+2}$ is described by 
$$
\begin{pmatrix}&&1\\&1&\\1&&\end{pmatrix}\cdot {}^t\mat_{\und{e}}(\phi_\bullet)^{-1}\cdot\begin{pmatrix}&&1\\&1&\\1&&\end{pmatrix}.
$$
As a consequence, $\FL(\rhobar^{\vee})=\FL(\rhobar)^{-1}$.
\end{rk}

\subsection{$p$-adic Hodge theory}
\label{sec:p-adic-hodge}

This section mainly consists of a review of some integral $p$-adic Hodge theory, although many of the results are not available 
in the literature in the form or generality that we need. 

In the first subsection (\S \ref{subsubBrModdd}) we define the categories of mod-$p$ objects we are going to work with (Breuil modules with descent data, \'etale $\varphi$-modules, etc.). Moreover we obtain a key result, Corollary \ref{corollary comparison1}, which provides a criterion for deciding when a given Breuil module with descent data and a Fontaine--Laffaille module have isomorphic Galois representations.

The second subsection (\S \ref{Linear awdd}) is of a more technical nature. On the one hand we make two of the functors from \S \ref{subsubBrModdd} (relating Breuil, Fontaine--Laffaille and \'etale $\varphi$-modules) more explicit. We also provide a useful change-of-basis result for a Breuil module with descent data.

All missing proofs of this section are contained in \S\ref{appendix: results}.
Our rationale is to state in this section all the results we need to prove our main results on the Galois side, and to relegate
technical details to the appendix.

\subsubsection{\textbf{Breuil modules with tame descent data}}
\label{subsubBrModdd}

Let $K'\subseteq K_0$ be a subfield containing $\Qp$.
The {residual Breuil ring} $\barS\defeq (k\otimes_{\Fp}\F)[u]/(u^{ep})$ is equipped with an action of $\mathrm{Gal}(K/K')$ by $k\otimes_{\Fp}\F$-semilinear automorphisms. Explicitly if 
 $g\in \mathrm{Gal}(K/K')$ and $a\in k\otimes_{\Fp}\F$, we have
$$
\widehat{g}(au^i)\defeq (g\otimes 1)(a)\cdot(\omega_{\varpi}(g)\otimes 1)^iu^i.
$$
If $\overline{\chi}:\mathrm{Gal}(K/K')\rightarrow \F\s$ is a character, we write $\barS_{\overline{\chi}}$ to denote the $\overline{\chi}$-isotypical component of $\barS$ for the action of $\Gal(K/K')$.

We recall that $\barS$ is equipped with an $k\otimes_{\Fp}\F$-linear derivation $N\defeq -u\frac{\partial}{\partial u}$ and with a semilinear Frobenius $\varphi$ defined by $u\mapsto u^p$ (semilinear with respect to the arithmetic Frobenius $\varphi\otimes 1$ on $k\otimes_{\Fp}\F$), which moreover commute with the action of $\Gal(K/K')$ on $\barS$.

Fix $r\in\{0,\dots,p-2\}$. A \emph{Breuil module of weight $r$ with descent data from $K$ to $K'$} is a quintuple $(\cM, \Fil^r\cM,\varphi_r, N, \{\widehat{g}\})$, where
\begin{enumerate}
	\item $\mathcal{M}$ is a finite free $\barS$-module;
	\item $\Fil^r\cM$ is an $\barS$-submodule of $\M$, verifying $u^{er}\cM\subseteq \Fil^r\cM$;
	\item a morphism $\varphi_{r}:\Fil^r\cM\to\cM$ which is $\varphi$-semilinear and whose image generates $\cM$;
	\item the operator $N:\M\to \M$ is $k\otimes_{\Fp}\F$-linear and satisfies certain axioms (see the beginning of Section 3.2 in \cite{EGH});
	\item an action of $\mathrm{Gal}(K/K')$  on $\cM$ by automorphisms $\widehat{g}$ which  are semilinear with respect to the Galois action on $\barS$ and which preserve $\Fil^r\cM$ and commute with $\varphi_r$ and $N$.\end{enumerate}	
A morphism of Breuil modules is an $\barS$-linear morphism which is compatible, in the evident sense, with the additional structures.

We write $\FBrModdd$ to denote the category of Breuil modules of weight $r$ with descent data from $K$ to $K'$; the field $K'$ will always be clear from the context (and will be specified in case of ambiguities). As we did for the coefficient ring $\barS$, given a character $\overline{\chi}:\mathrm{Gal}(K/K')\rightarrow \F\s$ we write $\cM_{\overline{\chi}}$, $(\Fil^r\cM)_{\overline{\chi}}$ to denote the $\overline{\chi}$-isotypical component of $\cM$, $\Fil^r\cM$ respectively.

We remark that the category $\FBrModdd$ is additive and admits kernels and cokernels (cf.\ \cite{caruso-fourier}, Th\'eor\`eme 4.2.4 and the Remarque following it). In particular a complex 
$$
0\rightarrow \cM_1\stackrel{f_1}{\rightarrow} \cM_2\stackrel{f_2}{\rightarrow} \cM_3\rightarrow 0
$$
in $\FBrModdd$ is \emph{exact} if the morphisms $f_i$ induce exact sequences on the underlying $\barS$-modules $\cM_j$ \emph{and} $\Fil^r\cM_j$  ($j\in\{0,1,2\}$). This endows $\FBrModdd$ with the structure of an exact category (see Proposition \ref{Proposition4Note Florian} below).

We recall that we have an exact, faithful, contravariant functor
\begin{align*}
\Tst^*: \FBrModdd &\rightarrow\mathrm{Rep}_{\F}(G_{K'})\\
\cM&\mapsto\Tst^*(\cM)\defeq \mathrm{Hom}(\cM,\widehat{A}),
\end{align*}
where $\widehat{A}$ is a certain period ring and homomorphisms respect all structures (cf.\ \S \ref{Appendix: Galois}).
We have $\dim_{\F}\,\Tst^*(\cM)=\rank_{\barS}\cM$ (cf.\ \cite{caruso-fourier}, Th\'eor\`eme 4.2.4 and the Remarque following it; see also \cite{EGH}, Lemma 3.2.2).

We will be mainly concerned with the following covariant version of $\Tst^*$: $\Tst^{r}(\cM)\defeq \left(\Tst^*(\cM)\right)^{\vee}\otimes \omega^r$. 
We remark that this is compatible with the notion of duality $\cM\mapsto\cM^*$ on $\FBrModdd$ recalled in \S\ref{subsection:dualities}, namely $\Tst^{r}(\cM)\cong \Tst^*(\cM^*)$, cf.\ the discussion before Corollary 3.2.9 in \cite{EGH}.

We now move to the categories of \'etale $\varphi$-modules. 
In the remainder of this subsection we will assume $K' = K_0$, i.e.\ we will only consider descent data from $K$ to $K_0$. 

We fix the field of norms $k((\underline{\varpi}))$ associated to a suitable Kummer extension $K_{\infty}$ of $K$ (see \S \ref{appendix: Modules} for its precise definition).
It is endowed with a $p$-power Frobenius endomorphism and with an action of $\Gal(K/K_0)$.

An \emph{\'etale $(\phi,k\otimes_{\Fp}\F((\underline{\varpi})))$-module with descent data} is the datum of a finite free $k\otimes_{\Fp}\F((\underline{\varpi}))$-module $\mathfrak{M}$ endowed with a semilinear injective Frobenius endomorphism $\phi:\mathfrak{M}\rightarrow\mathfrak{M}$ and a semilinear action of $\Gal(K/K_0)$, commuting with $\phi$. 
We write $\F\text{-}\mathfrak{Mod}_{k((\underline{\varpi})),\,\mathrm{dd}}$ to denote the category of \'etale $(\phi,k\otimes_{\Fp}\F((\underline{\varpi})))$-modules with descent data. In the special case when $e = 1$ this category is denoted by 
$\F\text{-}\mathfrak{Mod}_{k((\und{p}))}$.
(In this case $k((\und{p}))$ is also the field of norms associated to a suitable Kummer extension $(K_0)_{\infty}/K_0$.)
We remark that $k((\und \varpi))/k((\und p))$ is a cyclic extension of degree $e$ with Galois group $\Gal(K/K_0)$,
and we write $k((\underline{p}))^s$ to denote a choice of separable closure of $k((\underline{p}))$ containing $k((\und \varpi))$. For details, see \S\ref{appendix: Modules}.

Finally, recall the category of Fontaine--Laffaille modules $\F\text{-}\mathcal{FL}^{[0,p-2]}$ over $k\otimes_{\fp}\F$ (defined in \S \ref{sec:prel:FL}) and  the category $\FBrModdd[r]$ of Breuil modules of weight $r$ with descent data from $K$ to $K_0$  (defined in \S \ref{subsubBrModdd}).

The relations between the categories introduced so far are summarized in the following proposition.
Its proof, as well as the definition of the functors $M_{k((\und{\varpi}))}$, $\mathcal{F}$, \dots\ can be found in Appendix~\ref{sec:appendix}.
\begin{prop}
\label{proposition comparison}
Let $0\leq r\leq p-2$. We have the following commutative diagram:
\begin{equation*}
\xymatrix@=5pc{
\FBrModdd[r]\ar^-{M_{k((\underline{\varpi}))}}[rr]\ar^-{\Tst^*}[d]&&
\F\text{-}\mathfrak{Mod}_{k((\underline{\varpi})),\,\mathrm{dd}}
\ar_-{\hspace{-1cm}\mathrm{Hom}(-,k((\underline{p}))^s)}^-{\cong}[dl]\\
\mathrm{Rep}_{\F}(G_{K_0})\ar^{\mathrm{res}}[r]&\mathrm{Rep}_{\F}(G_{(K_0)_{\infty}})&\\
\F\text{-}\mathcal{FL}^{[0,p-2]}\ar_{\mathrm{T}^*_{\mathrm{cris}}}[u]\ar^-{\cF}[rr]&&
\F\text{-}\mathfrak{Mod}_{k((\underline{p}))}
\ar^-{\hspace{-1cm}\mathrm{Hom}(-,k((\underline{p}))^s)}_-{\cong}[ul]
\ar_-{-\otimes_{k((\underline{p}))}k((\underline{\varpi}))}^-\cong[uu]
}
\end{equation*}
where the descent data is from $K$ to $K_0$. Moreover, the functor $\mathrm{res}\circ \Tcris^*$ is fully faithful.
\end{prop}

We record the following immediate, yet crucial, corollary.
\begin{corollary}
\label{corollary comparison1}
Let $0 \le r\leq p-2$ and let $\cM$, $M$ be objects in $\FBrModdd$ and $\F\text{-}\mathcal{FL}^{[0,p-2]}$ respectively.  Assume that $\Tst^*(\cM)$ is Fontaine--Laffaille.
If 
\begin{equation*}
M_{k((\underline{\varpi}))}(\cM)\cong {\mathcal{F}}(M)\otimes_{k((\underline{p}))}k((\underline{\varpi}))
\end{equation*}
then one has an isomorphism of $G_{K_0}$-representations
$$
\Tst^*(\cM)\cong \Tcris^*(M).
$$
\end{corollary}

Let us explain the role that Corollary \ref{corollary comparison1} plays in the proof of our main theorem on the Galois side
(Theorem~\ref{main theorem Galois}). Thus suppose that $\rho$ is potentially semistable of Hodge--Tate weights $\{-2,-1,0\}$
with reduction $\rhobar$ that is maximally non-split and generic as in Definition~\ref{genericitycondition}.  Associated to
$\rho$ is a strongly divisible module $\wh\cM$, whose reduction $\cM\in\FBrModdd$ has Galois representation $\rhobar$.
Corollary \ref{corollary comparison1} will allow us to compute the Fontaine--Laffaille module $M$ associated to $\rhobar$ --
and hence the invariant $\FL(\rhobar)$ -- in terms of $\wh\cM$, so in terms of $\rho$.  
In fact, we will start with $\cM$ in the top left corner of the diagram and then go around it in a clockwise sense:
first computing $M_{k((\underline{\varpi}))}(\cM)$, then descending it to $\F\text{-}\mathfrak{Mod}_{k((\underline{p}))}$
and finally to $\F\text{-}\mathcal{FL}^{[0,p-2]}$. 
To obtain the precise result,
Theorem~\ref{main theorem Galois}, we will moreover need more information about $\wh\cM$, and this will be obtained in
\S\ref{sec:filtrations}.

\subsubsection{\textbf{Linear algebra with descent data}}
\label{Linear awdd}

We continue to assume that $K' = K_0$.  It will be convenient to introduce bases $\und e$ (resp.\ generating sets $\und f$)
of a Breuil module $\cM$ with descent data (resp.\ of $\Fil^r \cM$) that are compatible with the action of $\Gal(K/K_0)$, and
to describe $\Fil^r \cM$ (resp.\ $\varphi_r$) by matrices with respect to $\und e$, $\und f$.  We then use this formalism to
make the functors $M_{k((\und{\varpi}))}$, $\mathcal{F}$ appearing in the diagram of Proposition \ref{proposition comparison}
more explicit, as well as we obtain a change-of-basis result for a Breuil module with descent data. Proofs can be found in \S\ref{appendix: results}.

\begin{df}
We say that a Breuil module $\cM\in \FBrModdd$ \emph{is of type} $\omega_{\varpi}^{a_0}\oplus\dots\oplus\omega_{\varpi}^{a_{n-1}}$ (where $a_i\in \Z$) if $\cM/u\cM\cong 
\oplus_{i=0}^{n-1}(\omega_{\varpi}^{a_i}\otimes 1)$ as $(k\otimes_{\Fp} \F)[\Gal(K/K_0)]$-modules.
Equivalently, $\cM$ has an $\barS$-basis $(e_0,\dots,e_{n-1})$ such that $\widehat{g}e_i=(\omega_{\varpi}^{a_i}(g)\otimes 1) e_i$ for all $i$ and all $g\in \Gal(K/K_0)$. We call such a basis a \emph{framed basis of $\cM$}.

If $\cM$ is of type $\omega_{\varpi}^{a_0}\oplus\dots\oplus\omega_{\varpi}^{a_{n-1}}$ we say that $(f_0,\dots,f_{n-1})$ is a \emph{framed system of generators of $\Fil^r\cM$} if $\Fil^r\cM=\sum_{i=0}^{n-1}\barS\cdot f_i$ and $\widehat{g}f_i=(\omega_{\varpi}^{p^{-1}a_i}(g)\otimes 1) f_i$ for all $i$ and all $g\in \Gal(K/K_0)$.
\end{df}

To justify the claim implicit in this definition, choose an $\barS/u$-basis $(\o e_1,\dots,\o e_{n-1})$ of $\cM/u\cM$ such that
$\widehat{g}\cdot\o e_i=(\omega_{\varpi}^{a_i}(g)\otimes 1) \o e_i$ for all $i$ and all $g\in \Gal(K/K_0)$. Since
$(k\otimes_{\Fp} \F)[\Gal(K/K_0)]$ is a semisimple (commutative) ring we can pick a $(k\otimes_{\Fp} \F)[\Gal(K/K_0)]$-linear splitting
of $\cM \onto \cM/u\cM$ and hence find $e_i \in \cM$ lifting $\o e_i$ such that $\widehat{g}e_i=(\omega_{\varpi}^{a_i}(g)\otimes 1) e_i$.
By Nakayama's lemma, the $e_i$ form an $\barS$-basis of $\cM$.

The notion of a framed basis (resp.\ a framed system of generators) depends on an ordering of the integers $a_i$. It will always be clear from the context which ordering we use.

\begin{lm}
\label{lemma generators dd}
If $\cM\in \FBrModdd$ is of type $\oplus_{i=0}^{n-1}\omega_{\varpi}^{a_i}$, then $\Fil^r\cM$ admits a framed system of generators. 
\end{lm}

\begin{df}
\label{definition linear algebra dd}
Let $\cM\in\FBrModdd$ be of type $\omega_{\varpi}^{a_0}\oplus\dots\oplus\omega_{\varpi}^{a_{n-1}}$. 
Let $\underline{e}$, $\underline{f}$
be a framed basis and a framed system of generators of $\cM$, $\Fil^r\cM$ respectively.
The \emph{matrix of the filtration}, with respect to $\underline{e},\underline{f}$, is the element $\mat_{\underline{e},\underline{f}}\left(\Fil^r\cM\right)\in \mathrm{M}_{n}(\barS)$ verifying
$$
\underline{f}=\underline{e}\cdot\mat_{\underline{e},\underline{f}}\left(\Fil^r\cM\right).
$$
We define the \emph{matrix of Frobenius}  with respect to $\underline{e}$, $\underline{f}$ as the element $\mathrm{Mat}_{\underline{e},\underline{f}}(\varphi_r)\in \mathrm{GL}_{n}(\barS)$ characterized by
$$
\varphi_r(\und{f})=\underline{e}\cdot \mathrm{Mat}_{\underline{e},\underline{f}}(\varphi_r).
$$
\end{df}

As we require $\underline{e},\,\underline{f}$ to be framed, the coefficients $\mat_{\underline{e},\underline{f}}\left(\Fil^r\cM\right)_{ij}$ verify the following conditions:
$$
\mat_{\underline{e},\underline{f}}\left(\Fil^r\cM\right)_{ij}\in \barS_{\omega_{\varpi}^{p^{-1}a_j-a_i}}.
$$
Concretely $\mat_{\underline{e},\underline{f}}\left(\Fil^r\cM\right)_{ij}=u^{[p^{-1}a_j-a_i]}s_{ij}$, where $[x]\in\{0,\dots,e-1\}$ is defined by $[x]\equiv x\,\,\mathrm{mod}\,e$ for $x\in \Z/e\Z$ and $s_{ij}\in \barS_{\omega_{\varpi}^{0}}=k\otimes_{\fp}\F[u^e]/(u^{ep})$.

On the other hand, $\mathrm{Mat}_{\underline{e},\underline{f}}(\varphi_r)\in \mathrm{GL}_{n}^{\Box}(\barS)$, where
$$
\GL_{n}^{\Box}(\barS)\defeq
\left\{A\in \GL_{n}(\barS) : A_{ij}\in \barS_{\omega_{\varpi}^{a_j-a_i}}\,\, \text{for\,all}\,\,0\leq i,j\leq n-1 \right\}.
$$

\begin{lm}
\label{lemma lawdd 1}
Let $\cM$ be a Breuil module of type $\omega_{\varpi}^{a_0}\oplus\dots\oplus\omega_{\varpi}^{a_{n-1}}$, and let $\underline{e}$ be a framed basis of $\cM$ and  $\underline{f}$ a framed system of generators of $\Fil^r\cM$, respectively.
Let $V\defeq \mat_{\underline{e},\underline{f}}\left(\Fil^r\cM\right)\in \mathrm{M}_{n}(\barS)$ and $A\defeq \mathrm{Mat}_{\underline{e},\underline{f}}(\varphi_r)\in
\mathrm{GL}_{n}^{\Box}(\barS)$.
Then there exists a basis $\und{\mathfrak{e}}$ of $M_{k((\underline{\varpi}))}(\cM^{\ast})$ with $\widehat{g}\cdot\mathfrak{e}_i=(\omega_{\varpi}^{-p^{-1}a_i}(g)\otimes 1) \mathfrak{e}_i$ for all $i$ and $g\in \Gal(K/K_0)$, such that $\mathrm{Mat}_{\underline{\mathfrak{e}}}(\phi)\in \mathrm{M}_{n}(k\otimes_{\Fp}\F[[\underline{\varpi}]])$ is given by any chosen lift of ${}^t{V}\cdot{}^t{A}^{-1}\in \mathrm{M}_n(\barS)$ via the morphism $k\otimes_{\Fp}\F[[\underline{\varpi}]]\onto \barS$ sending $\sum \lambda_i\und{\varpi}^i$ to $\sum \lambda_i u^i$ and such that $\big(\mathrm{Mat}_{\underline{\mathfrak{e}}}(\phi)\big)_{ij}\in (k\otimes_{\Fp} \F[[\und{\varpi}]])_{\omega_{\varpi}^{p^{-1}a_i-a_j}}$ for all $0\leq i,j\leq n-1$.
\end{lm}

\begin{lm}
\label{lemma lawdd 2}
Let $M\in\F\text{-}\mathcal{FL}^{[0,p-2]}$ be a rank $n$ Fontaine--Laffaille module with parallel Hodge--Tate weights $0\leq m_{0}\leq \dots\leq m_{n-1}\leq p-2$ \(counted with multiplicities\).

Let $\underline{e}=(e_{0},\dots,e_{n-1})$ be a $k\otimes_{\fp}\F$-basis of $M_i$ that is compatible with the Hodge filtration $\mathrm{Fil}^{\bullet}M$, and let $F\in \GL_{n}(k\otimes_{\fp}\F)$ be the associated matrix of the Frobenius $\phi_{\bullet}:\gr^{\bullet}M\isoto M$.  

Then there exists a basis $\underline{\mathfrak{e}}$ of $\mathfrak{M}\defeq \cF(M)$ in $\F\text{-}\mathfrak{Mod}_{k((\underline{p}))}$ such that the Frobenius $\phi$ on $\mathfrak{M}$ is described by
$$
\mathrm{Mat}_{\underline{\mathfrak{e}}}(\phi)=\mathrm{Diag}(
\underline{p}^{m_{0}}\dots\underline{p}^{m_{n-1}})F.
$$
\end{lm}

\begin{lm}
\label{lemma base change matrix}
Let $\cM\in\FBrModdd$ be a Breuil module of type $\omega_{\varpi}^{a_0}\oplus\dots\oplus\omega_{\varpi}^{a_{n-1}}$ and let $\underline{e}$, $\underline{f}$ be a framed basis of $\cM$ and a framed system of generators of $\Fil^r\cM$ respectively.

Write $V\defeq\mat_{\underline{e},\underline{f}}\left(\Fil^r\cM\right)$, $A\defeq \mathrm{Mat}_{\underline{e},\underline{f}}(\varphi_r)$.
Assume that there exist elements $V'\in \mathrm{M}_{n}(\barS)$, $B\in \GL_n(\barS)$ such that $V'_{ij}\in\barS_{\omega_{\varpi}^{p^{-1}a_j-a_i}}$, $B_{ij}\in \barS_{\omega_{\varpi}^{p^{-1}(a_j-a_i)}}$ and
\begin{equation}
\label{equation base change}
A V'\equiv VB\,\,\mathrm{mod}\,u^{e(r+1)}.
\end{equation}
Then $\underline{e}'\defeq \underline{e}\cdot A$
is a framed basis of $\cM$, $\underline{f}'\defeq \underline{e}'\cdot V'$ is a framed system of generators for $\Fil^r\cM$, and  $\mat_{\und{e}',\und{f}'}(\varphi_r)=\varphi(B)$.
\end{lm}

\subsection{Breuil submodules and Galois representations}
\label{sec:breuil-submod}

In this subsection we discuss
some preliminaries on subobjects and quotients in the category $\FBrModdd[r]$. Even though these notions are presumably well-known to the experts, we did not find a suitable reference in the literature. The main result, Proposition \ref{Proposition7Note Florian}, is a slight improvement of a result of Caruso (\cite{caruso-fourier}).
All missing proofs of this subsection are contained in \S\ref{appendix: results-fil}.

In what follows, we let $\barS_k\defeq k[u]/u^{ep}$.
\begin{df}
\label{definition:Breuil submodule}
Let $\cM$ be an object in $\FBrModdd$. An $\barS$-submodule $\cN\subseteq\cM$ is said to be a \emph{Breuil submodule} if $\cN$ fulfills the following conditions:
\begin{enumerate}
	\item $\cN$ is an $\barS_k$-direct summand of $\cM$;
	\item $N(\cN)\subseteq \cN$ and $\widehat{g}(\cN)\subseteq \cN$ for all $g\in\Gal(K/K')$;
	\item $\varphi_r(\cN\cap \Fil^r\cM)\subseteq \cN$.
\end{enumerate}
\end{df}

\begin{lm}
\label{Lemma2Note Florian}
Let $\cM$ be an object in $\FBrModdd$ and let $\cN\subseteq\cM$  be a Breuil submodule.
Then the $\barS$-modules $\cN$,  $\cM/\cN$ with their induced structures are objects of $\FBrModdd$ and the sequence
$$
0\rightarrow \cN\rightarrow \cM\rightarrow \cM/\cN\rightarrow 0
$$
is exact in $\FBrModdd$.
Conversely, given an exact sequence 
$$
0\rightarrow \cM_1\stackrel{f}{\rightarrow} \cM\rightarrow \cM_2\rightarrow 0
$$ 
in $\FBrModdd$, then $f(\cM_1) \subseteq \cM$ is a Breuil submodule.
\end{lm}

An immediate consequence of Lemma \ref{Lemma2Note Florian} is that the notion of Breuil submodule is transitive:
\begin{lm}
\label{Lemma3Note Florian}
Let $\cM$ be an object in $\FBrModdd$.
\begin{enumerate}
\item If $\cM_1\subseteq \cM$ and $\cM_2\subseteq \cM_1$ are Breuil submodules, then so is $\cM_2\subseteq \cM$.
\item Let $\cM_1$, $\cM_2$ be Breuil submodules of $\cM$, and assume that $\cM_2\subseteq \cM_1$. Then $\cM_2$ is a Breuil submodule of $\cM_1$, and the Breuil module structures on $\cM_2$ inherited from $\cM_1$ and $\cM$ coincide. Similarly, $\cM_1/\cM_2 \subseteq \cM/\cM_2$ is a Breuil submodule, and the Breuil module structures on $\cM_1/\cM_2$, as a  Breuil submodule of $\cM/\cM_2$ and as a quotient of $\cM_1$, coincide.
\end{enumerate}
\end{lm}

\begin{prop}
\label{Proposition4Note Florian}
With the above notion of exact sequence, the category $\FBrModdd$ is an exact category in the sense of \cite{keller} and  $\Tst^{r}$ is an exact functor. 
\end{prop}

We can now state a crucial result relating Breuil submodules of $\cM$ and subrepresentations of $\Tst^{r}(\cM)$.
It improves \cite{caruso-fourier}, Proposition 2.2.5 (cf.\ also \cite{EGH}, Corollary 3.2.9).

\begin{prop}
\label{Proposition7Note Florian}
Let $\cM$ be an object in $\FBrModdd$. 
There is a natural order-preserving bijection
\begin{equation*}
\Theta:\left\{\text{Breuil submodules in}\,\,\cM\right\}\stackrel{\sim}{\longrightarrow}
\left\{\text{$G_{K'}$-subrepresentations of}\,\,\Tst^{r}(\cM)\right\}
\end{equation*}
sending $\cN\subseteq\cM$ to the image of $\Tst^{r}(\cN)\into\Tst^{r}(\cM)$. Moreover, if $\cM_2 \subseteq \cM_1$ are Breuil submodules of $\cM$, then
$\Theta(\cM_1)/\Theta(\cM_2) \cong \Tst^{r}(\cM_1/\cM_2)$.
\end{prop}

\subsection{On filtrations}
\label{sec:filtrations}

Our goal is to give a rough classification of the filtration on certain strongly divisible modules $\wh\cM$ of rank 3 whose associated $p$-adic Galois representation 
$\rho$ has maximally
non-split and generic reduction $\rhobar$. The idea is to start with a simpler analysis in characteristic $p$ (\S \ref{sec:fil-breuil-mod}), using that the mod $p$ reduction $\cM$ of $\wh\cM$ (a Breuil module) has associated
mod $p$ Galois representation $\rhobar$. The resulting Corollary \ref{corollary filtration Breuil module} concerning $\Fil^2 \cM$ helps us in our analysis
of $\Fil^2 \wh\cM$ in \S\ref{subsubsection Filtration str div}. We obtain our classification in Proposition \ref{Shape filtration str div modules}.
This will be a key input into our main local result on the Galois side, Theorem \ref{main theorem Galois}.

Our conventions regarding strongly divisible modules will be explained at the start of \S\ref{subsubsection Filtration str div}.
\emph{For the remainder of Section~\ref{sec:local-galois-side} we will assume that $e = p^f -1$.}

\subsubsection{\textbf{Filtration on Breuil modules}}
\label{sec:fil-breuil-mod}

We now obtain the first structure results for Breuil modules with descent data giving rise to maximally non-split $\rhobar$.
\begin{prop}
\label{mainsection2.2}
Let $K_0=\qp$ and suppose that $\cM\in\FBrModdd[2]$ is of type $\omega^{a_0}\oplus\omega^{a_1}\oplus\omega^{a_2}$. Assume that $\rhobar\defeq \Tst^{2}(\cM)$ is maximally non-split and generic as in Definition \ref{genericitycondition}.

There is a framed basis $\underline{e}=(e_0,e_1,e_2)$ of $\cM$ and a framed system of generators $\underline{f}=(f_0,f_1,f_2)$ for $\Fil^2\cM$ such that the coordinates of the elements in $\underline{f}$ with respect to $\und e$ are described as follows:
\begin{align*}
f_0&=\,\columnvct{u^e\\ \mu u^{e-(a_1-a_0)}\\\nu u^{e-(a_2-a_0)}},& f_1&=\,\columnvct{0\\u^e\\\lambda u^{e-(a_2-a_1)}},& f_2&=\,\columnvct{0\\0\\u^e },
\end{align*}
where $\lambda,\mu,\nu\in\F$ verify moreover $\lambda\mu\neq 0$.
\end{prop}

The proof of Proposition \ref{mainsection2.2} will occupy the remainder of this subsection.
We start with the following lemma. It gives a concrete criterion for the Galois representation associated to a rank two Breuil module to split as a direct sum of characters.

\begin{lm}
\label{lemmasemisimple}
Assume $K_0=\qp$ \(so $e=p-1$\) and let
$$
0\rightarrow \cM_1\rightarrow \cM\rightarrow \cM_2\rightarrow0
$$
be an extension of rank one objects in $\FBrModdd[2]$. 
For each $i\in\{1,2\}$, assume that $\cM_i$ is of type $\omega^{a_i}$, with $\omega^{a_1}\not\cong \omega^{a_2}$ and suppose $\Fil^2\cM_i=u^{\delta_i e}\cM_i$ for $0\leq \delta_1\leq\delta_2\leq 2$. 
Finally, assume that $\rhobar\defeq \Tst^{2}(\cM)$ is Fontaine--Laffaille \(possibly after a twist\).
\begin{enumerate}
	\item If the extension of $\barS$-modules
\begin{equation}
\label{split filtration}
0\rightarrow \Fil^2\cM_1\rightarrow \Fil^2\cM\rightarrow \Fil^2\cM_2\rightarrow0
\end{equation}
splits, then $\rhobar$ splits as a direct sum of two characters.
	\item If $\delta_1=1,\,\delta_2=2$ then $\rhobar$ splits as a direct sum of two characters.
\end{enumerate}
\end{lm}
\begin{proof}
By assumption, $\cM$ is of type $\omega^{a_1}\oplus\omega^{a_2}$.
Let $\underline{e}\defeq (e_1,e_2)$ be a framed basis of $\cM$ such that $\cM_1=\barS\cdot e_1$.
In what follows, we define $[a_2-a_1]\in\{1,\dots,e-1\}$  by the requirement $[a_2-a_1]\equiv a_2-a_1$ modulo $e$.

From the exact sequence (\ref{split filtration}) we can find a framed system of generators of $\Fil^2\cM$ of the form $f_1=u^{\delta_1e}e_1$, $f_2=u^{\delta_2e}e_2+xu^{[a_2-a_1]}e_1$, where $x\in \barS_{\omega^0}$. 
Moreover, we have
$$
A\defeq \mathrm{Mat}_{\underline{e},\underline{f}}(\varphi_2)=\maqdue{\alpha}{\gamma u^{[a_2-a_1]}}{}{\beta}\in\GL_2(\barS)
$$
for some $\alpha,\beta\in\barS^{\times}_{\omega^0}$ and $\gamma\in \barS_{\omega^0}$.

Assume that the sequence (\ref{split filtration}) splits. Then we can fix a splitting $s:\Fil^2\cM_2\rightarrow \Fil^2\cM$ which we can moreover assume to be $\Gal(K/K_0)$-equivariant (by averaging). Let $\overline{e}_2\in \left(\cM_2\right)_{\omega^{a_2}}$ be a generator of $\cM_2$. Then $s(u^{e\delta_2}\overline{e}_2)$ is killed by $u^{e(p-\delta_2)}$ and hence of the form $u^{\delta_2e}e_2$ for some $e_2\in \cM$.  Without loss of generality, $e_2\in \cM_{\omega^{a_2}}$. It is now easy to see that $(e_1,e_2)$ is a framed basis of $\cM$  and that $(f_1,f_2)=(u^{\delta_1e}e_1,u^{\delta_2e}e_2)$ is a framed system of generators of $\Fil^2\cM$. In other words, we can take $x=0$ above.

From the obvious matrix equality
$$
A\maqdue{u^{\delta_1e}}{}{}{u^{\delta_2e}}=\maqdue{u^{\delta_1e}}{}{}{u^{\delta_2e}} \underbrace{\maqdue{\alpha}{\gamma u^{(\delta_2-\delta_1)e+[a_2-a_1]}}{}{\beta}}_{\defeq B}
$$
and Lemma \ref{lemma base change matrix}, we deduce that $\underline{e}'\defeq \underline{e}\cdot A$ is a basis of $\cM$ such that $\cM_1=\barS\cdot e_1'$, that $\underline{f}'\defeq (u^{\delta_1e}e_1',u^{\delta_2e}e_2')$ is a system of generators for $\Fil^2\cM$ and that 
$$
\mathrm{Mat}_{\underline{e}',\underline{f}'}(\varphi_2)=\varphi(B)=\maqdue{\varphi(\alpha)}{\varphi(\gamma)u^{[a_2-a_1]p}}{}{\varphi(\beta)}.
$$
As $[a_2-a_1]>0$ we can further assume that $\varphi(\gamma)=0$, up to re-iterating the procedure above.
Therefore $M_{\Fp((\underline{\varpi}))}(\cM)$ splits into a direct sum of rank one $(\varphi,\F((\und{\varpi})))$-modules (as can be immediately checked by Lemma \ref{lemma lawdd 1}), hence $\rhobar\vert_{G_{(\Qp)_{\infty}}}$ splits as a direct sum of two characters.
As $\rhobar$ is Fontaine--Laffaille, we deduce from Proposition \ref{proposition comparison} that $\rhobar$ splits into a direct sum of two characters, completing the proof of case $\mathrm{(i)}$.

Let us assume that $\delta_1=1,\,\delta_2=2$. We then have $f_2=u^{2e}e_2+xu^{[a_2-a_1]}e_1$ and by adding a multiple of $f_1$ we may assume $x\in \F$. 
If $x\neq0$ it follows that $\Fil^2\cM=\left\langle u^{[a_2-a_1]}e_1, u^{3e-[a_2-a_1]}e_2\right\rangle$, so $u^{2e}e_2\not\in\Fil^2\cM$, contradiction.
Therefore $x=0$ and we are in the situation of $\mathrm{(i)}$.
\end{proof}

\begin{remark}
The analogous statement of Lemma \ref{lemmasemisimple} holds when $K_0\neq \Qp$. This will be described in detail in a particular case in the proof of Proposition \ref{weight elimination FL}, when $[K_0:\Qp]=2$ and developed in further generality in \cite{MP}, Lemma 3.2.
\end{remark}

We now make use of our knowledge of rank one Breuil modules with descent data and their associated Galois representations
to start to understand $\Fil^2\cM$ in Proposition~\ref{mainsection2.2}.

\begin{lm}
\label{lemmadevissage} Assume $K_0=\qp$ and let $\rhobar:\GQp\rightarrow \GL_3(\F)$ be maximally non-split and generic as in Definition \ref{genericitycondition}.
Let $\cM\in\FBrModdd[2]$ be of type
${\omega}^{a_0}\oplus{\omega}^{a_1}\oplus{\omega}^{a_2}$ and such that 
$\Tst^{2}(\cM)\cong \rhobar$ and write $0=\cM_3\subsetneq \cM_2\subsetneq \cM_1\subsetneq\cM_0\defeq \cM$ to denote the Breuil submodule filtration on $\cM$ deduced from the socle filtration on $\rhobar$ \(cf.\ Proposition \ref{Proposition7Note Florian}\). 

Then for each $i\in\{0,1,2\}$ the rank one quotient $\cM_{i}/\cM_{i+1}\in \FBrModdd[2]$  is of type $\omega^{a_i}$ and its filtration is described by
$$
\Fil^2\left(\cM_{i}/\cM_{i+1}\right)=u^e\left(\cM_{i}/\cM_{i+1}\right).
$$
\end{lm}
\begin{proof}

By Proposition \ref{Proposition7Note Florian},  there exists a permutation $\sigma\in S_3$ such that for any $i\in\{0,1,2\}$ the rank one module $\cM_i/\cM_{i+1}$ is of type $\omega^{a_{\sigma(i)}}$.
This implies, by means of \cite{EGH}, Lemma 3.3.2, that there exists $\delta_i\in\{0,1,2\}$ such that
$\Fil^{2}\left(\cM_i/\cM_{i+1}\right)=u^{e\delta_i}(\cM_i/\cM_{i+1})$ and $\Tst^{2}(\cM_i/\cM_{i+1})\vert_{I_{\Qp}}\cong \omega^{a_{\sigma(i)}+\delta_i}$.

On the other hand 
we have $\Tst^{2}(\cM_i/\cM_{i+1})\vert_{I_{\Qp}}\cong \omega^{a_i+1}$ for all $i$, by the definition of the filtration $\{\cM_{i}\}_{i}$.
As $\rhobar$ is generic (cf.\ Definition \ref{genericitycondition}) we conclude that $\sigma=\mathrm{id}$ and that $\delta_i = 1$ for all $i$.
\end{proof}

\begin{proof}[Proof of Proposition \ref{mainsection2.2}]
Let $0\subsetneq\cM_2\subsetneq\cM_1\subsetneq\cM_0\defeq \cM$ be the filtration by Breuil submodules on $\cM$, obtained from the socle filtration on $\rhobar$ via Proposition \ref{Proposition7Note Florian}.
From Lemma \ref{lemmadevissage} we obtain a framed basis $(e_0,e_1,e_2)$ of $\cM$ with $e_i\in \cM_i$ such that
\begin{equation}
\label{filtration subquotient}
\Fil^2\left(\cM_i/\cM_{i+1}\right)=\langle u^e \overline{e}_i\rangle_{\barS}
\end{equation}
for $i\in\{0,1,2\}$ (with the obvious notation for the elements $\overline{e}_i$).

As the descent data acts semisimply on $\cM_1$, from the exact sequence 
$$
0\rightarrow \Fil^2\cM_2\rightarrow\Fil^2\cM_1\rightarrow\Fil^2\left(\cM_1/\cM_2\right)\rightarrow0
$$ 
we see that
\begin{equation*}
\Fil^2\cM_1=\langle u^ee_2,u^ee_1+s e_2\rangle_{\barS}\quad\text{for\,\,some}\quad s\in \barS_{\omega^{a_1-a_2}}=u^{e-(a_2-a_1)}\barS_{\omega^0}.
\end{equation*}
Without loss of generality, we can assume $s=\lambda u^{e-(a_2-a_1)}$ for some $\lambda\in\F$ and, by Lemma \ref{lemmasemisimple} and the non-splitness assumption, we moreover have $\lambda\neq0$.

In a completely similar fashion we obtain 
\begin{equation*}
\Fil^2\cM=\langle u^ee_2,\ u^ee_1+\lambda u^{e-(a_2-a_1)} e_2,\  u^e e_0+s_1 e_1+s_2 e_2\rangle_{\barS},
\end{equation*}
where $s_1\in \barS_{\omega^{a_0-a_1}},\, s_2\in \barS_{\omega^{a_0-a_2}}$.
As above, we can assume without loss of generality that $s_1= \mu u^{e-(a_1-a_0)}$ and $s_2=\nu u^{e-(a_2-a_0)}$ and we furthermore deduce from Lemma \ref{lemmasemisimple} (applied to $\cM/\cM_2$)  that $\nu\neq0$.
\end{proof}

The following immediate corollary of Proposition~\ref{mainsection2.2} will play an important role in describing the filtration on certain strongly divisible $\oe$-modules, see \S\ref{subsubsection Filtration str div}.

\begin{corollary}
\label{corollary filtration Breuil module}
Let $\cM$ be a Breuil module and $\lambda,\nu,\mu\in\F$ as in the statement of Proposition \ref{mainsection2.2}.

The elementary divisors for $\cM/\Fil^2\cM$ as an $\F[u]$-module are described by one of the following possibilities:
\begin{enumerate}
\item if $\nu(\nu-\lambda\mu)\neq 0$, by $(u^{e-(a_2-a_0)},u^e,u^{e+(a_2-a_0)})$;
\item if $\nu-\lambda\mu=0$, by $(u^{e-(a_2-a_0)},u^{e+(a_2-a_1)},u^{e+(a_1-a_0)})$;
\item if $\nu=0$, by $(u^{e-(a_2-a_1)},u^{e-(a_1-a_0)},u^{e+(a_2-a_0)})$.
\end{enumerate}
We also have:
\begin{enumerate}
\setcounter{enumi}{3}
\item $\Fil^2\cM\subseteq u^{e-(a_2-a_0)}\cM$; 
\item $\left(\Fil^2\cM\right)_{\omega^{a_0}}\not\subseteq u^{e}\cM$;
\item $\left(\Fil^2\cM\cap u^e\cM\right)_{\omega^{a_0}}\subseteq u^{2e-(a_2-a_0)}\cM$;
\item $\left(\Fil^2\cM\right)_{\omega^{a_2}}\subseteq u^{e}\cM$.
\end{enumerate}
\end{corollary}

\subsubsection{\textbf{Filtration on strongly divisible modules}}
\label{subsubsection Filtration str div}

In this section we pursue the analysis started in \S\ref{sec:fil-breuil-mod}.
The main result of this section is Proposition \ref{Shape filtration str div modules}.

As in \S \ref{subsubBrModdd}, we let $K'\subseteq K_0$ be a subfield containing $\Qp$.
The ring $S_{W(k)}$ (cf.\ \cite{breuil-griffiths}, \S 4.1) is defined as the $p$-adic completion of $W(k)[u,\frac{u^{ie}}{i!}]_{i\in\N}$.
The ring $S_{W(k)}$ is endowed with a descending filtration $\{\Fil^iS_{W(k)}\}_{i\in\N}$, a semilinear Frobenius $\varphi$, a $W(k)$-linear derivation $N$, and with a Galois action by $W(k)$-algebra endomorphisms defined by $\widehat{g}(u)=\teich{\omega}_{\varpi}(g)u$ for any $g\in\mathrm{Gal}(K/K')$.
In particular, the action of any $g\in \mathrm{Gal}(K/K')$ preserves the filtration and commutes with the Frobenius and the monodromy on $S_{W(k)}$. 
By extension of scalars,  the ring $S\defeq S_{W(k)}\otimes_{\Zp} \oe$ is endowed with the additional structures inherited from $S_{W(k)}$.
Note in particular that we have a natural map $S\onto\barS$, defined as the reduction modulo $(\varpi_E,\Fil^p S)$.
For more details, see \cite{EGH}, \S 3.1.

Fix $r\in\{0,\dots,p-2\}$. We refer to \cite{EGH}, \S 3.1 for the definition of the category $\OEModdd$ of strongly divisible $\oe$-modules of weight $r$ and with descent data from $K$ to $K'$. 
The objects are certain quintuples $(\widehat{\cM}, \Fil^r\widehat{\cM}, \varphi, N, \{\widehat{g}\})$, where
$\widehat{\cM}$ is a finite free $S$-module.
There is a contravariant functor $\Tst^{*,K'}:\OEModdd\rightarrow \mathrm{Rep}_{\oe}(G_{K'})$, which by a theorem of T.\ Liu (\cite{liu-breuil}, Theorem 2.3.5, cf.\ also \cite{EGH}, Proposition 3.1.4) provides an equivalence of categories of $\OEModdd$ with the category of $G_{K'}$-stable $\oe$-lattices in finite dimensional $E$-representations of $G_{K'}$ that become semistable over $K$ and have Hodge--Tate weights in $[-r,0]$.
As in the case of Breuil modules, we consider its covariant version defined by $\Tst^{K',r}(\wh\cM)\defeq \big(\Tst^{*,K'}(\wh\cM)\big)^{\vee}\otimes \varepsilon^r$.

We also recall that if $\widehat{\cM}\in\OEModdd$ the base change $\widehat{\cM}\otimes_S \barS$ is naturally an object of $\FBrModdd$ and one has a natural isomorphism $\Tst^{*, K'}(\widehat{\cM})\otimes_{\oe}\F\cong
\Tst^*(\widehat{\cM}\otimes_S \barS)$ of $\F[G_{K'}]$-modules.

We write $\mathrm{Mod}_E^{\mathrm{w.a.}}(\varphi,N,K/K')$  for the category of weakly admissible filtered $(\varphi,N, K/K', E)$-modules (see e.g.\ \cite{EGH}, \S 3.1).
In particular, the underlying objects are finite free $K_0\otimes_{\Qp}E$-modules.
We recall the contravariant equivalence of categories $\Dst^{*,K'}: \mathrm{Rep}_E^{\text{$K$-st}}(G_{K'})\rightarrow \mathrm{Mod}_E^{\mathrm{w.a.}}(\varphi,N,K/K')$, where $\mathrm{Rep}_E^{\text{$K$-st}}(G_{K'})$ denotes the category of finite dimensional $E$-representations of $G_{K'}$ that become semistable over $K$.
If $\rho\in \mathrm{Rep}_E^{\text{$K$-$\mathrm{cris}$}}(G_{K'})$ has Hodge--Tate weights in $[-r,0]$, we define $\Dst^{K',r}(\rho)\defeq \Dst^{*,K'}\left(\rho^{\vee}\otimes \varepsilon^r\right)$.

As in the mod $p$ setting, given an $\oe[\Gal(K/K')]$-module $X$ and a character $\chi:\Gal(K/K')\rightarrow \oe\s$ we let $X_\chi$ denote its $\chi$-isotypical component.

We first require two lemmas.

\begin{lm} 
\label{lemma str div dd}
Assume that $K'=\Qp$ and that the $p$-adic Galois representation $\Tst^{\Qp,r}(\widehat{\cM})\otimes_{\oe}E$ has inertial type $\oplus_{i=0}^{n-1}\teich{\omega}_f^{a_i}$.
Then $\widehat{\cM}\otimes_S\barS\in \FBrModdd$ is of type $\oplus_{i=0}^{n-1}\omega_{\varpi}^{a_i}$.
\end{lm}
\begin{proof}
The assumption that $K'=\Qp$ implies that the multi-set $\{\teich{\omega}_f^{a_i}\}_{i=0}^{n-1}$ (and hence the multi-set $\{\omega_{\varpi}^{a_i}\}_{i=0}^{n-1}$) is stable under the action of the $p$-power Frobenius. Together with \cite{EGH}, Proposition 3.3.1, this is all we need to use the argument at the beginning of \cite{EGH}, proof of Theorem 3.3.13 to construct the required basis for $\widehat{\cM}\otimes_S\barS$.
\end{proof}

\begin{lm}
\label{lemma shape filtration}
Assume $K_0=\Qp$ and let $\rho:G_{\Qp}\rightarrow \GL_3(E)$ be a Galois representation, becoming semistable over $K$, with Hodge--Tate weights $\{-2,-1,0\}$. Let $\widehat{\cM}\in \OEModdd[2]$ be such that $\Tst^{\Qp,2}(\widehat{\cM})\otimes_{\oe}E\cong \rho$.

Let $X\defeq \left(\frac{\Fil^2\widehat{\cM}}{\Fil^2S\cdot\widehat{\cM}}\right)\otimes_{\oe}E$ and 
$Y\defeq \left(\frac{\Fil^1S\cdot\widehat{\cM}}{\Fil^2S\cdot\widehat{\cM}}\right)\otimes_{\oe}E$.
Then for any character $\chi:\Gal(K/\Qp)\rightarrow E\s$ we have
\begin{enumerate}
	\item $\dim_E X_{\chi}=3$;
	\item	$\dim_E (X\cap Y)_{\chi}=2$.
\end{enumerate}
Moreover, multiplication by $ u\in S$ induces an isomorphism
$X_{\chi}\stackrel{\sim}{\longrightarrow}X_{\chi \widetilde{\omega}}$.
\end{lm}
\begin{proof}
Let $D\defeq \Dst^{*,\Qp}(\rho)$.
As $\rho\vert_{G_K}$ is semistable, with Hodge--Tate weights $\{-2,-1,0\}$, the $E$-linear spaces $\gr^i(\Fil^{\bullet}D_K)$ are at most one dimensional and they are non-zero if and only if $i\in\{0,1,2\}$.

Let $\cD\defeq D\otimes_{E}S_{E}$.
We then have $X=\Fil^2\cD/\left(\Fil^2S_E\cdot\cD\right)$, $Y=\left(\Fil^1S_E\cdot \cD\right)/\left(\Fil^2S_E\cdot \cD\right)$ and by the analogue with $E$-coefficients of \cite{breuil-griffiths}, Proposition A.4, we deduce that
\begin{equation*}
	\Fil^1\cD=(\Fil^1S_E) \widehat{f}_0\oplus S_E\widehat{f}_1\oplus S_E\widehat{f}_2,
	\quad\Fil^2\cD=(\Fil^2S_E) \widehat{f}_0\oplus (\Fil^1S_E)\widehat{f}_1\oplus S_E\widehat{f}_2,
\end{equation*}
for some $S_E$-basis $\widehat{f}_i$ of $\cD$.

From the $E$-linear isomorphism
$$
\frac{S}{\Fil^mS}\cong \bigoplus_{i=0}^{m-1}\bigoplus_{j=0}^{e-1}\langle u^jE(u)^i\rangle_E
$$
for $m\leq p$, we deduce that $\dim_E X=3e$, $\dim_E (X\cap Y)=2e$.
We now note that $\Gal(K/\Qp)$ acts semisimply and that multiplication by $u$ defines an $E$-linear automorphism on
$S_E/\Fil^pS_E$ which cyclically permutes isotypical components. 
The result follows.
\end{proof}

We can now state the main result of this section:

\begin{prop}
\label{Shape filtration str div modules}
Assume that $K_0=\Qp$ and let $\rho:G_{\Qp}\rightarrow \GL_3(\oe)$ be a semistable Galois representation, becoming crystalline over $K$, 
with inertial type $\teich{\omega}^{a_0}\oplus\teich{\omega}^{a_1}\oplus\teich{\omega}^{a_2}$ and Hodge--Tate weights $\{-2,-1,0\}$.
Let $\widehat{\cM}\in\OEModdd[2]$ be a strongly divisible $\oe$-module such that $\Tst^{\Qp,2}(\widehat{\cM})\cong\rho$. Assume that $\overline{\rho}:G_{\Qp}\rightarrow \GL_3(\F)$ is maximally non-split and generic as in Definition \ref{genericitycondition}.

There exist a framed basis $\widehat{\und{e}}=(\widehat{e}_0,\widehat{e}_1,\widehat{e}_2)$ of $\widehat{\cM}$, with respect to  $\teich{\omega}^{a_0}\oplus\teich{\omega}^{a_1}\oplus\teich{\omega}^{a_2}$,
and a framed system of generators $\widehat{\und{f}}=(\widehat{f}_0,\widehat{f}_1,\widehat{f}_2)$ for $\Fil^2\widehat{\cM}/\big(\Fil^p S\cdot\widehat{\cM}\big)$, whose coordinates with respect to $\wh{\und e}$ are described by one of the following possibilities.

\emph{\textbf{Case A:}}
$$
\widehat{f}_0=\,\columnvct{\alpha\\0\\u^{e-(a_2-a_0)}},\,\,
\widehat{f}_1=\,E(u)\columnvct{0\\1\\0},\,\,\widehat{f}_2=\,\columnvct{(p+E(u))u^{a_2-a_0}\\0\\-\frac{p^2}{\alpha}},
$$
where $\alpha\in \oe$ with $0<\ord(\alpha)<2$.

\emph{\textbf{Case B}:}
$$
\widehat{f}_0=\,\columnvct{\alpha\\0\\u^{e-(a_2-a_0)}},\,\,
\widehat{f}_1=\,E(u)\columnvct{u^{a_1-a_0}\\-\frac{p}{\beta}\\0},\,\,
\widehat{f}_2=\,\columnvct{\beta u^{a_2-a_0}\\E(u)u^{a_2-a_1}\\-\frac{p\beta}{\alpha}},
$$
where $\alpha,\beta\in \oe$ with $0<\ord(\beta)<1$ and $0<\ord(\alpha)<\ord(\beta)+1$.

\emph{\textbf{Case C}:}
$$
\widehat{f}_0=\,\columnvct{-\frac{p\alpha}{\beta}\\u^{e-(a_1-a_0)}\\\frac{p}{\beta}u^{e-(a_2-a_0)}},\,\,
\widehat{f}_1=\,\columnvct{-\alpha u^{a_1-a_0}\\-{\beta}\\u^{e-(a_2-a_1)}},\,\,
\widehat{f}_2=\,\columnvct{(p+E(u))u^{a_2-a_0}\\ \frac{p\beta}{\alpha}u^{a_2-a_1}\\\frac{p^2}{\alpha}},
$$
where $\alpha,\beta\in \oe$ with $0<\ord(\beta)<1$ and $0<\ord(\alpha)<\ord(\beta)+1$.
\end{prop}

Here, by a framed basis $\und{\widehat{e}}$ of $\widehat{\cM}$ we mean an $S$-basis such that $\widehat{g}\cdot\widehat{e}_i=\teich{\omega}(g)^{a_i}\widehat{e}_i$ for all $0\leq i\leq 2$, $g\in \Gal(K/K_0)$.
Similarly, we mean that $\und{\widehat{f}}$ consists of elements of $\Fil^2\widehat{\cM}$ that generate $\Fil^2\widehat{\cM}/(\Fil^p S\cdot \widehat{\cM})$ as $S$-module and such that $\widehat{g}\cdot\widehat{f}_i=\teich{\omega}(g)^{a_i}\widehat{f}_i$ for all $0\leq i\leq 2$, $g\in \Gal(K/K_0)$.
Finally, $\ord$ is the valuation of $E$ normalized by $\ord(p)=1$.

\begin{proof}
Let $\underline{e}\defeq (e_0,e_1,e_2)$ be a framed basis of $\widehat{\cM}$.
We write the elements of $\widehat{\cM}$ in terms of coordinates with respect to $\und{e}$. Moreover, we let $\cM\defeq \widehat{\cM}\otimes_S \barS$, define $\cD\defeq\widehat{\cM}\otimes_{\oe} S_E$, and set $X=\big(\Fil^2\cM/\Fil^2S\cdot \cM\big)\otimes_{\oe} E$ as in Lemma \ref{lemma shape filtration}.
In particular $\cM\in \FBrModdd[2]$ is of type ${\omega}^{a_0}\oplus{\omega}^{a_1}\oplus{\omega}^{a_2}$ by Lemma \ref{lemma str div dd}.

By Lemma \ref{lemma shape filtration} we have a non-zero element $f_0\in \big(\Fil^2\widehat{\cM}/\Fil^2S\cdot \widehat{\cM}\big)_{\teich{\omega}^{a_0}}$
of the form
$$
f_0=\columnvct{x\\ y u^{e-(a_1-a_0)}\\ z u^{e-(a_2-a_0)}}+E(u)\columnvct{x'\\ y' u^{e-(a_1-a_0)}\\ z' u^{e-(a_2-a_0)}},
$$
where $x,y,z,x',y',z'\in\oe$. As $\Fil^2\widehat{\cM}/\Fil^2S\cdot \widehat{\cM}$ is $\varpi_E$-torsion free, we may assume that one of $x,y,z,x',y',z'$ is in $\oe\s$.
If $x,y,z\in \varpi_E\oe$ for all choices of $f_0$, we get a contradiction from Corollary \ref{corollary filtration Breuil module}(v). 
On the other hand by Corollary \ref{corollary filtration Breuil module}(iv) we necessarily have $x\in \varpi_E\oe$, hence we may assume that $y\in \oe\s$ or $z\in\oe\s$.

\emph{Case 1: assume $z\in \oe\s$.}
Let us define the element $e_2'\in\widehat{\cM}$ by
$$
e_2'=\columnvct{x'u^{a_2-a_0}\\(y+y'E(u))u^{a_2-a_1}\\z+z'E(u)}.
$$
As $z\in \oe\s$ we deduce by Nakayama's lemma that $\underline{e}'\defeq (e_0,e_1,e_2')$ is again a
framed basis of $\widehat{\cM}$. By letting $\alpha \defeq x+px'$ we therefore have the following coordinates for $f_0$ in the basis $\underline{e}'$:
$$
f_0=\columnvct{\alpha\\0\\u^{e-(a_2-a_0)}},
$$
where $\ord(\alpha)>0$. (Recall that $E(u)=u^e+p$.) 
From now on we use the basis $\underline{e}'$ to write the coordinates of the elements in $\widehat{\cM}$.

By Lemma \ref{lemma shape filtration} we easily deduce that $\left(\Fil^2\widehat{\cM}/\Fil^2S\cdot\widehat{\cM}\right)_{\teich{\omega}^{a_0}}$ equals
$$
\left\langle\columnvct{\alpha\\0\\u^{e-(a_2-a_0)}},\,E(u)\columnvct{\alpha\\0\\u^{e-(a_2-a_0)}},\,
E(u)\columnvct{\beta\\\gamma u^{e-(a_1-a_0)}\\0}\right\rangle_{\oe},
$$
where $\beta, \gamma \in\oe$ are such that either $\beta\not\equiv 0$ or $\gamma\not\equiv 0$ modulo $\varpi_E$.
Moreover, by Corollary \ref{corollary filtration Breuil module}(vi) we necessarily have $\beta\equiv 0$ so that, without loss of generality, we may assume $\gamma=1$.
By Lemma \ref{lemma shape filtration} we obtain
$$
X_{\teich{\omega}^{a_1}}=
\left\langle\columnvct{\alpha u^{a_1-a_0}\\0\\u^{e-(a_2-a_1)}},\,E(u)\columnvct{\alpha u^{a_1-a_0}\\0\\u^{e-(a_2-a_1)}},\,
E(u)\columnvct{\beta u^{a_1-a_0}\\ -p \\0}\right\rangle_{E}
$$
and we need to further distinguish two subcases according to the valuation of $\beta\in\oe$.

\emph{Case 1a: assume $\beta\in p\oe$.} Then the element $e_1'$ defined by
$$
e_1'=\columnvct{-\frac{\beta}{p}u^{a_1-a_0}\\1\\0}
$$
is in $\left(\Fil^2\widehat{\cM}/\Fil^2S\cdot\widehat{\cM}\right)_{\teich{\omega}^{a_1}}$ and the family $\underline{e}''\defeq (e_0,e_1',e_2')$ is again a framed basis of $\widehat{\cM}$. Until the end of Case 1a we use the basis $\underline{e}''$ to write the coordinates of the elements in $\widehat{\cM}$.

Note that
$$
\frac{1}{\alpha}u^{a_2-a_0}(p+E(u))f_0=\columnvct{(p+E(u))u^{a_2-a_0}\\0\\-\frac{p^2}{\alpha}}\in X_{\teich{\omega}^{a_2}}.
$$
By Corollary \ref{corollary filtration Breuil module}(vii) we necessarily have $\frac{p^2}{\alpha}\in\varpi_E\oe$ and therefore
\begin{equation*}
\left\langle\columnvct{\alpha\\0\\u^{e-(a_2-a_0)}},\,E(u)\columnvct{0\\1\\0},\,
\columnvct{(p+E(u))u^{a_2-a_0}\\0\\-\frac{p^2}{\alpha}}\right\rangle_{S}\subseteq
\Fil^2\widehat{\cM}/\Fil^2S\cdot\widehat{\cM}.\nonumber
\end{equation*}
By Nakayama's lemma, noticing that the left-hand side surjects onto $\Fil^2\cM/u^{2e}\cM$ (e.g.\ by Corollary
\ref{corollary filtration Breuil module}) we conclude that the inclusion is indeed an equality. 
We also see that the elementary divisors of $\cM/\Fil^2\cM$ are those described by Corollary \ref{corollary filtration Breuil module}(i).
In fact, we can repeat the Nakayama argument to deduce that
the same three elements generate the $S$-module $\Fil^2\widehat{\cM}/\Fil^pS\cdot\widehat{\cM}$, since their images generate $\Fil^2 \cM$,
so they form a framed system of generators of this module. (Alternatively there is a direct argument using $\Fil^2 S = (E(u)^2,\Fil^p S)$.)
We therefore land in Case A. 

\emph{Case 1b: assume $p\in \beta\varpi_E\oe$.}
As above we deduce 
$$
E(u)\columnvct{u^{a_1-a_0}\\-\frac{p}{\beta}\\0}\in \left(\Fil^2\widehat{\cM}/\Fil^2S\cdot\widehat{\cM}\right)_{\teich{\omega}^{a_1}}
$$
and 
\begin{equation*}
\hspace{-2cm}\left\langle\columnvct{\alpha u^{a_2-a_0}\\0\\-p}+E(u)\columnvct{0\\0\\1},\,E(u)\columnvct{\alpha u^{a_2-a_0}\\0\\-p},
E(u)\columnvct{u^{a_2-a_0}\\-\frac{p}{\beta}u^{a_2-a_1}\\0}\right\rangle_{\oe}\subseteq
\left(\Fil^2\widehat{\cM}/\Fil^2S\cdot\widehat{\cM}\right)_{\teich{\omega}^{a_2}}.
\end{equation*}
In particular 
$$
\columnvct{\beta u^{a_2-a_0}\\u^{a_2-a_1}E(u)\\-\frac{p\beta}{\alpha}}\in X_{\teich{\omega}^{a_2}}
$$
and, by Corollary \ref{corollary filtration Breuil module}(vii) we necessarily have $\beta\in \varpi_E\oe$ and $p\beta\in\alpha\varpi_E\oe$.

We thus obtain
\begin{equation*}
\hspace{-.5cm}\left\langle\columnvct{\alpha\\0\\u^{e-(a_2-a_0)}},\,E(u)\columnvct{u^{a_1-a_0}\\-\frac{p}{\beta}\\0},\,
\columnvct{\beta u^{a_2-a_0}\\u^{a_2-a_1}E(u)\\-\frac{p\beta}{\alpha}}\right\rangle_{S}\subseteq
\Fil^2\widehat{\cM}/\Fil^2S\cdot\widehat{\cM}.
\end{equation*}
It follows, as in Case 1a, that the inclusion is actually an equality and we land in Case B. We also see that the elementary divisors of $\cM/\Fil^2\cM$ are those described by Corollary \ref{corollary filtration Breuil module}(ii).

\emph{Case 2: assume $y\in \oe\s$ and $z\in \varpi_E\oe$.} We note that the image of $f_0$ in $(\Fil^2\cM)_{\omega^{a_0}}$ is divisible by $u^{e-(a_1-a_0)}$ but not by $u^e$.  From Proposition \ref{mainsection2.2} we deduce that the elementary divisors of $\cM/\Fil^2\cM$ are those described by Corollary \ref{corollary filtration Breuil module}(iii). 
It is then not hard to see that Case B of the proposition applies to the dual strongly divisible module $\widehat{\cM}^*$ (\cite{xavier-thesis}, Chapitre V) upon interchanging $(a_0,a_1,a_2)$ by $(-a_2,-a_1,-a_0)$. We arrive at Case C by dualizing.
Alternatively, Case 2 can be treated directly with some effort.
\end{proof}

\begin{corollary}
\label{corollary Shape filtration str div modules}
Let ${\rho}$ and $\widehat{\cM}$ be as in Proposition \ref{Shape filtration str div modules}.
Write $\lambda_i$ for the Frobenius eigenvalue on $\Dst^{\Qp,2}(\rho)^{I_{\Qp} = \teich{\omega}^{a_i}}$. 
Then
\begin{equation*}
(\ord(\lambda_0),\ord(\lambda_1),\ord(\lambda_2))=\left\{\begin{matrix}
(2-\ord(\alpha),1,\ord(\alpha))\\
(2-\ord(\alpha),\ord(\beta),1+\ord(\alpha)-\ord(\beta))\\
(1+\ord(\beta)-\ord(\alpha),2-\ord(\beta),\ord(\alpha))
\end{matrix}\right.
\end{equation*}
according to whether we are in Case A, B, or C of Proposition \ref{Shape filtration str div modules}.
In particular, $\ord(\lambda_i) \in (0,2)$ for all $i$ and $\sum_{i=0}^2 \ord(\lambda_i) = 3$.
\end{corollary}

\begin{proof}
  We first note that the elements $\varphi_2(\widehat{f}_i)$ span $\widehat{\cM}$ by Nakayama's lemma, as $\varphi_2(\Fil^p
  S\cdot \wh\cM) \subset p^{p-2} \wh M \subset \mathfrak{m}_S \wh M$ and $p > 2$.

From the proof of \cite{EGH}, Proposition 3.1.4 we have 
$\Dst^{\Qp,2}(\rho)\cong \widehat{\cM}[\frac{1}{p}]\otimes_{S_{\Qp},s_0}\Qp$ where $s_0: S_{\Qp}\rightarrow \Qp$ is defined by ``$u\mapsto 0$''.
We have $\varphi_2(\widehat{f}_j)=\sum_{i=0}^2\widehat{\alpha}_{ij}\widehat{e}_i$ where $\widehat{\alpha}_{ij}\in S_{\teich{\omega}^{a_j-a_i}}$. Hence, as the elements $\varphi_2(\widehat{f}_i)$ span $\widehat{\cM}$, we see that $\widehat{\alpha}_{ii}\in S\s$ for all $0\leq i\leq 2$.
On the other hand note that $S_{\teich{\omega}^i}\subseteq \ker(s_0)$ for all $i\not\equiv 0\,\mod{p-1}$.
Let $e_i,\,f_i$ denote the images of $\widehat{e}_i,\,\widehat{f}_i$ in $\Dst^{\Qp,2}(\rho)$.
In Case A we now deduce that $\frac{\varphi}{p^2}(f_0)=\mu e_0$ where $\mu\defeq s_0(\widehat{\alpha}_{00})\in \oe\s$. In other words, $\alpha\varphi(e_0)=p^2\mu e_0$ so that $\lambda_0=p^2\mu\alpha^{-1}$ with $\ord(\lambda_0)=2-\ord(\alpha)$.
The other cases are deduced similarly.
\end{proof}

\subsection{From strongly divisible modules to Fontaine--Laffaille modules}
\label{Str div to FL}
The aim of this section is to explicitly determine the Fontaine--Laffaille module associated to the mod-$p$ reduction of a potentially semistable lift of a maximally non-split and generic Galois representation $\rhobar$. 
Using Proposition \ref{Shape filtration str div modules} we compute the \'etale $\varphi$-modules associated to the mod-$p$ reduction of strongly divisible modules, letting us compare the Frobenius eigenvalues on the weakly admissible module of certain potentially semistable lifts of $\rhobar$ with the Fontaine--Laffaille invariant $\FL(\rhobar)$.

Let $\red:\mathbb{P}^1(\oe)\rightarrow \mathbb{P}^1(\F)$  denote the specialization map. 
The main local result on the Galois side is the following theorem. It will be an important input for proving our global Theorem~\ref{main local global}.

\begin{thm}
\label{main theorem Galois}
Let $\rho: \GQp\rightarrow \mathrm{GL}_3(\oe)$ be a potentially semistable $p$-adic Galois representation of Hodge--Tate weights $\{-2,-1,0\}$ and inertial type 
$\teich{\omega}^{a_0}\oplus\teich{\omega}^{a_1}\oplus\teich{\omega}^{a_2}$.
Assume that the residual representation $\rhobar : \GQp\rightarrow \mathrm{GL}_3(\F)$ is maximally non-split and generic as in Definition \ref{genericitycondition}.
Let $\lambda_i\in \oe$ be the Frobenius eigenvalue on $\Dst^{\Qp,2}(\rho)^{I_{\Qp} = \teich{\omega}^{a_i}}$.
Then the Fontaine--Laffaille invariant of $\rhobar$ is given by:
$$
\FL(\rhobar)=\red\left(p\lambda_1^{-1}\right).
$$
\end{thm}

We first determine the filtration and the Frobenius action on the Breuil modules obtained as the base change of strongly divisible $\oe$-modules corresponding to $\oe$-lattices inside $\rho$.

\begin{prop}
\label{prop main 1}
Let $\widehat{\cM}$, $\widehat{\und{e}}$, $\widehat{\und{f}}$ be as in the statement of Proposition \ref{Shape filtration str div modules}.
Define $\alpha_i\in \F^{\times}$ via the condition $[\varphi_2(\widehat{f}_i)]=\alpha_i [\widehat{e}_i]$ in $\widehat{\cM}/\mathfrak{m}_S\widehat{\cM}$ and let $\cM\in \FBrModdd[2]$ be the base change of $\widehat{\cM}$ via $S\onto\barS$.
Then $\cM$ is of type ${\omega}^{a_0}\oplus{\omega}^{a_1}\oplus{\omega}^{a_2}$ and there exist a framed basis $\underline{e}$ of $\cM$ and a framed system of generators $\underline{f}$ for $\Fil^2\cM$ such that one of the following hold:
\begin{enumerate}
\item $\Fil^2\widehat{\cM}$ is as in Case ${A}$ of Proposition \ref{Shape filtration str div modules} and
\begin{equation*}
\mat_{\underline{e},\underline{f}}\left(\Fil^2\cM\right)=\maq{c_{00}u^e }{c_{01}u^{e+(a_1-a_0)}}{u^{e+(a_2-a_0)}}{c_{10}u^{e-(a_1-a_0)}}{u^e}{}{u^{e-(a_2-a_0)}}{}{};\ 
\mathrm{Mat}_{\underline{e},\underline{f}}(\varphi_2)=\mathrm{Diag}(\alpha_0,\alpha_1,\alpha_2)
\end{equation*}
for some $c_{ij}\in \F$.
\item $\Fil^2\widehat{\cM}$ is as in Case ${B}$ of Proposition \ref{Shape filtration str div modules} and
\begin{equation*}
\mat_{\underline{e},\underline{f}}\left(\Fil^2\cM\right)=\maq{c_{00}u^e}{u^{e+(a_1-a_0)}}{}{c_{10}u^{e-(a_1-a_0)}}{}{u^{e+(a_2-a_1)}}{u^{e-(a_2-a_0)}}{}{};\
\mathrm{Mat}_{\underline{e},\underline{f}}(\varphi_2)=\mathrm{Diag}(\alpha_1,\alpha_2,\alpha_0)
\end{equation*}
for some $c_{ij}\in \F$.
\item $\Fil^2\widehat{\cM}$ is as in Case $C$ of Proposition \ref{Shape filtration str div modules} and
\begin{equation*}
\mat_{\underline{e},\underline{f}}\left(\Fil^2\cM\right)=\maq{c_{00}u^e}{c_{01}u^{e+(a_1-a_0)}}{u^{e+(a_2-a_0)}}{u^{e-(a_1-a_0)}}{}{}{}{u^{e-(a_2-a_1)}}{};\
\mathrm{Mat}_{\underline{e},\underline{f}}(\varphi_2)=\mathrm{Diag}(\alpha_2,\alpha_0,\alpha_1)
\end{equation*}
for some $c_{ij}\in \F$.
\end{enumerate}
\end{prop}
\begin{proof}
We let $a\defeq a_2-a_0$, $b\defeq a_1-a_0$ and define $\und{e}_0\defeq \widehat{e}\otimes1$, $\und{f}_0\defeq \widehat{f}\otimes1$. 
Then $\cM$ is of type ${\omega}^{a_0}\oplus{\omega}^{a_1}\oplus{\omega}^{a_2}$, $\und{e}_0$ is a framed basis of $\cM$ and $\und{f}_0$ is a framed system of generators of $\Fil^2\cM$.
We let 
\begin{align*}
V_0&\defeq \mat_{\und{e}_0,\und{f}_0}(\Fil^2\cM),& 
A_0&\defeq \mathrm{Mat}_{\und{e}_0,\und{f}_0}(\varphi_2).
\end{align*} 
Note in particular that $\alpha_i\equiv(A_0)_{ii}\,\,\mathrm{mod}\,u^e$ for $i=0,1,2$, by construction.

We now treat separately each of the three cases $\mathrm{(i)}$, $\mathrm{(ii)}$ and $\mathrm{(iii)}$. 

(i) From Case A of Proposition \ref{Shape filtration str div modules} we deduce that
$$
V_0=\mat_{\und{e}_0,\und{f}_0}(\Fil^2\cM)=\maq{}{}{u^{e+a}}{}{u^e}{}{u^{e-a}}{}{}.
$$ 
We can write $V_{0}^{\adj}=-u^e W_0$ for some 
$W_0\in \mathrm{M}_{3}(\barS)$, which is uniquely defined modulo $u^{e(p-1)}$; in fact, we can and will take $W_0 = V_0$.
We claim that there exist $b_{00},\,b_{01},\,b_{10}\in \F$ such that 
\begin{equation}
\label{equations matrix}
W_0\cdot A_0\cdot 
\underbrace{\maq{b_{00}u^e}{b_{01}u^{e+b}}{u^{e+a}}{b_{10}u^{e-b}}{u^e}{}{u^{e-a}}{}{}}_{\defeq V_1}
= u^{2e}B_0
\end{equation}
for some $B_0\in \GL_{3}^{\Box}(\barS)$.
Indeed, if (\ref{equations matrix}) holds, one obtains 
$$
V_0B_0\equiv A_0V_1\,\,\mathrm{mod}\,\,u^{3e}
$$
(as $p\geq 5$) and by Lemma \ref{lemma base change matrix} we deduce that $V_1=\mat_{\underline{e}_1,\underline{f}_1}(\Fil^2\cM)$ is a matrix for $\Fil^2\cM$ with respect to the basis $\underline{e}_1\defeq \underline{e}_0\cdot A_0$ and a system of generators $\underline{f}_1$ of $\Fil^2\cM$, and that $A_1\defeq\mathrm{Mat}_{\underline{e}_1,\underline{f}_1}(\varphi_2)=\varphi(B_0)$ describes the Frobenius action on $\cM$ with respect to $(\und{e}_1,\und{f}_1)$.

To prove the claim, for $0\leq i,j\leq 2$ let us write $(A_0)_{ij}=u^{[a_j-a_i]}\alpha_{ij}$ with $\alpha_{ij}\in\F[u^{e}]/u^{ep}$. 
A computation, using the explicit description of the element $W_0\in \mathrm{M}_{3}^{\Box}(\barS)$, gives: 
$$
W_0\cdot A_0\cdot V_1\equiv \maq{\alpha_{22}u^{2e}}{}{}{xu^{e-b}}{\alpha_{11}u^{2e}}{}{yu^{e-a}}{zu^{e-(a-b)}}{\alpha_{00}u^{2e}}\,\,
\mathrm{mod}\,u^{3e},
$$
where $x,y,z\in \F[u^{e}]/u^{ep}$ are defined by
\begin{align*}
x&=(\alpha_{11}b_{10}+\alpha_{12})u^e+\alpha_{10}b_{00}u^{2e},&y&=(\alpha_{00}b_{00}+\alpha_{01}b_{10}+\alpha_{02})u^e\\
z&=(\alpha_{00}b_{01}+\alpha_{01})u^e.&
\end{align*}
As $\alpha_{ij}\in \barS_{\omega^0}$ and $\alpha_{ii}\in\barS^{\times}$ we deduce the existence of unique $b_{00},\,b_{10},\,b_{01}\in \F$ such that (\ref{equations matrix}) holds and, what is more important, 
$$
B_0\equiv \maq{\alpha_2}{}{}{\beta_{10}u^{e-b}}{\alpha_1}{}{\beta_{20}u^{e-a}}{\beta_{21}u^{e-(a-b)}}{\alpha_0}\,\,\mathrm{mod}\,u^e
$$
for some $\beta_{ij}\in\F$. 
In particular, we obtain 
\begin{align*}
A_1\defeq \mathrm{Mat}_{\underline{e}_1,\underline{f}_1}(\varphi_2)&=\maq{\alpha_2}{}{}{\beta_{10}u^{p(e-b)}}{\alpha_1}{}
{\beta_{20}u^{p(e-a)}}{\beta_{21}u^{p(e-(a-b))}}{\alpha_0}\\
&\equiv \mathrm{Diag}(\alpha_2,\alpha_1,\alpha_0)\,\,\mathrm{mod}\,u^{3e}
\end{align*}
by the genericity assumption $(\ref{genericityG})$.

Writing again $V_1^{\adj}= -u^eW_1$ for some $W_1\in \mathrm{M}_{3}(\barS)$ (uniquely defined modulo $u^{e(p-1)}$) we obtain
$$
W_1\cdot A_1\equiv \maq{}{}{\alpha_0u^{e+a}}{}{u^e\alpha_1}{-\alpha_0b_{10}u^{e+(a-b)}}
{\alpha_2u^{e-a}}{-\alpha_1b_{01}u^{e-(a-b)}}{\alpha_0(b_{10}b_{01}-b_{00})u^e}\,\,\mathrm{mod}\, u^{3e}
$$
and a computation shows that there exist uniquely determined elements $c_{00},\,c_{10},\,c_{01}\in\F$ such that 
\begin{equation}
\label{equations matrix 1}
W_1\cdot A_1\cdot 
\underbrace{\maq{c_{00}u^e}{c_{01}u^{e+b}}{u^{e+a}}{c_{10}u^{e-b}}{u^e}{}{u^{e-a}}{}{}}_{\defeq V_2}
= u^{2e}B_1,
\end{equation}
where now
$$
B_1\equiv\mathrm{Diag}(\alpha_0,\alpha_1,\alpha_2)\,\,\mathrm{mod}\,u^e.
$$ 
By Lemma \ref{lemma base change matrix} one has $V_2=\mat_{\underline{e}_2,\underline{f}_2}(\Fil^2\cM)$, where $\underline{e}_2\defeq \underline{e}_1\cdot A_1$, $\underline{f}_2$ is a system of generators for $\Fil^2\cM$ and the Frobenius action is now given by:
$$
A_2\defeq \mathrm{Mat}_{\underline{e}_2, \underline{f}_2}(\varphi_2)=\varphi(B_1)=\mathrm{Diag}(\alpha_0,\alpha_1,\alpha_2).
$$

(ii) The proof is similar to (i) and we content ourselves with giving the general argument, leaving the computational details to the scrupulous reader.

From Proposition \ref{Shape filtration str div modules} the filtration for $\Fil^2\cM$ now has the form
$$
V_0\defeq\maq{}{u^{e+b}}{}{}{}{u^{e+(a-b)}}{u^{e-a}}{}{}.
$$ 

Writing $V_{0}^{\adj}=u^e W_0$ for some $W_0\in \mathrm{M}_{3}(\barS)$ we claim that there exist $b_{00},\,b_{10}\in \F$ and $B_0\in \GL_{3}^{\Box}(\barS)$ such that 
\begin{equation}
\label{equations matrix1}
W_0\cdot A_0\cdot 
\underbrace{\maq{b_{00}u^e}{u^{e+b}}{}{b_{10}u^{e-b}}{}{u^{e+(a-b)}}{u^{e-a}}{}{}}_{\defeq V_1}
= u^{2e}B_0.
\end{equation}

Indeed, using the same notation as in case $\mathrm{(i)}$, an easy computation now gives: 
$$
W_0\cdot A_0\cdot V_1\equiv \maq{\alpha_{22}u^{2e}}{}{}{xu^{e-b}}{\alpha_{00}u^{2e}}
{\alpha_{01}u^{2e+(a-b)}}
{yu^{e-a}}{\alpha_{10}u^{3e-(a-b)}}{\alpha_{11}u^{2e}}\,\,
\mathrm{mod}\,u^{3e},
$$
where $x,y\in \F[u^{e}]/u^{ep}$ are defined by
\begin{align*}
x&=(\alpha_{00}b_{00}+\alpha_{01}b_{10}+\alpha_{02})u^e,&y&=(\alpha_{11}b_{10}+\alpha_{12})u^{e}+\alpha_{10}b_{00}u^{2e}.
\end{align*}
As $\alpha_{ii}\in\barS^{\times}$ we find unique $b_{00},\,b_{10}\in\F$ such that (\ref{equations matrix1}) holds and moreover
$$
B_0\equiv \maq{\alpha_{2}}{}{}{\beta_{10}u^{e-b}}{\alpha_0}{\beta_{12}u^{a-b}}{\beta_{20}u^{e-a}}{\beta_{21}u^{e-(a-b)}}{\alpha_1}\,\,\mathrm{mod}\,u^e
$$
for some $\beta_{ij}\in\F$. 
Note that
\begin{equation*}
A_1\defeq \mathrm{Mat}_{\underline{e}_1,\underline{f}_1}(\varphi_2)\equiv
\mathrm{Diag}(\alpha_2,\alpha_0,\alpha_1)\,\,\mathrm{mod}\,u^{3e}
\end{equation*}
by the genericity assumption.

By means of Lemma \ref{lemma base change matrix}  we can iterate the procedure: one finds an explicit element $V_2=\mat_{\underline{e}_2,\underline{f}_2}(\Fil^2\cM)\in \mathrm{M}_{3}(\barS)$, as in the statement of case $\mathrm{(ii)}$, which provides the filtration $\Fil^2\cM$ with respect to the basis $\underline{e}_2\defeq \underline{e}_1\cdot A_1$ and which verifies
$$
\mathrm{Mat}_{\underline{e}_2,\underline{f}_2}(\varphi_2)=\mathrm{Diag}(\alpha_1,\alpha_2,\alpha_0).
$$ 

(iii) This is analogous to (ii). More succinctly, it can be obtained from $\mathrm{(ii)}$ by duality.
\end{proof}

Thanks to Proposition \ref{prop main 1} it is now straightforward to compute the Fontaine--Laffaille invariant of $\rhobar$.

\begin{prop}
\label{prop main 2}
Let $\widehat{\cM}$, $\alpha_i\in \F\s$ be as in Proposition \ref{prop main 1}.
Let $\rho\defeq \Tst^{\Qp,2}(\widehat{\cM})$. The one of the following possibilities holds.
\begin{enumerate}
	\item $\Fil^2\widehat{\cM}$ is as in Case ${A}$ of Proposition \ref{Shape filtration str div modules} and $\FL(\rhobar)=\alpha_1^{-1}$;
\item $\Fil^2\widehat{\cM}$ is as in Case ${B}$ of Proposition \ref{Shape filtration str div modules} and $\FL(\rhobar)=0$;
\item $\Fil^2\widehat{\cM}$ is as in Case ${C}$ of Proposition \ref{Shape filtration str div modules} and $\FL(\rhobar)=\infty$.
\end{enumerate}
\end{prop}
\begin{proof}
By twisting $\rho$ by $\teich{\omega}^{-a_0}$ we may assume that $a_0=0$. Let $a\defeq a_2$, $b\defeq a_1$.

(i) 
Let $\cM\in \FBrModdd[2]$ be the base change of $\widehat{\cM}$ via $S\onto\barS$.
By Proposition \ref{prop main 1} and Lemma \ref{lemma lawdd 1}, the $\phi$-action on the 
$(\phi,\F((\underline{\varpi})))$-module 
$\mathfrak{M}\defeq M_{\Fp((\underline{\varpi}))}(\cM^{\ast})$
is given by
$$
\mathrm{Mat}_{\underline{\mathfrak{e}}}(\phi)=
\maq{c_{00}\underline{\varpi}^e}{c_{01}\underline{\varpi}^{e-b}}{\alpha_2^{-1}\underline{\varpi}^{e-a}}
{c_{10}\underline{\varpi}^{e+b}}{\alpha_1^{-1}\underline{\varpi}^e}{}
{\alpha_0^{-1}\underline{\varpi}^{e+a}}{}{}
$$
for a framed basis $\underline{\mathfrak{e}}$ (for the dual type $\omega^0\oplus\omega^{-b}\oplus\omega^{-a}$) of $\mathfrak{M}$ and some $c_{ij}\in\F$ (not those of Lemma \ref{lemma lawdd 1}).
In the basis $\underline{\mathfrak{e}}'\defeq \underline{\mathfrak{e}}\cdot \mathrm{Diag}(1,\underline{\varpi}^{b},\underline{\varpi}^{a})$, which is $\Gal(K/K_0)$-invariant, one obtains
$$
\mathrm{Mat}_{\underline{\mathfrak{e}}'}(\phi)=
\maq{c_{00}\underline{\varpi}^{e}}{c_{01}\underline{\varpi}^{e(b+1)}}{\alpha_2^{-1}\underline{\varpi}^{e(a+1)}}
{c_{10}\underline{\varpi}^{e}}{\alpha_1^{-1}\underline{\varpi}^{e(b+1)}}{}
{\alpha_0^{-1}\underline{\varpi}^{e}}{}{}
$$
proving that $\mathfrak{M}$ is the base change to $\F((\underline{\varpi}))$ of the $(\phi,\F((\underline{p})))$-module $\mathfrak{M}_0$ described by
$$
\mat_{\und{\mathfrak{e}}_0}(\phi)=\underbrace{\maq{c_{00}}{c_{01}}{\alpha_2^{-1}}{c_{10}}{\alpha_1^{-1}}{}{\alpha_0^{-1}}{}{}}_{\defeq F}\mathrm{Diag}(\underline{p},\underline{p}^{b+1},\underline{p}^{a+1})
$$
for some basis $\und{\mathfrak{e}}_0$ of $\mathfrak{M}_0$.
In the basis $\und{\mathfrak{e}}'_0\defeq \und{\mathfrak{e}}_0\cdot F$ the matrix of $\phi$ becomes $\mathrm{Diag}(\underline{p},\underline{p}^{b+1},\underline{p}^{a+1})F$.
We can now construct a Fontaine--Laffaille module $M\in \F\text{-}\mathcal{FL}^{[0,p-2]}$ with Hodge--Tate weights $\{1,b+1,a+1\}$ endowed with a basis $\underline{e}$ which is compatible with the Hodge filtration and satisfies $\mathrm{Mat}_{\underline{e}}(\phi_{\bullet})=F$.
By Lemma \ref{lemma lawdd 2} we deduce that $\mathfrak{M}_0\cong{\mathcal{F}}(M)$.

On the other hand $\Tcris^*(M)\cong\Tst^*(\cM^*)\cong{\rhobar}$ by Corollary \ref{corollary comparison1}. 
As $\rhobar$ is maximally non-split, there exists a basis $\underline{e}'$ of $M$ compatible with the Hodge filtration and such that $\mat_{\und{e}'}(\phi_{\bullet})$ is described as in Corollary \ref{corinvFL}. In other words $\und{e}$, $\und{e}'$ are related via a lower triangular unipotent matrix $A\in \GL_3(\F)$ in such a way that
\begin{equation*}
\label{fontaine-laffaille comparison}
\mathrm{Mat}_{\underline{e}'}(\phi_{\bullet})=A\cdot F=\maq{\mu_0}{x}{z}{}{\mu_1}{y}{}{}{\mu_2}
\end{equation*}
for some $x,y,\mu_i\in\F\s$, $z\in\F$. It follows that 
\begin{equation*}
\FL(\rhobar)=\frac{\det\begin{pmatrix}x&z\\\mu_1&y\end{pmatrix}}{-z}=
\frac{\det\begin{pmatrix}c_{01}&\alpha_2^{-1}\\\alpha_1^{-1}&\end{pmatrix}}{-\alpha_2^{-1}}=\alpha_1^{-1}.
\end{equation*}

(ii)
The proof is similar to case (i).
One checks that
$\mathfrak{M}\defeq M_{\Fp((\underline{\varpi}))}(\cM^{\ast})\cong\mathfrak{M}_0\otimes_{\F_p((\underline{p}))}\F_p((\underline{\varpi}))$, where the $(\phi,\F((\underline{p})))$-module $\mathfrak{M}_0$ is described by 
$$
\underbrace{\maq{c_{00}}{c_{01}}{\alpha_0^{-1}}{\alpha_1^{-1}}{}{}{}{\alpha_2^{-1}}{}}_{\defeq F}\mathrm{Diag}(\underline{p},\underline{p}^{b+1},\underline{p}^{a+1}).
$$
We deduce by Lemma \ref{lemma lawdd 2} that $\mathfrak{M}_0\cong{\mathcal{F}}(M)$ for the rank $3$ Fontaine--Laffaille module 
$M\in \F\text{-}\mathcal{FL}^{[0,p-2]}$  with Hodge--Tate weights $\{1,b+1,a+1\}$ and 
$\mathrm{Mat}_{\underline{e}}(\phi_{\bullet})=F$
for some basis $\underline{e}$ of $M$ that is compatible with the Hodge filtration. 
As in (i) we deduce that $\FL(\rhobar)=0$.

(iii) Follows by a similar argument.
\end{proof}

\begin{proof}[Proof of Theorem \ref{main theorem Galois}]
Let  $\widehat{\cM}\in \OEModdd[2]$ be such that $\Tst^{\Qp,2}(\widehat{\cM})\cong \rho$ and let $\und{\widehat{e}}$, $\und{\widehat{f}}$ be as in the statement of Proposition \ref{Shape filtration str div modules}.

In Case A, from the proof of Corollary \ref{corollary Shape filtration str div modules} we see that $\lambda_1=ps_0(\widehat{\alpha}_{11})$, where $\widehat{\alpha}_{ij}\in S$ are defined by $\varphi_2(\widehat{f}_j)=\sum_{i=0}^2\widehat{\alpha}_{ij}\widehat{e}_i$, and $p\lambda_1^{-1}\in \oe\s$. Hence $\red(p\lambda_1^{-1})=\overline{s_0(\widehat{\alpha}_{11})}^{-1}\in \F\s$. We conclude by Proposition \ref{prop main 2}, as $\overline{s_0(\widehat{\alpha}_{11})}=\alpha_1$ in the notation of that proposition.
In Case B we conclude by Proposition \ref{prop main 2}, noting that $\red(p\lambda_1^{-1})=0$, as $\ord(\lambda_1)<1$.
Case C is similar.
\end{proof}

\subsection{Elimination of Galois types}
\label{sectionWE}
The aim of this section is to perform a \emph{type elimination} for a maximally non-split and generic representation $\rhobar: G_{\Qp}\rightarrow \GL_3(\F)$ (cf.\ Definition \ref{genericitycondition}), by means of integral $p$-adic Hodge theory.
In other words we show for most tame inertial types $\tau: I_{\Qp}\rightarrow \GL_3(E)$ (of niveau 1 or 2) that $\rhobar$ cannot have a potentially crystalline lift with Hodge--Tate weights $\{-2,-1,0\}$ and type $\tau$.

The main results of this section (Propositions \ref{proposition weight elimination niveau 1}, \ref{proposition weight elimination niveau 2} and \ref{weight elimination FL}) will be used to prove results about Serre weights in \S \ref{sec:Weight elimination}. 

\subsubsection{\textbf{Niveau 1 types}}
\label{sec:niveau-1-types}

We first eliminate niveau 1 types. In particular, we have $K_0=\Qp$, $e=p-1$, and $K=\Qp(\zeta_p)$.
The result is the following:

\begin{prop}
\label{proposition weight elimination niveau 1}
Let $\rho:G_{\Qp}\rightarrow \GL_3(\oe)$ be a $p$-adic Galois representation, becoming semistable over $K$ with Hodge--Tate weights $\{-2,-1,0\}$. Assume that $\rhobar$ is maximally non-split and generic as in Definition \ref{genericitycondition}. 
Then
$$
\mathrm{WD}(\rho)\vert_{I_{\Qp}}\cong \teich{\omega}^{i}\oplus\teich{\omega}^{j}\oplus\teich{\omega}^{k},
$$
where either $(i,j,k)\equiv (a_2,a_1,a_0)$ modulo $p-1$ or $(i,j,k)\equiv (a_2-1,a_1,a_0+1)$ modulo $p-1$.
\end{prop}

Here and in the following the notation $x\equiv y$ mod $p-1$ for $x,y\in \Z^3$ is shorthand for $x_i\equiv y_i$ mod $p-1$ for all $i$.

\begin{proof}
After twisting we may assume $a_0=0$ and set $a\defeq a_2$, $b\defeq a_1$.
Let $\widehat{\cM}$ be a strongly divisible $\oe$-module such that $\Tst^{\Qp,2}(\widehat{\cM})=\rho$ and let $\cM$ be the Breuil module obtained by base change to $\barS$. 
By Proposition \ref{Proposition7Note Florian} $\cM$ admits a unique filtration by Breuil submodules, with graded pieces described by
$$
\cM:\cM_2\textbf{---}\cM_1\textbf{---}\cM_0
$$
(where the notation means that $\cM_2$ is a submodule and $\cM_0$ a quotient)
and the rank one Breuil submodules $\cM_i$ verify 
$\Tst^{2}\left(\cM_i\right)\vert_{I_{\Qp}}\cong \omega^{a_i+1}$ for $i=0,1,2$. 
It follows from \cite{EGH}, Lemma 3.3.2 that $\cM_i$ is of type $\omega^{k_i}$, where 
$k_2=a+1-\delta_2$, $k_1=b+1-\delta_1$, $k_0=1-\delta_0$ for some $\delta_i\in\{0,1,2\}$  and
$$
\Fil^2\cM_2=u^{\delta_2e}\cM_2,\,\,\Fil^2\cM_1=u^{\delta_1e}\cM_1,\,\,\Fil^2\cM_0=u^{\delta_0e}\cM_0.
$$
On the other hand we see that $\mathrm{WD}(\rho)\vert_{I_{\Qp}}\cong \teich{\omega}^{k_0}\oplus \teich{\omega}^{k_1}\oplus \teich{\omega}^{k_2}$ (by Lemma \ref{lemma str div dd}) and hence $\det\rho=\teich{\omega}^{k_0+k_1+k_2}\varepsilon^3$, as $\rho$ has Hodge--Tate weights $\{-2,-1,0\}$. It follows that $\delta_0+\delta_1+\delta_2=3$. 

We now use Lemma \ref{lemmasemisimple} to conclude the proof. 
In what follows, we write $\rhobar_{21}$ to denote the unique 2-dimensional sub-representation of $\rhobar$, and, similarly, $\rhobar_{10}$ to denote its unique 2-dimensional quotient.
If $(\delta_2,\delta_1)=(1,2)$ then $\rhobar_{21}$ splits, by Lemma \ref{lemmasemisimple}(ii).
If $\delta_1=0$ then $\Fil^2\left(\cM_1\textbf{---}\cM_0\right)$ splits as an $\barS$-module (as $\Fil^2\cM_1$ is an injective $\barS$-module, cf.\ Lemma \ref{Lemma1Note Florian}), hence $\rhobar_{10}$ splits by Lemma \ref{lemmasemisimple}(i). For the same reason, if $\delta_2=0$ then $\rhobar_{21}$ splits.

Therefore $(\delta_0,\delta_1,\delta_2)=(1,1,1)$ or  $(\delta_0,\delta_1,\delta_2)=(0,1,2)$, as claimed.
\end{proof}

\subsubsection{\textbf{Niveau 2 types}}
\label{sec:niveau-2-types}

We now eliminate niveau 2 types.
In particular, we have $k=\fpn{2}$, $K_0=\mathbb{Q}_{p^2}$ is the degree two unramified extension of $\Qp$, $e=p^2-1$, and $K=\mathbb{Q}_{p^2}(\sqrt[e]{-p})$.

\begin{lm}
\label{rank one niveau 2}
Let $\cM$ be a rank one object in $\FBrModdd[2]$, with descent data from $K$ to $\Qp$. Then there exists a generator $m\in\cM$ such that:
\begin{enumerate}
	\item $\Fil^2\cM=u^{r(p-1)}\cM$, where $0\leq r\leq 2(p+1)$;
	\item there exists $k\in\Z$ verifying $k+pr\equiv 0$ mod $p+1$ and such that $\cM$ is of type $\omega_{\varpi}^k$ and $m$ is a framed basis for $\cM$;
	\item $N(m)=0$.
\end{enumerate}
Moreover, one has
$$
\Tst^{2}(\cM)\vert_{I_{\Qp}}=\omega_2^{k+pr}.
$$
\end{lm}
\begin{proof}
Let $\sigma\in\mathrm{Gal}(K/\Qp)$ be the lift of the arithmetic Frobenius in $\mathrm{Gal}(k/\fp)$ that fixes $\varpi$, hence $\sigma^2=1$.
Then $\sigma$ interchanges the two primitive idempotents $\epsilon_0$, $\epsilon_1$, and one has $\sigma\circ g=g^p\circ\sigma$ for all $g\in\mathrm{Gal}(K/K_0)$.

By \cite{EGH}, Lemma 3.3.2 there is a generator $m\in\cM$ such that
\begin{enumerate}
	\item $\epsilon_i\Fil^2\cM=u^{r_ie}\epsilon_i\cM$ for $i=0,1$;
	\item $\widehat{g}(m)=(\omega_{\varpi}(g)^{k_0}\otimes 1)\epsilon_0m+(\omega_{\varpi}(g)^{k_1}\otimes 1)\epsilon_1m$ for all $g\in\mathrm{Gal}(K/K_0)$,
\end{enumerate}
where the integers $r_i, k_i$ verify $0\leq r_i\leq 2e$, $k_i\equiv p(k_{i-1}+r_{i-1})$ mod $e$.

As ${\sigma}$ exchanges the idempotents $\epsilon_i\in k\otimes_{\fp}\F$ we deduce that $r\defeq r_0=r_1$.
Let us consider $\cM/\left(u\cdot \cM\right)$. It is a rank one module over $k\otimes_{\fp}\F$, endowed with a semilinear action of $\mathrm{Gal}(K/\Qp)$. In particular $\sigma \overline{m}=\mu \overline m$ for some $\mu\in (k\otimes_{\fp}\F)\s$ verifying $\mu\sigma(\mu)=1$.
We therefore deduce, by (ii) above, that for all $g\in \Gal(K/K_0)$,
\begin{align*}
\widehat{\sigma}\circ\widehat{g}(\overline{m})&=\widehat{\sigma}\left(
(\omega_{\varpi}^{k_0}(g)\otimes1)\epsilon_0\overline m+(\omega_{\varpi}^{k_1}(g)\otimes1)\epsilon_1\overline m\right)\\
&=\left[(\omega_{\varpi}^{pk_0}(g)\otimes1)\epsilon_1\mu \overline m+
(\omega_{\varpi}^{pk_1}(g)\otimes1)\epsilon_0\mu\overline m\right]
\end{align*}
and
\begin{equation*}
\widehat{g^p}\circ\widehat{\sigma}(\overline{m})=\mu\left[(\omega_{\varpi}^{pk_0}(g)\otimes1)\epsilon_0 \overline m+
(\omega_{\varpi}^{pk_1}(g)\otimes1)\epsilon_1\overline m\right].
\end{equation*}
As $\sigma\circ g=g^p\circ\sigma$ we conclude that $k\defeq k_0\equiv k_1$ mod $e$.
We finally obtain $k\equiv p(k+r)$ mod $e$, so that $r\equiv 0$ mod $p-1$ and $k+\frac{pr}{p-1}\equiv 0$ mod $p+1$.
The last statement on $\Tst^{2}(\cM)$ follows again from \cite{EGH}, Lemma 3.3.2 and the previous results on $k_i,r_i$.
\end{proof}

For the remainder of this section, all Breuil modules $\cM\in\FBrModdd[2]$ have descent data from $K$ to $\Qp$.
\begin{prop}
\label{proposition weight elimination niveau 2} 
Let $\rho:G_{\Qp}\rightarrow \GL_3(\oe)$ be a $p$-adic Galois representation, becoming semistable over $K$ with Hodge--Tate weights $\{-2,-1,0\}$. 
Assume that $\rhobar$ is maximally non-split and generic as in Definition \ref{genericitycondition} and that the inertial type of $\rho$ has niveau two, i.e.\ $\rho$ does not become semistable over $\Qp(\zeta_p)$.
Then
$$
\mathrm{WD}(\rho)\vert_{I_{\Qp}}\cong \teich{\omega}^x\oplus\teich{\omega}_2^y\oplus\teich{\omega}_2^{py},
$$
where $x,y\in \Z$, $y\not\equiv0\,\mathrm{mod}\,{(p+1)}$ and the pair $(x,y)\in\Z^2$ verifies one of the following possibilities:
\begin{enumerate}
	\item	$x\equiv a_0-\delta\,\,\mathrm{mod}\,p-1$ and $y\equiv a_2+pa_1+\delta-\epsilon(p-1)\,\,\mathrm{mod}\,e$;
	\item    $x\equiv a_1-\delta\,\,\mathrm{mod}\,p-1$ and $y\equiv a_2+pa_0+\delta-\epsilon(p-1)\,\,\mathrm{mod}\,e$;
	\item     $x\equiv a_2-\delta\,\,\mathrm{mod}\,p-1$ and $y\equiv a_1+pa_0+\delta-\epsilon(p-1)\,\,\mathrm{mod}\,e$,
\end{enumerate}
where $\epsilon\in\{0,1\}$, $\delta\in\Z$ are such that $\delta+\epsilon\in\{0,1\}$.
\end{prop}

\begin{proof}As in the proof of Proposition \ref{proposition weight elimination niveau 1} we may assume $a_0=0$ and set $a\defeq a_2$, $b\defeq a_1$.
Let $\widehat{\cM}$ be a strongly divisible $\oe$-module  such that $\Tst^{\Qp,2}(\widehat{\cM})\cong\rho$; let $\cM$ be the Breuil module obtained by base change to $\barS$.
Then $\cM$ is of type ${\omega}_{\varpi}^{x(p+1)}\oplus{\omega}_{\varpi}^y\oplus{\omega}_{\varpi}^{py}$ for some $x,y\in\Z$, $y\not\equiv 0\,\,\mathrm{mod}\,p+1$ (by Lemma \ref{lemma str div dd}).
As $\rho$ has Hodge--Tate weights $\{-2,-1,0\}$ we moreover have $\det\rho=\teich{\omega}^{x+y}\varepsilon^{3}$.

By Proposition \ref{Proposition7Note Florian}, $\cM$ has a unique filtration by Breuil submodules with graded pieces described by
$$
\cM: \cM_2\textbf{---}\cM_1\textbf{---}\cM_0,
$$
where 
$\Tst^{2}(\cM_i)\vert_{I_{\Qp}}=\omega^{a_i+1}$. 
We conclude by Lemma \ref{rank one niveau 2} that one of the following three possibilities holds (up to interchanging $y$ and $py$): 
\begin{enumerate}
	\item the Breuil modules $(\cM_2,\cM_1,\cM_0)$ are respectively of type $(\omega_{\varpi}^y,\omega_{\varpi}^{py},\omega_{\varpi}^{x(p+1)})$,  where
$x\equiv 1-\delta\,\,\mathrm{mod}\,(p-1)$,  $y\equiv (a+1)(p+1)-rp\,\,\mathrm{mod}\,e$ and $py\equiv (p+1)(b+1)-sp\,\,\mathrm{mod}\,e$;
	\item the Breuil modules $(\cM_2,\cM_1,\cM_0)$ are respectively of type $(\omega_{\varpi}^{y}, \omega_{\varpi}^{x(p+1)},\omega_{\varpi}^{py})$, where 
$x\equiv b+1-\delta\,\,\mathrm{mod}\,(p-1)$,  $y\equiv (a+1)(p+1)-rp\,\,\mathrm{mod}\,e$ and $py\equiv (p+1)-sp\,\,\mathrm{mod}\,(p-1)$;
	\item the Breuil modules $(\cM_2,\cM_1,\cM_0)$ are respectively of type $(\omega_{\varpi}^{x(p+1)},\omega_{\varpi}^{y},\omega_{\varpi}^{py})$, where
$x\equiv a+1-\delta\,\,\mathrm{mod}\,(p-1)$,  $y\equiv (b+1)(p+1)-rp\,\,\mathrm{mod}\,e$ and $py\equiv (p+1)-sp\,\,\mathrm{mod}\,e$,
\end{enumerate}
where the integers $\delta,r,s$ verify $\delta\in\{0,1,2\}$, $0\leq r,s\leq 2(p+1)$ and $y+pr$, $y+s\equiv 0$ modulo $p+1$. Moreover, $r\not\equiv0 \mod{p+1}$, as $\rho$ does not become semistable over $\Qp(\zeta_p)$.

Assume we are in case $\mathrm{(i)}$. We deduce that $pr-s\equiv (a-b)(p+1)\,\,\mathrm{mod}\,e$ hence $r+s\equiv 0$ modulo $p+1$. Provided the restrictions on $r,s$, we obtain $r+s=\alpha(p+1)$ for some $\alpha\in\{1,2,3\}$ and hence $r\equiv \alpha+(a-b)\,\,\mathrm{mod}\,p-1$.
Considering the condition $\det{\rho}\vert_{I_{\Qp}}\equiv \omega^{x+y}\omega^3$, we see that the integers $\alpha$, $\delta$ verify moreover the relation $\delta+\alpha=3$ (as $p>3$).

As $2< a-b< p-3$, the condition $0< r,s< 2(p+1)$ implies that $r\in\{1+(a-b), 2+(a-b),p+1+(a-b), p+2+(a-b)\}$ and hence $(x,y)\in\left\{(-1,a+pb+1),(0,a+pb-(p-1)),(0,a+pb),(1,a+pb-p)\right\}.$
This concludes the analysis of case $\mathrm{(i)}$. 

Cases $\mathrm{(ii)}$ and $\mathrm{(iii)}$ are strictly analogous. We remark that in case $\mathrm{(ii)},\,\mathrm{(iii)}$ we obtain the condition $r\equiv \alpha+\sigma$ for $\sigma\in\{a,b\}$; as $2< \sigma< p-3$ this provides us with $r\in\{\sigma+1+\epsilon+\zeta p:\epsilon,\zeta\in\{0,1\}\}$.
\end{proof}

\begin{prop}
\label{weight elimination FL}
Let $\rho:G_{\Qp}\rightarrow \GL_3(\oe)$ be a potentially semistable $p$-adic Galois representation with Hodge--Tate weights $\{-2,-1,0\}$. 
Assume that $\rhobar$ is maximally non-split and generic as in Definition \ref{genericitycondition}. 

\emph{(i)} If 
$
\mathrm{WD}(\rho)\vert_{I_{\Qp}}\cong \teich{\omega}^{a_0}\oplus\teich{\omega}_2^{a_2+1+p(a_1-1)}\oplus\teich{\omega}_2^{a_1-1+p(a_2+1)},
$
then $\FL(\rhobar)=\infty$.

\emph{(ii)} If 
$
\mathrm{WD}(\rho)\vert_{I_{\Qp}}\cong \teich{\omega}^{a_2}\oplus\teich{\omega}_2^{a_0-1+p(a_1+1)}\oplus\teich{\omega}_2^{a_1+1+p(a_0-1)},
$
then $\FL(\rhobar)=0$.
\end{prop}

In order to prove Proposition \ref{weight elimination FL} we follow the procedure used in Section \ref{Str div to FL}: we first diagonalize the Frobenius action on the Breuil module $\cM$ and, subsequently, we compute $\FL(\rhobar)$ by extracting the Fontaine--Laffaille module from the datum of $\cM$ (via the results of Section \ref{sec:p-adic-hodge}).

The proof will occupy the remainder of this section.

\begin{proof}
We first note that (ii) follows from (i), replacing $\rho$ by $\rho^\vee$ and $(a_0,a_1,a_2)$ by $(-a_2,-a_1,-a_0)$ and keeping in mind Remark \ref{rmk FL duality}.

We will now explain the proof of (i). By twisting we may assume that $a_0=0$. We set $a\defeq a_2$, $b\defeq a_1$ and define:
\begin{align*}
k_1&\defeq a+1+p(b-1),&k_2&\defeq b-1+p(a+1),&\\
r_1&\defeq a-b+2,
&r_2&\defeq 2p-(a-b),&
\end{align*}
so that 
$\tau\defeq \mathrm{WD}(\rho)\vert_{I_{\Qp}}=
\teich{\omega}^0\oplus\teich{\omega}_2^{k_1}\oplus\teich{\omega}_2^{k_2}$. Note in particular that $0<k_1,k_2<e$ as well as the following relations, which will be useful for the computations below:
\begin{align*}
r_1(p-1)&<e,&(p-1)(r_1+r_2)&=2e,\\
k_1+r_1(p-1)&=k_2,&k_2+r_2(p-1)&=k_1+2e.
\end{align*}

Let $\widehat{\cM}\in\OEModdd[2]$ be a strongly divisible $\oe$-module such that $\Tst^{\Qp,2}(\widehat{\cM})=\rho$ and let $\cM\in \FBrModdd[2]$ be the base change of $\widehat{\cM}$ via $S\onto\barS$. In particular by Lemma \ref{lemma str div dd}, $\cM$ is of type ${\omega}_{\varpi}^0\oplus{\omega}_{\varpi}^{k_1}\oplus{\omega}_{\varpi}^{k_2}$.
The first step of the proof is to describe the filtration and the Frobenius action on $\cM$.
\begin{claim}
\label{diagonalize Frobenius N2}
There is a framed basis $\underline{e}$ of $\cM$ and a framed system of generators $\underline{f}$ for $\Fil^2\cM$ such that
\begin{align*}
\mat_{\underline{e},\underline{f}}(\Fil^2\cM)&=
\maq{u^e}{}{}{}{u^{r_1(p-1)}}{z}{yu^{2e-k_2}}{}{u^{r_2(p-1)}},&
\mathrm{Mat}_{\underline{e},\underline{f}}(\varphi_2)&=\maq{\lambda}{}{}{}{\mu}{}{}{}{\nu},
\end{align*}
where $y,z\in k\otimes_{\fp}\F$ and $\lambda,\mu,\nu\in (k\otimes_{\fp}\F)\s$.
\end{claim}
\begin{proof}
From Proposition \ref{Proposition7Note Florian}, the proof of Proposition \ref{proposition weight elimination niveau 2}, and genericity, we see that $\cM$ has a unique filtration by Breuil submodules
$$
\cM: \cM_{1}\textbf{---} \cM_{2}\textbf{---} \cM_0,
$$
where $\cM_{1}$ (resp.\ $\cM_{2}$, resp.\ $\cM_0$) is a rank one Breuil module of type $\omega_{\varpi}^{k_1}$ (resp.\ $\omega_{\varpi}^{k_2}$, resp.\ $\omega_{\varpi}^0$), and filtration described by $\Fil^2\cM_{1}=u^{r_1(p-1)}\cM_{1}$ (resp.\ $\Fil^2\cM_{2}=u^{r_2(p-1)}\cM_{1}$, resp.\ $\Fil^2\cM_{0}=u^{e}\cM_{0}$).

We therefore deduce the existence of a framed basis $\underline{e}_0\defeq(e_0,e_1,e_2)$, and a framed system of generators $\underline{f}_0=(f_0,f_1,f_2)$ of $\Fil^2\cM$ such that
\begin{align*}
V_0&\defeq \mat_{\und{e}_0,\und{f}_0}\left(\Fil^2\cM\right)=\maq{u^e}{}{}{x_0u^{e-{k_1}}}{u^{r_1(p-1)}}{z_0}{y_0'u^{e-k_2}}{}{u^{r_2(p-1)}},\\
A_0&\defeq \mathrm{Mat}_{\underline{e}_0,\und{f}_0}(\varphi_2)=
\maq{\lambda_0}{}{}{\alpha_0u^{e-k_1}}{\mu_0}{\gamma_0u^{r_1(p-1)}}{\beta_0u^{e-k_2}}{}{\nu_0},
\end{align*}
where $x_0,y_0',z_0,\alpha_0, \beta_0,\gamma_0\in \barS_{\omega_{\varpi}^0}$ and $\lambda_0,\mu_0,\nu_0\in  \barS^{\times}_{\omega_{\varpi}^0}$.

\begin{claim}
\label{splitting lemma niveau 2}
We have $y_0'=u^ey_0$ for some $y_0\in\barS_{\omega_{\varpi}^0}$.
\end{claim}
\begin{proof}
Let 
$
\cN: \cM_{2}\textbf{---}\cM_0
$
denote the rank two quotient of $\cM$.
By the previous analysis of $\cM$, we deduce that there is a basis $(\overline{e}_0,\overline{e}_2)$ of $\cN$ and a pair $(\overline f_0,\overline f_1)$ of elements of $\Fil^2\cN$ such that
\begin{equation*}
\Fil^2\cN=\left\langle \overline f_0\defeq\columnvct{u^e\\y_0'u^{e-k_2}},\,\overline f_1\defeq \columnvct{0\\u^{r_2(p-1)}}\right\rangle_{\barS}
\end{equation*}
(in terms of the basis $(\overline e_0,\overline e_2)$).
If $y_0'\not\equiv 0$ modulo $u$, then it is easy to see that $u^{2e}\overline e_0\not\in \Fil^2\cN$, contradicting the Breuil module axioms.
Therefore the claim follows.
\end{proof}

We now complete the procedure to diagonalize the Frobenius action. An elementary computation shows that there exist $y_1,z_1,\zeta_1\in \barS_{\omega_{\varpi}^0}$ such that the following equality holds:
$$
A_0 \maq{u^e}{}{}{}{u^{r_1(p-1)}}{z_1}{y_1u^{2e-k_2}}{}{u^{r_2(p-1)}}\equiv
V_0
\underbrace{\maq{\lambda_0}{}{}{\zeta_1u^{e-k_2}}{\mu_0}{}{}{}{\nu_0}}_{\defeq B_0}\,\,\mathrm{modulo}\,\,u^{3e}.
$$
Indeed, it suffices to take $y_1\equiv \nu^{-1}_0(\lambda_0 y_0-\beta_0)$ modulo $u^{2e}$, $z_1\equiv  \mu_0^{-1}(\nu_0 z_0-\gamma_0u^{2e})$ modulo $u^{3e}$ and
$\zeta_1\equiv -\lambda_0 x_0+(\alpha_0+y_1\gamma_0)u^e$ modulo $u^{3e}$.

We use now Lemma \ref{lemma base change matrix} to deduce that, in the basis $\underline{e}_1\defeq\und{e}_0\cdot A_0$, the filtration and the Frobenius action can respectively be described as:
\begin{align*}
V_1&\defeq\maq{u^e}{}{}{}{u^{r_1(p-1)}}{z_1}{y_1u^{2e-k_2}}{}{u^{r_2(p-1)}},&
A_1&\defeq
\maq{\varphi(\lambda_0)}{}{}{\varphi(\zeta_1)u^{p(e-k_2)}}{\varphi(\mu_0)}{}{}{}{\varphi(\nu_0)},
\end{align*}
where $y_1,z_1\in \barS_{\omega_{\varpi}^0}$ and 
$A_1=\varphi(B_0)$.

We write $\und{f}_1\defeq \und{e}_1\cdot V_1$. As $r_1(p-1)<e$, $(3e-k_2)-r_2(p-1)=e-k_1>r_1(p-1)$, we can find a matrix in $\GL_3^\Box(\barS)$ of the form
$$
C\defeq \begin{pmatrix}
1&&\\c_{10}u^{e-k_1-r_1(p-1)} &1&c_{12}u^{e-r_1(p-1)}\\c_{20}u^{3e-k_2-r_2(p-1)} &&1
\end{pmatrix}
$$
such that, relative to $\und{e}_1$ and $\und{f}_1'\defeq \und{f}_1\cdot C$, the filtration and the Frobenius action can be described by
\begin{equation*}
V_1'\defeq \maq{u^e}{}{}{}{u^{r_1(p-1)}}{z_1'}{y_1'u^{2e-k_2}}{}{u^{r_2(p-1)}},\,\,
A_1'\defeq
\maq{\varphi(\lambda_0)}{}{}{a'_{10}u^{p(e-k_2)}}{\varphi(\mu_0)}{a'_{12}u^{p(e-r_1(p-1))}}{a'_{20}u^{p(e-k_1)}}{}{\varphi(\nu_0)},
\end{equation*}
where $y_1',z_1'\in k\otimes_{\Fp}\F$ and $a'_{ij}\in \barS_{\omega_{\varpi}^0}$.
We can write $a'_{10}u^{p(e-k_2)}=\alpha_1'u^{e-k_1}$, $a'_{20}u^{p(e-k_1)}=\beta_1'u^{e-k_2}$,
$a'_{12}u^{p(e-r_1(p-1))}=\gamma_1'u^{r_1(p-1)}$ with $\alpha_1',\beta_1',\gamma_1'\in u^{2e}\barS_{\omega_{\varpi}^0}$ by genericity.
Writing $\lambda_1'\defeq \varphi(\lambda_0)$, $\mu_1'\defeq \varphi(\mu_0)$, $\nu_1'\defeq \varphi(\nu_0)$, which all lie in $k\otimes_{\Fp} \F$,  we see that $(V_1', A_1')$ is of the same form as $(V_0, A_0)$, so we can apply the above diagonalization procedure to obtain $(V_2,A_2)$ as in the statement of the claim.
\end{proof}

We now compute the Fontaine--Laffaille invariant of the Galois representation associated to $\cM$.
\begin{claim}
\label{fontaine-laffaille N2}
With the hypotheses and notation of Claim \ref{diagonalize Frobenius N2} we have
$$
\FL(\Tst^{2}(\cM))=\infty.
$$
\end{claim}
\begin{proof} 
Consider the $(\phi,k\otimes_{\Fp}\F((\underline{\varpi})))$-module  $\mathfrak{M}\defeq M_{k((\underline{\varpi}))}(\cM^{\ast})$.
From Claim \ref{diagonalize Frobenius N2} and Lemma \ref{lemma lawdd 1} $\mathfrak{M}$ is endowed with a framed basis $\underline{\mathfrak{e}}$ (for the type $\omega_\varpi^0\oplus\omega_\varpi^{-k_1}\oplus\omega_\varpi^{-k_2}$) such that 
$$
\mathrm{Mat}_{\underline{\mathfrak{e}}}(\phi)=
\maq{\lambda^{-1}\underline{\varpi}^e}{y\nu^{-1}\underline{\varpi}^{2e-k_2}}{}{}{\nu^{-1}\underline{\varpi}^{r_2(p-1)}}{z\mu^{-1}}{}{}{\mu^{-1}\underline{\varpi}^{r_1(p-1)}},
$$
where $y,z\in k\otimes_{\fp}\F$, $\lambda,\mu,\nu\in (k\otimes_{\fp}\F)\s$.

By considering the $\Gal(K/K_0)$-invariant basis $\underline{\mathfrak{e}}'\defeq \underline{\mathfrak{e}}\cdot \mathrm{Diag}(1,\underline{\varpi}^{k_1},\underline{\varpi}^{k_2})$ one deduces that $\mathfrak{M}\cong\mathfrak{M}_0\otimes_{\F((\underline{p}))}\F((\underline{\varpi}))$, where the $(\phi,k\otimes_{\Fp}\F((\und{p})))$-module $\mathfrak{M}_0$ is described (in some basis $\underline{\mathfrak{e}}_0$) by
\begin{equation}
\label{explicit FL niveau 2}
\mathrm{Mat}_{\und{\mathfrak{e}}_0}(\phi)={\maq{\lambda^{-1}}{y\nu^{-1}}{}{}{\nu^{-1}}{z\mu^{-1}}{}{}{\mu^{-1}}}\mathrm{Diag}(\underline{p},\underline{p}^{b+1},\underline{p}^{a+1}).
\end{equation}
From (\ref{explicit FL niveau 2}), Lemma \ref{lemma lawdd 2} (and a change of basis over $k\otimes_{\Fp}\F$ as in the proof of Proposition \ref{prop main 2}) we deduce that the Frobenius action on the Fontaine--Laffaille module $M_k$ associated to $\rhobar\vert_{G_{K_0}}$ is given by
$$
F_k\defeq\mathrm{Mat}_{\underline{e}}(\phi_{\bullet})=\varphi\maq{\lambda^{-1}}{y\nu^{-1}}{}{}{\nu^{-1}}{z\mu^{-1}}{}{}{\mu^{-1}}\in\GL_3(k\otimes_{\fp}\F)
$$
for some basis $\underline{e}$ of $M_k$ compatible with the Hodge filtration on $M_k$ (and $\varphi$ denotes the relative Frobenius on $k\otimes_{\Fp} \F$).
By Corollary \ref{corollary comparison1} we have $\Tcris^*(M_k)\cong \rhobar|_{G_{K_0}}$. In particular, 
we have $M_k\cong M\otimes_{\Fp}k$ for some Fontaine--Laffaille module $M\in\F\text{-}\mathcal{FL}^{[p-2]}$ such that $\Tcris^*(M)\cong\rhobar$ as $G_{\Qp}$-representations.

As $p>2$, $\rhobar|_{G_{K_0}}$ is still maximally non-split, hence $M_k$ is endowed with a unique filtration  $0\subsetneq M_{k,0}\subsetneq M_{k,1}\subsetneq M_{k,2}=M_k$, obtained from the cosocle filtration on $\rhobar\vert_{G_{K_0}}$. By unicity, we have $M_{k,i}= M_i\otimes_{\fp}k$ for all $i=0,1,2$, where $0\subsetneq M_{0}\subsetneq M_{1}\subsetneq M_{2}=M$ is the filtration on $M$ obtained via $\Tcris^*$ from the cosocle filtration on $\rhobar$. 

As $\Fil^{1}M_k\cap M_{k,0}$, $\Fil^{b+1}M_k\cap M_{k,1}$, $\Fil^{a+1}M_k\cap M_{k,2}$ are all  free of rank one over $k\otimes_{\fp}\F$ it follows that any change of basis of $M_k$ preserving  both the Hodge and the submodule filtration needs to be diagonal.

Therefore, by letting $F\in\GL_3(\F)$ denote the matrix of the Frobenius action on $M$, we deduce that $F_{02}=0$ if and only if $(F_k)_{02}=0$, which shows that $\FL(\rhobar)=\infty$, as claimed.
\end{proof}

From the Claims \ref{diagonalize Frobenius N2}, \ref{fontaine-laffaille N2}, the proof of Proposition \ref{weight elimination FL}(i) follows.
\end{proof}

\section{The local automorphic side}
\label{sec:LAS}

In this section we analyze the $\GL_3$-side. In particular, we establish Proposition~\ref{intro:main mod p GL3} of the introduction.
We have \cite{JantzenRAG} as a main reference for the notation and terminology.

Let $G\defeq {\GL_n}_{/\zp}$ and let $T$ be the maximal split torus consisting of diagonal matrices.
We let $\Phi$ denote the set of roots with respect to $T$ and $B\supseteq T$ the Borel subgroup of upper-triangular matrices. We write $\Delta \subseteq \Phi$ for the simple roots associated to the pair $(B,T)$ and $\langle \cdot,\cdot\rangle$ for the natural pairing on $X^{\ast}(T)\times X_{\ast}(T)$. The character and cocharacter groups  are identified with $\Z^n$ in the usual way. We let $\le$ denote the dominance order on $X^\ast (T)$.
Finally, we let $\overline{G},\,\overline{B},\dots$ denote the base change of $G,\,B,\dots$ via $\Zp\onto\Fp$. 

The Weyl group $W_G$ of $G$ is canonically isomorphic to the Weyl group of $\overline{G}$.
We write $w_0\in W_G$ for the longest element.

For any dominant character $\lambda\in X^{\ast}(\overline{T})$ let 
\begin{equation*}
H^0(\lambda)\defeq \left(\Ind_{\overline{B}}^{\overline{G}}w_0\lambda\right)^{\mathrm{alg}} \otimes_{\fp} \F
\end{equation*}
be the associated dual Weyl module. It is an algebraic representation of $\overline{G}$ (or more precisely of $\overline{G}_{/\F}$) and we write $F(\lambda)\defeq \mathrm{soc}_{\overline{G}}\left(H^0(\lambda)\right)$ for its irreducible socle.
If the weight $\lambda$ is moreover $p$-restricted, i.e.\ if $0\leq \langle \lambda,\alpha^{\vee}\rangle\leq p-1$ for all $\alpha \in \Delta$, then $F(\lambda)$ is irreducible as $\overline{G}(\fp)$-representation (see for example \cite{florian-duke}, Corollary 3.17). If $Q$ is an algebraic representation of $\overline{G}$ and $\nu \in X^{\ast}(T)$ we denote by $Q_\nu$ the $\nu$-weight space for the action of $\overline{T}$ on $Q$, as usual.

Note that if $n=1$ and $a\in \Z$ then $F(a)\cong \omega^a$, where $\omega$ is seen as a character of $\GL_1(\Fp)$ via the Artin reciprocity map (so $\omega$ is the inclusion $\fp\s \to \F\s$).

Let us specialize to the case $n=3$. We set $\alpha_1\defeq (1,-1,0)$, $\alpha_{2}\defeq (0,1,-1)$, so that $\Delta = \{\alpha_1,\alpha_2\}$, and we fix the following lifts $\dot{s}_1,\,\dot{s}_2\in G(\zp)$ of the simple reflections $s_1,s_2\in W_G$
corresponding to $\alpha_1$, $\alpha_2$:
\begin{equation*}
{\dot{s}}_1\defeq \text{\tiny$\maq{}{1}{}{1}{}{}{}{}{1}$},\quad
{\dot{s}}_2\defeq \text{\tiny$\maq{1}{}{}{}{}{1}{}{1}{}$}.
\end{equation*}
Note that $\dot{w}_0\defeq \dot{s}_1\dot{s}_2\dot{s}_1$ is then a lift of $w_0\in W_G$.

Finally, for each $\alpha\in\Phi$ we denote by $U_{\alpha}$ the associated root subgroup and by $u_{\alpha}:\mathbb{G}_a\stackrel{\sim}{\longrightarrow}U_{\alpha}$ an isomorphism as in \cite{JantzenRAG}, II.1.2. Explicitly we will take:
\begin{equation*}
u_{\alpha_1}(x)=\text{\tiny$\maq{1}{x}{}{}{1}{}{}{}{1}$},\
u_{\alpha_2}(x)=\text{\tiny$\maq{1}{}{}{}{1}{x}{}{}{1}$},\
u_{\alpha_1+\alpha_2}(x)=\text{\tiny$\maq{1}{}{x}{}{1}{}{}{}{1}$}.
\end{equation*}

Let $K \defeq G(\Zp)$. (This should not cause any confusion with the field $K$ of Sections~\ref{sec:local-galois-side} and \ref{sec:appendix}.)
We write $I$ for the Iwahori subgroup of $K$ which is the preimage of $\overline{B}(\fp)$ under the reduction map $K \onto \o G(\fp)$ and 
similarly $I_1\leq I$ for the pro-$p$ Iwahori subgroup which is the preimage of $\overline{U}(\fp)$.
If $V$ is a representation of $K$ over $\cO_E$ and $a_i \in \Z$ we write 
\[ V^{I,(a_2,a_1,a_0)} \defeq \Hom_I(\cO_E(\teich{\omega}^{a_2}\otimes\teich{\omega}^{a_1}\otimes\teich{\omega}^{a_0}), V), \] the largest $\cO_E$-submodule of $V$ on which $I$ acts via 
$\teich{\omega}^{a_2}\otimes\teich{\omega}^{a_1}\otimes\teich{\omega}^{a_0}$. In particular, this makes sense when $V$ is killed by $\varpi_E$.
Also, if $V$ is a representation of $\overline G(\fp)$ over $\F$ we write $V^{\overline T(\fp),(a_2,a_1,a_0)}$ to denote the largest subspace on which $\overline T(\fp)$ acts as $F(a_2)\otimes F(a_1)\otimes F(a_0)$, and we note that $V^{I,(a_2,a_1,a_0)} \subseteq V^{\overline T(\fp),(a_2,a_1,a_0)}$.

\subsection{Group algebra operators} \label{mod p group algebra operators}

We introduce certain mod-$p$ group algebra operators for $\overline{G}(\fp)$ and study their effect on extensions of $\overline{G}(\fp)$-representations.

Let $(a_2,a_1,a_0)\in\Z^3$ be a triple verifying the genericity condition (\ref{genericityG}).
Note that in this case the weight $(a_2,a_1,a_0)\in \Z^3$ is in particular $p$-restricted. 
We define the following elements of $\F[\overline{G}(\fp)]$:
\begin{equation}
\label{operators char p}
\begin{split}
  S&\defeq
  \underset{x,y,z\in\fp}{\sum}x^{p-(a_2-a_0)}z^{p-(a_1-a_0)}\text{\tiny$\maq{1}{x}{y}{}{1}{z}{}{}{1}$}
  \dot{w}_0,
    \\
  S'&\defeq
  \underset{x,y,z\in\fp}{\sum}x^{p-(a_2-a_1)}z^{p-(a_2-a_0)}\text{\tiny$\maq{1}{x}{y}{}{1}{z}{}{}{1}$}
  \dot{w}_0.
  \end{split}
\end{equation}
In other words,
$S=\underset{x,y,z\in\fp}{\sum}x^{p-(a_2-a_0)}z^{p-(a_1-a_0)}u_{\alpha_2}(z)u_{\alpha_1+\alpha_2}(y)u_{\alpha_1}(x)\dot{w}_0$ and similarly for $S'$.

The following property of the operators $S,\,S'$ will be crucial for us.
\begin{prop}
\label{main mod p GL3}
Let $(a_2,a_1,a_0)\in\Z^3$ be a triple satisfying \(\ref{genericityG}\). Consider the associated operators $S,\,S'\in\F[\overline{G}(\fp)]$.
\begin{enumerate}
	\item There is a unique non-split extension of irreducible $\overline{G}(\fp)$-representations
$$
0\rightarrow F(a_2-1,a_1,a_0+1)\rightarrow V\rightarrow F(a_1+p-1,a_2,a_0)\rightarrow0
$$
and $S$ induces an isomorphism $S:\,V^{I,(a_1,a_2,a_0)}\stackrel{\sim}{\longrightarrow}V^{I,(a_2-1,a_1,a_0+1)}$
of one-dimensional vector spaces.
\item There is a unique non-split extension of irreducible $\overline{G}(\fp)$-representations
$$
0\rightarrow F(a_2-1,a_1,a_0+1)\rightarrow V\rightarrow F(a_2,a_0,a_1-p+1)\rightarrow0
$$
and $S'$ induces an isomorphism $S':\,V^{I,(a_2,a_0,a_1)}\stackrel{\sim}{\longrightarrow}V^{I,(a_2-1,a_1,a_0+1)}$ of one-dimensional vector spaces.
\end{enumerate}
\end{prop}

More generally, one could ask (also for more general groups):

\begin{question}\label{qu:group-algebra-oper}
  Suppose that $V$ is an $\overline{G}(\fp)$-representation that is generated by an $I$-eigenvector $v \in V^{I,(b_2,b_1,b_0)}$.   Given an $I$-eigenvector $v' \in V^{I,(c_2,c_1,c_0)}$,   find an explicit element of $\F[\overline{G}(\fp)]$ sending $v$ to $v'$.
\end{question}

By Frobenius reciprocity, $V$ is a quotient of the principal series $\Ind_{\overline{B}(\fp)}^{\overline{G}(\fp)}\big(F(b_2)\otimes F(b_1)\otimes F(b_0)\big)$,
and an interesting special case of the above question is when $V$ has irreducible socle that is moreover generated by $v'$.
Proposition \ref{main mod p GL3} answers this in some instances. (See also the proof of \cite{BD}, Proposition 2.6.1 for $\GL_2$.)

The proof of Proposition \ref{main mod p GL3} requires a certain amount of preliminaries and will occupy the rest of this section.

Let $\overline{P}\supseteq \overline{B}$ denote the standard parabolic associated to the simple root $\alpha_1$. We consider the equivariant monomorphism
\begin{equation*}
\tau\defeq{\Ind_{\overline{P}(\fp)}^{\overline{G}(\fp)}\big(F(a_1+p-1,a_2)\otimes F(a_0)\big)}\into
{\Ind_{\overline{B}(\fp)}^{\overline{G}(\fp)}\big(F(a_2)\otimes F(a_1)\otimes F(a_0)\big)}\defeq\sigma
\end{equation*}
obtained by parabolically inducing the injection $F(a_1+p-1,a_2)\into \Ind_{B_2(\fp)}^{\GL_2(\fp)}\big(F(a_2)\otimes F(a_1)\big)$ coming from Frobenius reciprocity (where $B_2 \subseteq \GL_2$ denotes the upper-triangular Borel subgroup).

As in \cite{EGH}, proof of Lemma 6.1.1, we check that
$$
\mathrm{soc}(\tau)\cong F(a_0+p-1,a_1,a_2-p+1),\ 
\mathrm{cosoc}(\tau)\cong F(a_1+p-1,a_2,a_0)
$$ 
and, again by \cite{EGH}, Lemma 6.1.1, we see that $\mathrm{rad}(\tau)/\mathrm{soc}(\tau)$ has length 4, with constituents
\begin{multline*}
F(a_1,a_0,a_2-p+1),\,\,F(a_0+p-2,a_2,a_1+1),\\F(a_2-1,a_1,a_0+1),\,\,F(a_1-1,a_0,a_2-p+2).
\end{multline*}
By \cite{Andersen87}, \S 4 we deduce that $\mathrm{rad}(\tau)/\mathrm{soc}(\tau)$ is semisimple, in particular the extension $V$ exists and we have a surjection $\tau\onto V$. Also, by \emph{loc.\ cit.}, the extension $V$ is unique, as
\begin{multline*}
\mathrm{Ext}^1_{\overline{G}(\fp)}(F(a_1+p-1,a_2,a_0),F(a_2-1,a_1,a_0+1))\\ \stackrel{\sim}{\longrightarrow}
\left(\mathrm{Ext}^1_{\mathrm{SL}_3(\fp)}(F(a_1+p-1,a_2,a_0),F(a_2-1,a_1,a_0+1))\right)^{\fp\s}
\end{multline*}
has dimension at most one.

We also see by \cite{florian-inv}, Lemma 2.3 that  $V^{I,(a_2-1,a_1,a_0+1)}$ and $V^{I,(a_1,a_2,a_0)}\cong \tau^{I,(a_1,a_2,a_0)}\cong\sigma^{I,(a_1,a_2,a_0)}$ are one-dimensional.

Let $W\defeq H^0(\lambda)$, where $\lambda\defeq (a_0+p-1,a_1,a_2-p+1)$.
We recall that $0 \to F(\lambda) \to W \to F(a_2-1,a_1,a_0+1) \to 0$ (\cite{florian-duke}, Proposition 3.18)
and moreover that this sequence is non-split even on the level of $\overline{G}(\fp)$-representations (\cite{bib:Humph_LieType}, Theorem 5.9).
The natural $\overline{G}(\fp)$-linear evaluation map 
$$
f : W=\left(\Ind_{\overline{B}}^{\overline{G}}(w_0\lambda)\right)^{\text{alg}}
\longrightarrow {\Ind_{\overline{B}(\fp)}^{\overline{G}(\fp)}(w_0\lambda)}=\sigma
$$
is injective. (To see this, it is enough to check that the restriction of $f$ to $\mathrm{soc}_{\overline{G}(\fp)}(W)\cong F(\lambda)$
is injective, and hence enough to check that the composite of $f$ followed by evaluation at $\dot w_0 \in \overline{G}(\fp)$ is injective on
the highest weight space of $\mathrm{soc}_{\overline{G}(\fp)}(W)$. The last statement is true by the proof of \cite{JantzenRAG}, Proposition II.2.2(a).)
As $\sigma$ is multiplicity free we obtain an injection $W\into \tau$.

\begin{lm}
\label{Lemma1 Nota1 Florian}
We have $S\left(\tau^{I,(a_1,a_2,a_0)}\right)\subseteq \mathrm{rad}(\tau)$.
\end{lm}
\begin{proof}
It is equivalent to show that the operator $S$ kills the highest weight space $F(\nu)^{I,(a_1,a_2,a_0)}=F(\nu)_{\nu}$, where $\nu\defeq(a_1+p-1,a_2,a_0)$.

Recall that if $Q$ is an algebraic $\overline G$-module, $v_{\nu}\in Q_{\nu}$ and $\alpha\in \Phi$ is a root then (cf.\ \cite{JantzenRAG}, II.1.19, (5) and (6))
$$
u_{\alpha}(t)v_{\nu}=\sum_{i\in \N}t^iv_{\nu+i\alpha},
$$
where $v_{\nu+i\alpha}\in Q_{\nu+i\alpha}$ (and $0^0 = 1$).

Applying this with $Q=F(\nu)$ we obtain
\begin{equation*}
S(v_{\nu})=\underset{i,j,k\geq 0}{\sum}\left(\underset{x,y,z\in\fp}{\sum}
x^{i+p-(a_2-a_0)}y^jz^{k+p-(a_1-a_0)}\right)v_{i,j,k},
\end{equation*}
where $v_{i,j,k} \in F(\nu)_{(a_0,a_2,a_1+p-1)+(i+j,-i+k,-j-k)}$.
We deduce that if the inner sum is non-zero for a fixed triple $(i,j,k)\in \N^3$, then we necessarily have $i\geq a_2-a_0-1$, $j\geq p-1$ and $k\geq a_1-a_0-1$. But then $(a_0,a_2,a_1+p-1)+(i+j,-i+k,-j-k)\not\leq (a_1+p-1,a_2,a_0)$ and the corresponding weight space in $F(\nu)$ is zero.
\end{proof}

\begin{remark}\label{rk:isotypic}
By an immediate computation one checks that $S\left(\tau^{I,(a_1,a_2,a_0)}\right) \subseteq \tau^{\overline{T}(\fp),(a_2-1,a_1,a_0+1)}$.
\end{remark}

\begin{lm}
\label{Lemma2 Nota1 Florian}
We have $S(\tau^{I,(a_1,a_2,a_0)})\subseteq W$.
\end{lm}
\begin{proof}
By Lemma \ref{Lemma1 Nota1 Florian} and Remark~\ref{rk:isotypic}, it is enough to show that the weight spaces of $F(a_1,a_0,a_2-p+1)$, $F(a_0+p-2,a_2,a_1+1)$ and of $F(a_1-1,a_0,a_2-p+2)$ do not afford the $\overline{T}(\fp)$-character $F(a_2-1)\otimes F(a_1)\otimes F(a_0+1)$.

By using the dominance order $\leq$ on $X^{\ast}(\overline{T})$ we see that the weights in $F(a_1,a_0,a_2-p+1)$ are of the form $(\alpha,\beta,\gamma)$ with $\alpha,\beta,\gamma\in\{a_2-p+1,\dots,a_1-1,a_1\}$. 
Hence $a_2-p<\alpha<a_2-1$ and $\alpha\not\equiv a_2-1$ modulo $(p-1)$.
With a similar argument we see that any weight $(\alpha,\beta,\gamma)$ appearing in $F(a_0+p-2,a_2,a_1+1)$, $F(a_1-1,a_0,a_2-p+2)$ satisfies $\beta\not\equiv a_1$ modulo $(p-1)$. 
\end{proof}

\begin{lm}
\label{Lemma3 Nota1 Florian}
We have $S(\tau^{I,(a_1,a_2,a_0)})\neq 0$.
\end{lm}
\begin{proof}
Consider $\tau\into \sigma$. Then $\tau^{I,(a_1,a_2,a_0)}$ consists of functions $f\in\sigma$ such that
\begin{enumerate}
	\item $\mathrm{supp}(f)\subseteq \overline{B}(\fp)\dot s_1 \overline{B}(\fp)$;
	\item $f$ is constant on $\dot s_1\text{\tiny$\maq{1}{\ast}{\ast}{}{1}{\ast}{}{}{1}$}$.
\end{enumerate}
Take any non-zero $f\in \tau^{I,(a_1,a_2,a_0)}$ and let $g\defeq\text{\tiny$\maq{}{}{1}{}{1}{-1}{1}{-1}{}$}$.
It then follows immediately from the definitions that
$$
S(f)(g)=\underset{x,y,z\in\fp}{\sum}x^{p-(a_2-a_0)}z^{p-(a_1-a_0)}f\text{\tiny$\maq{1}{}{}{z-1}{1}{}{y-z}{x-1}{1}$}.
$$
We deduce from (i) that the term inside the sum is non-zero iff $x=1$, $z=y$, $z\neq 1$.
Hence,
\begin{align*}
S(f)(g)&=\underset{z\in\fp\s}{\sum}(z+1)^{p-(a_1-a_0)}f\text{\tiny$\maq{1}{}{}{z}{1}{}{}{}{1}$}\\
&=(-1)^{a_2}\underset{z\in\fp\s}{\sum}(z+1)^{p-(a_1-a_0)}z^{-(a_2-a_1)}f(\dot s_1),
\end{align*}
where the last equality follows from $\mathrm{(i)}$, $\mathrm{(ii)}$ via the relation
$$
\text{\tiny$\maq{1}{}{}{z}{1}{}{}{}{1}$}=\text{\tiny$\maq{-z^{-1}}{1}{}{}{z}{}{}{}{1}$}\dot{s}_1
\text{\tiny$\maq{1}{z^{-1}}{}{}{1}{}{}{}{1}$}.
$$
Note that the term inside the sum is a linear combination of monomials $z^i$ with $\vert i\vert<p-1$. We easily deduce:
\begin{equation*}
S(f)(g)=(-1)^{a_2-1}\binom{p-(a_1-a_0)}{a_2-a_1}f(\dot{s}_1)\neq 0.\qedhere
\end{equation*}
\end{proof}

\begin{lm}
\label{Lemma4 Nota1 Florian}
We have $S\left(V^{I,(a_1,a_2,a_0)}\right)\subseteq V^{I,(a_2-1,a_1,a_0+1)}$.
\end{lm}
\begin{proof}
This is clear from Lemma \ref{Lemma1 Nota1 Florian} and Remark~\ref{rk:isotypic}, as 
(just as in Lemma~\ref{Lemma2 Nota1 Florian}) we have $\mathrm{soc}(V)^{\overline{T}(\fp),(a_2-1,a_1,a_0+1)} = \mathrm{soc}(V)^{\overline{U}(\fp)}$.
\end{proof}

\begin{proof}[Proof of Proposition \ref{main mod p GL3}]
We start by proving (i). Recall that we defined $W\defeq H^0(\lambda)$ with $\lambda\defeq (a_0+p-1,a_1,a_2-p+1)$.
By Lemma \ref{Lemma4 Nota1 Florian} it is enough to show that $S(V^{I,(a_1,a_2,a_0)})\neq0$ or, equivalently (cf.\ Lemma \ref{Lemma2 Nota1 Florian}), that $S(\tau^{I,(a_1,a_2,a_0)})\not\subseteq \mathrm{soc}(W)$.

Let us define $X\defeq \underset{y\in\fp}{\sum}y^{p-2}\text{\tiny$\maq{1}{}{y}{}{1}{}{}{}{1}$}$. 
We note that $X$ kills $S(\tau^{I,(a_1,a_2,a_0)})$ as the latter is easily checked to be fixed under $\text{\tiny$\maq{1}{}{\fp}{}{1}{}{}{}{1}$}$; hence, to prove (i), it is enough to show that $X$ acts injectively on $\mathrm{soc}(W)^{\overline{T}(\fp),(a_2-1,a_1,a_0+1)}$.

As in Lemma~\ref{Lemma2 Nota1 Florian} we see that $W^{\overline{T}(\fp),(a_2-1,a_1,a_0+1)} = W_\mu$, where
$\mu\defeq (a_2-1,a_1,a_0+1)$. Let $\alpha\defeq (1,0,-1)\in\Phi$. By \cite{JantzenRAG}, II.1.19, (5) and (6) we have for $v_\mu\in W_{\mu}$:
\begin{equation*}
u_{\alpha}(t)v_{\mu}=\underset{i\geq 0}{\sum}t^i\left(X_{\alpha,i}v_{\mu}\right),
\end{equation*}
where $X_{\alpha,i}\in\mathrm{Dist}(\overline U_{\alpha})$ and $X_{\alpha,i}v_{\mu}\in W_{\mu+i\alpha}$
(see \cite{JantzenRAG}, II.1.12 for the definition of the operators $X_{\alpha,i}$). As $\mu+i\alpha\not\leq \lambda$ if $i\ge p$ we deduce that $X=-X_{\alpha,1}$ on $W_{\mu}$.

To prove the proposition we will show that $X_{\alpha,1}:W_{\mu}\onto W_{\mu+\alpha}$ is surjective, with a one-dimensional kernel not contained in $\mathrm{soc}(W)$.

Let $\overline{M}$ denote the Levi subgroup $\text{\tiny$\maq{\ast}{}{\ast}{}{\ast}{}{\ast}{}{\ast}$}\subseteq \overline{G}$. Any irreducible constituent of the restriction $W\vert_{\overline{M}}$ has highest weight $(\alpha,\beta,\gamma)$ with $\alpha,\beta,\gamma\in\{a_2-p+1,\dots,a_0+p-1\}$. As $a_0+p>a_0+p-1$, the Linkage Principle \cite{JantzenRAG}, II.6.17 and \cite{JantzenRAG}, II.2.12(1) yield
$$
W\vert_{\overline{M}}\cong F^{\oplus m}\oplus W',
$$
where $F$ is the irreducible $\overline{M}$-representation with highest weight $\mu$ and no irreducible constituent of $W'$ has highest weight $\mu$. 

\medskip

\noindent\textbf{Claim:}
The operator $X_{\alpha,1}$ induces an isomorphism $W'_{\mu}\stackrel{\sim}{\longrightarrow}{W}'_{\mu+\alpha}$.

\begin{proof}[Proof of the claim]
By d\'evissage we may assume that $W'$ is irreducible. Let $\nu\in X^{\ast}(\overline{T})$ be its highest weight, so in particular $\nu\neq\mu$. Then $W'\into H^{0}_{\overline{M}}(\nu)$, where $H^{0}_{\overline{M}}(\nu)$ denotes the dual Weyl module for $\overline{M}$ of highest weight $\nu$. As $\mathrm{Dist}(\overline M)$ preserves $W'\into H^{0}_{\overline{M}}(\nu)$ it is enough to show that $X_{\alpha,1}: H^0_{\overline{M}}(\nu)_{\mu}\rightarrow H^{0}_{\overline{M}}(\nu)_{\mu+\alpha}$ and $X_{-\alpha,1}: H^0_{\overline{M}}(\nu)_{\mu+\alpha}\rightarrow H^{0}_{\overline{M}}(\nu)_{\mu}$ are both injective.

Since $M^{\mathrm{der}}\cong \mathrm{SL}_2$ and $\text{\tiny$\maqdue{1}{t}{}{1}$}X^{i}Y^j=X^iY^j+tjX^{i+1}Y^{j-1}+\dots$, we see that $X_{\alpha,1}: H^0_{\overline{M}}(\nu)_{\mu}\rightarrow H^{0}_{\overline{M}}(\nu)_{\mu+\alpha}$ is injective provided that $\frac{1}{2}\left(\langle \nu,\alpha^{\vee}\rangle-\langle \mu,\alpha^{\vee}\rangle\right)\not\equiv0$ modulo $p$. 
(This is because one has $H^0_{\overline{M}}(\nu)\vert_{\overline{M}^{\mathrm{der}}}\cong \Sym^{\langle\nu,\alpha^{\vee}\rangle}(\mathrm{Std})$ and $0\neq X^{\frac 12 (\langle\nu,\alpha^\vee\rangle+\langle\mu,\alpha^\vee\rangle)}Y^{\frac 12 (\langle\nu,\alpha^\vee\rangle-\langle\mu,\alpha^\vee\rangle)}\in H^0_{\overline{M}}(\nu)_{\mu}$.)
Similarly, $X_{-\alpha,1}: H^0_{\overline{M}}(\nu)_{\mu+\alpha}\rightarrow H^{0}_{\overline{M}}(\nu)_{\mu}$ is injective provided that $\frac{1}{2}\left(\langle \nu,\alpha^{\vee}\rangle+\langle \mu+\alpha,\alpha^{\vee}\rangle\right)\not\equiv0$ modulo $p$.
As $\alpha^\vee$ is dominant and $\nu \le \lambda$ (both relative to $\overline M$) we have
$$
\langle \nu,\alpha^{\vee}\rangle\leq \langle \lambda,\alpha^{\vee}\rangle=2(p-1)-(a_2-a_0).
$$
Similarly, if $H^0_{\overline{M}}(\nu)_{\mu}\neq0$ or $H^0_{\overline{M}}(\nu)_{\mu+\alpha}\neq0$ then
$$
\langle \nu,\alpha^{\vee}\rangle>\langle \mu,\alpha^{\vee}\rangle=a_2-a_0-2
$$
(here we used that $\nu\neq\mu$). Now we can deduce the above two congruences, proving the claim.
\end{proof}

From the claim we deduce that the morphism $X_{\alpha,1}:W_{\mu}\rightarrow W_{\mu+\alpha}$ is surjective, with kernel of dimension $m$. The Kostant multiplicity formula shows that
\begin{align*}
m&=\dim_{\F} W_{\mu}-\dim_{\F} W_{\mu+\alpha}\\
&=(p+1-(a_2-a_0))-(p-(a_2-a_0))=1.
\end{align*}
Moreover, $F(\mu)\cong \mathrm{cosoc}(W)$ contains the irreducible representation of $\overline{M}$ of highest weight $\mu$.
Hence $\mathrm{soc}(W) \subseteq W'$ and $X$ is indeed injective on $\mathrm{soc}(W)_{\mu}$.
This proves part (i) of the proposition.

For statement (ii) it is now enough to apply the automorphism $g\mapsto\text{\tiny$\maq{}{}{1}{}{-1}{}{1}{}{}$}\cdot {}^tg^{-1}\cdot \text{\tiny$\maq{}{}{1}{}{-1}{}{1}{}{}$}$ to part (i) and relabel $(-a_0,-a_1,-a_2)$ as $(a_2,a_1,a_0)$.
\end{proof}

\begin{rk}
  Similarly, there exists a unique non-split extension of irreducible $\overline{G}(\fp)$-representations
$$
0\rightarrow F(a_2-1,a_1,a_0+1)\rightarrow V\rightarrow F(a_0+p-1,a_1,a_2-p+1)\rightarrow0,
$$
namely $V$ is the Weyl module of highest weight $(a_0+p-1,a_1,a_2-p+1)$.
The following operator induces an isomorphism $V^{I,(a_0,a_1,a_2)}\stackrel{\sim}{\longrightarrow}V^{I,(a_2-1,a_1,a_0+1)}$
of one-dimensional vector spaces:
\begin{equation*}
  \underset{x,y,z\in\fp}{\sum}\sum_{i=0}^{p-(a_2-a_0)} \frac{(-1)^i}{i+1} \binom{p-(a_1-a_0)}{p-(a_2-a_0)-i} (xz)^{p-1-i}y^{i+1}\text{\tiny$\maq{1}{x}{y}{}{1}{z}{}{}{1}$}
  \dot{w}_0.
\end{equation*}
We skip the somewhat lengthy proof, since we will not need it.
\end{rk}

We conclude this section with a simple but important lemma concerning the action of $U_p$-operators on $G(\Qp)$-representations over $\F$. We define
\begin{equation}\label{eq:6}
  \begin{aligned}
    U_1 &\defeq \underset{x,y\in\zp/p\zp}{\sum}\text{\tiny$\maq {p^{-1}}{}{}{}{1}{}{}{}{1}$} \text{\tiny$\maq{1}{}{}{x}{1}{}{y}{}{1}$}, \\
    U_2 &\defeq \underset{x,y\in\zp/p\zp}{\sum}\text{\tiny$\maq {p^{-1}}{}{}{}{p^{-1}}{}{}{}{1}$}
    \text{\tiny$\maq{1}{}{}{}{1}{}{x}{y}{1}$},
  \end{aligned}
\end{equation}
so that for any $G(\Qp)$-representation $\sigma$ over $\cO_E$ and any $(a_2,a_1,a_0)\in \Z^3$, both $U_1$ and $U_2$ preserve $\sigma^{I,(a_2,a_1,a_0)}$.

\begin{lm}
\label{lemma killed by U1 U2}
Let $(a_2,a_1,a_0)\in \Z^3$ be a triple with $a_2-a_1>0$, $a_1-a_0>0$ and $a_2-a_0<p-1$.  Let $\o\tau \defeq\Ind_I^K\big(\omega^{a_1}\otimes\omega^{a_2}\otimes\omega^{a_0}\big)$, and suppose that
$\sigma$ is a representation of $G(\Qp)$ over $\F$. 

Then $\Hom_K(\o\tau, \sigma)[U_i] = \Hom_K(\o\tau/M_i, \sigma)$ for $i \in \{1, 2\}$, where $M_1$ \(resp.\ $M_2$\) is the minimal subrepresentation
of $\o\tau$ containing $F(a_2,a_0,a_1-p+1)$ \(resp.\ $F(a_0+p-1,a_1,a_2-p+1)$\) as subquotient. 
Concretely, both $M_i$ are indecomposable of length 2 with socle $F(a_0+p-1,a_2,a_1)$.
\end{lm}

\begin{proof}
Note that $\o\tau^{I,(a_1,a_2,a_0)}$ is spanned by a function $f$ with $\supp(f)=I$ and $f(1)=1$. Consider the case $i = 2$, and suppose $\varphi \in \Hom_K(\o\tau, \sigma)$. Then $U_2 \varphi = 0$ iff $U_2(\varphi(f)) = 0$.
We note that $\Pi^2 U_2 = U_2'$,
where $U_2'\defeq \underset{x,y\in\Fp}{\sum}\text{\tiny$\maq {1}{x}{y}{}{1}{}{}{}{1}$} \text{\tiny$\maq{}{}{1}{1}{}{}{}{1}{}$}$, sending $\rho^{I,(a_1,a_2,a_0)}$ to $\rho^{I,(a_0,a_1,a_2)}$ for any $K$-representation $\rho$ over $\F$. 
We have that $U_2'f\neq 0$, as $(U_2'f)\text{\tiny$\maq{}{1}{}{}{}{1}{1}{}{}$}=1$.
Hence $\langle K\cdot U_2'f\rangle_\F$ is a quotient of $\o\tau' \defeq \Ind_I^K\big(\omega^{a_0}\otimes\omega^{a_1}\otimes\omega^{a_2}\big)$ that injects into $\o\tau$, so 
as the cosocle of $\o\tau'$ is $F(a_0+p-1,a_1,a_2-p+1)$ and $\o\tau$ is multiplicity-free, we deduce that $\langle K\cdot U_2'f\rangle_\F = M_2$. By 
\cite{le}, Proposition 2.2.2, $M_2$ is uniserial of shape $F(a_0+p-1,a_2,a_1)\textbf{---}F(a_0+p-1,a_1,a_2-p+1)$. 
We have that
$$
U_2\varphi=0\Longleftrightarrow U_2'(\varphi(f))=0\Longleftrightarrow \varphi(U_2'f)=0\Longleftrightarrow \langle K\cdot U_2'f\rangle_\F\subseteq \ker \varphi.
$$
The proof for the case $i = 1$ is analogous.
\end{proof}

\begin{corollary}
\label{cor killed by U1 U2}
Let $(a_2,a_1,a_0)\in \Z^3$ be a triple with $a_2-a_1>0$, $a_1-a_0>0$ and $a_2-a_0<p-1$.
Suppose $\sigma$ is a representation of $G(\Qp)$ over $\F$.
\begin{enumerate}
\item If $v\in \sigma^{I,(a_1,a_2,a_0)}$ then $U_2v=0$ if and only if $\langle Kv\rangle_{\F}$ does not contain
  $F(a_0+p-1,a_1,a_2-p+1)$ as subquotient.
\item If $v\in \sigma^{I,(a_2,a_0,a_1)}$ then $U_1v=0$ if and only if $\langle Kv\rangle_{\F}$ does not contain
  $F(a_0+p-1,a_1,a_2-p+1)$ as subquotient.
\end{enumerate}
\end{corollary}
\begin{proof}
(i) In the notation of the proof of Lemma~\ref{lemma killed by U1 U2}, by Frobenius reciprocity there is a unique $\varphi \in \Hom_K(\o\tau, \sigma)$
such that $v = \varphi(f)$. The claim follows, as $\langle Kv\rangle_{\F} = \im \varphi$ (since $\langle Kf\rangle_\F=\o\tau$).
The proof of (ii) is analogous.
\end{proof}

\subsection{Jacobi sums and characteristic zero principal series}

In this section we consider certain group algebra operators with $\oe$-coefficients and study their effect on principal series representations.

We fix once and for all a triple $(a_2,a_1,a_0)\in\Z^3$ verifying $(\ref{genericityG})$. 
We define the following elements in the group algebra $\oe[\overline G(\Fp)]$, lifting the mod-$p$ operators introduced in Section \ref{mod p group algebra operators}:
\begin{equation}
\label{operators char 0}
\begin{split}
  \widehat{S}&\defeq
  \underset{x,y,z\in\fp}{\sum}\teich{x}^{p-(a_2-a_0)}\teich{z}^{p-(a_1-a_0)}\text{\tiny$\maq{1}{{x}}{{y}}{}{1}{{z}}{}{}{1}$}\dot{w}_0,
    \\
  \widehat{S}'&\defeq
  \underset{x,y,z\in\fp}{\sum}\teich{x}^{p-(a_2-a_1)}\teich{z}^{p-(a_2-a_0)}\text{\tiny$\maq{1}{{x}}{{y}}{}{1}{{z}}{}{}{1}$}\dot{w}_0.
  \end{split}
\end{equation}

The following relation between the actions of $\widehat{S},\,\widehat{S}'$ on certain principal series representations will be crucial later.
To lighten notation, we set $\Pi\defeq\text{\tiny$\maq{}{1}{}{}{}{1}{p}{}{}$}\in N_{G(\qp)}(I)$.

\begin{prop}
\label{main char 0 GL3}
Let $(a_2,a_1,a_0)\in\Z^3$ be a triple verifying $(\ref{genericityG})$. Let $\pi_p\defeq\Ind_{B(\Qp)}^{G(\Qp)}\big(\chi_1\otimes\chi_{2}\otimes\chi_{0}\big)$
be a principal series representation, where the smooth characters $\chi_{i}:\Q_p\s\rightarrow E\s$ verify $\chi_{i}\vert_{\zp\s}=\teich{\omega}^{a_i}$ for $i\in\{0,1,2\}$. 

On the one-dimensional isotypical component $\pi_p^{I,(a_1,a_2,a_0)}$ we have
\begin{equation}
\label{equation hecke}
\widehat{S}' \circ \Pi =p\chi_1(p)\kappa\, \widehat{S},
\end{equation}
where $\kappa\in \zp\s$ verifies 
$\kappa\equiv (-1)^{a_1-a_0}\cdot\frac{a_2-a_1}{a_1-a_0}\,\,\mod{p}$ and is independent of the scalars $\chi_i(p)$.
\end{prop}
\begin{rk}
Note that equation (\ref{equation hecke}) makes sense, as both $\widehat{S}$ and $\widehat{S}'$ act on $\pi^{K(1)}$ for any smooth $K$-representation $\pi$, where $K(1)\defeq \ker\big(K\onto\overline{G}(\Fp)\big)$.
\end{rk}

The proof of Proposition \ref{main char 0 GL3} relies on certain direct manipulations on Jacobi sums and will occupy the rest of this section.
Let us pick a non-zero element $\widehat{v}$ in the one-dimensional isotypical component $\pi_p^{I,(a_1,a_2,a_0)}$. In particular, $\widehat{v}$ is $K(1)$-fixed. We note that
$\mathrm{supp}(\wh v)= B(\Qp)I$ and $\wh v(1)\ne 0$.

\begin{lm}
\label{lem1 char 0 GL3}
We have
\begin{equation*}
\Pi\cdot \widehat{v}=\chi_1(p)\,\underset{\lambda,\mu\in\Fp}{\sum}\text{\tiny$\maq{{\lambda}}{1}{}{{\mu}}{}{1}{1}{}{}$}\cdot\widehat{v}.
\end{equation*}
\end{lm}
\begin{proof}
Noting that $\Pi U_1 = \underset{\lambda,\mu\in\Fp}{\sum}\text{\tiny$\maq{{\lambda}}{1}{}{{\mu}}{}{1}{1}{}{}$}$, it suffices to show that
$U_1 \cdot \wh v = \chi_1(p)^{-1}\cdot\wh v$. As $U_1 \cdot \wh v$ lies in the one-dimensional space $\pi_p^{I,(a_1,a_2,a_0)}$
and $\widehat{v}(1) \ne 0$, it is enough to observe that
\begin{equation*}
(U_1\widehat{v})(1)=
\underset{\lambda,\mu\in \Zp/p\Zp}{\sum}\widehat{v}\left(\text{\tiny$\maq{p^{-1}}{}{}{}{1}{}{}{}{1}$}
\text{\tiny$\maq{1}{}{}{\lambda}{1}{}{\mu}{}{1}$}\right)=\chi_1(p)^{-1}\cdot\widehat{v}(1).\qedhere
\end{equation*}
\end{proof}

We now compute the action of the operator $\widehat{S}'$ on the element $\Pi\widehat{v}$. 
\begin{lm}
\label{lem2 char 0 GL3}
We have
\begin{equation*}
\hspace{-1cm}
\widehat{S}'(\Pi\widehat{v})=(-1)^{a_2-a_0}\chi_1(p)\underset{\lambda,\mu,x,y,z\in\fp}{\sum}
\teich{\lambda}^{a_2-a_1}\teich{\mu}^{p-1-(a_1-a_0)} \teich{(x-\mu)}^{p-(a_2-a_1)}\teich{(z-\lambda)}^{p-(a_2-a_0)}
\text{\tiny$\maq{1}{{x}}{{y}}{}{1}{{z}}{}{}{1}$}\dot{w}_0\widehat{v}.
\end{equation*}
\end{lm}
\begin{proof}
In order to ease notation, we write $a\defeq a_2-a_0$, $b\defeq a_1-a_0$.
An immediate computation using Lemma \ref{lem1 char 0 GL3} gives
\begin{equation}
\label{bigsum}
\widehat{S}'(\Pi\widehat{v})=\chi_1(p)\underset{\lambda,\mu,x,y,z\in\fp}{\sum}
\teich{x}^{p-(a-b)}\teich{z}^{p-a}
\underbrace{\text{\tiny$\maq{1+{y \lambda}+{x \mu}}{{y}}{{x}}{{z \lambda}+{\mu}}{{z}}{1}{{\lambda}}{1}{}$}}_{\defeq A_{\lambda,\mu}}\widehat{v}.
\end{equation}

We now break up the sum (\ref{bigsum}) according to the values of $\lambda,\mu\in\fp$.

\emph{{Case 1:} we have $\lambda\mu=0$.} 
Suppose first that $\lambda=\mu=0$. Then one has
$$
A_{0,0}=
\text{\tiny$\maq{1}{}{}{}{1}{{z}}{}{}{1}$}\dot{s}_2\text{\tiny$\maq{1}{{y}}{{x}}{}{1}{}{}{}{1}$}
$$
and therefore, using that $\widehat{v}$ is $I_1$-fixed, we obtain
\begin{align*}
\underset{x,y,z\in\fp}{\sum}\teich{x}^{p-(a-b)}\teich{z}^{p-a}
A_{0,0} \widehat{v}&=
\left(\underset{z\in\fp}{\sum}\teich{z}^{p-a}\text{\tiny$\maq{1}{}{}{}{1}{{z}}{}{}{1}$}\dot{s}_2\widehat{v}\right){\left(\underset{x,y\in\fp}{\sum}\teich{x}^{p-(a-b)}\right)}\\
&=0.
\end{align*}
In an analogous fashion, if $\lambda \mu \ne 0$ we have
\begin{align*}
A_{\lambda,0}&=
\text{\tiny$\maq{1}{}{{y}+ {\lambda}^{-1}}{}{1}{{z}}{}{}{1}$}\dot{s}_2\dot{s}_1\text{\tiny$\maq{\lambda}{1}{}{}{-{\lambda}^{-1}}{{x}}{}{}{1}$}, \\
A_{0,\mu}&=
\text{\tiny$\maq{1}{{\mu}^{-1}+{x}}{{y}-{z}({\mu}^{-1}+{x})}{}{1}{}{}{}{1}$}\dot{s}_1\dot{s}_2\text{\tiny$\maq{\mu}{{z}}{1}{}{1}{}{}{}{-\mu^{-1}}$},
\end{align*} 
so, using that $\widehat{v}$ is an eigenvector for the $I$-action, we obtain 
\begin{align*}
\underset{x,y,z\in\fp}{\sum}\,\underset{\lambda\in\fp\s}{\sum}
\teich{x}^{p-(a-b)}\teich{z}^{p-a} A_{\lambda,0}\widehat{v}&=0,
\\
\underset{x,y,z\in\fp}{\sum}\,\underset{\mu\in\fp\s}{\sum}
\teich{x}^{p-(a-b)}\teich{z}^{p-a} A_{0,\mu}\widehat{v}&=0
\end{align*}
(in the first equation, the sum over $x$ is zero; in the second equation the sum over $z$ is zero, provided
one first replaces $y$ by $y' \defeq {y}-{z}({\mu}^{-1}+{x})$).

\emph{{Case 2:} we have $\lambda\mu\neq 0$.}
As before, a direct computation gives
\begin{equation*}
A_{\lambda,\mu}=
\text{\tiny$\maq{1}{{\mu}^{-1}+{x}}{{y}+ {x\mu}{\lambda}^{-1}+{\lambda}^{-1}}{}{1}{{z}+{\mu}{\lambda}^{-1}}{}{}{1}$}\dot{w}_0\text{\tiny$\maq{\lambda}{1}{}{}{-{\mu}{\lambda}^{-1}}{1}{}{}{-{\mu}^{-1}}$}
\end{equation*}
and, recalling that $\wh v$ lies in $\pi_p^{I,(a_1,a_2,a_0)}$, we obtain
\begin{multline*}
\underset{x,y,z\in\fp}{\sum}\,\underset{\lambda,\mu\in\fp\s}{\sum}
\teich{x}^{p-(a-b)}\teich{z}^{p-a}A_{\lambda,\mu}\widehat{v}=\\
\qquad=\underset{x,y,z\in\fp}{\sum}\,\underset{\lambda,\mu\in\fp\s}{\sum}
(-1)^{a} \teich{\lambda}^{-(a-b)}\teich{\mu}^{a}\teich{x}^{p-(a-b)}\teich{z}^{p-a}
\text{\tiny$\maq{1}{{\mu^{-1}+x}}{{y}+ {x\mu}{\lambda}^{-1}+{\lambda}^{-1}}{}{1}{{z+\mu \lambda^{-1}}}{}{}{1}$}\dot{w}_0
\widehat{v}.
\end{multline*}

We obtain the desired result by putting the two cases together and using the change of variables
$(x',y',z',\lambda',\mu') \defeq (x+{\mu^{-1}, {y}+ {x\mu}{\lambda}^{-1}+\lambda}^{-1}, z+\mu \lambda^{-1}, \mu \lambda^{-1}, \mu^{-1})$.
\end{proof}

\begin{proof}[Proof of Proposition \ref{main char 0 GL3}]
Note that in the expression of Lemma \ref{lem2 char 0 GL3} we may assume that $x \ne 0$ (resp.\ $z \ne 0$)
since otherwise the sum over $\mu$ (resp.\ $\lambda)$ vanishes. Using the change of variables
$(\lambda',\mu') \defeq (\lambda z^{-1}, \mu x^{-1})$ we obtain
\begin{equation*}
\widehat{S}'(\Pi\widehat{v})=(-1)^{a}\chi_1(p)\kappa_1\kappa_2 \widehat{S}(\widehat{v}),
\end{equation*}
where $\kappa_1,\,\kappa_2\in\zp$ are the following Jacobi sums:
\begin{align*}
\kappa_1&\defeq\underset{\mu\in\fp}{\sum}\teich{\mu}^{p-1-b} \teich{(1-\mu)}^{p-(a-b)},&
\kappa_2&\defeq\underset{\lambda\in\fp}{\sum}\teich{\lambda}^{a-b} \teich{(1-\lambda)}^{p-a}.
\end{align*}
By Stickelberger's theorem (see e.g.\ \cite{BD}, Th\'eor\`eme 2.5.1) 
one has
$\ord(\kappa_1)=0,\,\ord(\kappa_2)=1$ and
\begin{align*}
\kappa_1&\equiv\frac{(p-1-b)!(p-(a-b))!}{(p-a)!},&
\frac{\kappa_2}{p}&\equiv(-1)\frac{(a-b)!(p-a)!}{(p-b)!}\mod{p}.
\end{align*}
Therefore 
\begin{equation*}
\widehat{S}'(\Pi\widehat{v})=(-1)^{a}\chi_1(p)\kappa\widehat{S}(\widehat{v}),
\end{equation*}
where $\ord(\kappa)=1$ and $\frac{\kappa}{p}\equiv -\frac{(p-(a-b))!(a-b)!}{p-b}\mod{p}$.

As $(p-(a-b))!(a-b)!\equiv (-1)^{p-(a-b)}(p-1)!(a-b)$ mod $p$ we finally obtain
$$
(-1)^{a}\frac{\kappa}{p}\equiv (-1)^{b}\frac{a-b}{b},
$$
and the proof of Proposition \ref{main char 0 GL3} is now complete.
\end{proof}

\begin{lm}
\label{lemma Up-operators}
Suppose that $\pi_p\defeq\Ind_{B(\Qp)}^{G(\Qp)}\big(\chi_1\otimes\chi_{2}\otimes\chi_{0}\big)$ is a principal series representation, where the smooth characters $\chi_{i}:\Q_p\s\rightarrow E\s$ verify $\chi_{i}\vert_{\zp\s}=\teich{\omega}^{a_i}$ for $i\in\{0,1,2\}$ and where $a_0,a_1,a_2$ are distinct modulo $p-1$.
On the one-dimensional isotypical component $\pi_p^{I,(a_1,a_2,a_0)}$ we have $U_1=\chi_1(p)^{-1}$ and $U_2=\chi_1(p)^{-1}\chi_2(p)^{-1}$.
\end{lm}
\begin{proof}
For $U_1$ this was observed in the proof of Lemma~\ref{lem1 char 0 GL3}.
The $U_2$-eigenvalue is computed similarly.
\end{proof}

\section{Local-global compatibility}
\label{sec:local-glob-comp}

In this section we establish most of our main results, namely Theorems \ref{intro:main local global}, \ref{intro:theorem weight elimination}, and 
\ref{intro:corollary weight elimination} of the introduction.

\subsection{The space of automorphic forms on certain unitary groups} 
\label{automorphic unitary groups}
Let $F/\Q$ be a CM field and $F^+$ its maximal totally real subfield. Assume that $F^+\neq \Q$, and that $p$ splits completely in $F$.
We write $c$ for the generator of $\Gal(F/F^+)$. For $w\nmid \infty$ (resp.\ $v\nmid \infty$) a place of $F$ (resp.\ $F^+$) we denote by $k_w$ (resp.\ $k_v$) the residue field of $F_w$ (resp.\ $F^+_v$).

We let $G_{/F^+}$ be a reductive group, which is an outer form of $\GL_3$ which splits over~$F$. We assume that $G(F^{+}_{v})\cong U_3(\R)$ for all $v|\infty$.  
By the argument of \cite{CHT}, \S3.3, $G$ admits a model $\cG$ over $\cO_{F^+}$ such that $\cG \times \cO_{F^+_v}$ is reductive for all places $v$ of $F^+$ that
split in $F$.
For any such place $v$ of $F^+$ and $w|v$ of $F$ we get an isomorphism $\iota_w:\,G(F^+_v)\stackrel{\sim}{\longrightarrow}\GL_3(F_w)$
which restricts moreover to an isomorphism $\iota_w:\,\cG(\cO_{F^+_v})\stackrel{\sim}{\longrightarrow}\GL_3(\cO_{F_w})$.

Let $F_p^+\defeq F^+\otimes_{\Q}\Qp$ and $\cO_{F^+,p}\defeq \cO_{F^+}\otimes_\Z\Zp$.
If $W$ is an $\cO_E$-module endowed with an action of $\cG(\cO_{F^+,p})$ and $U\leq G(\bA_{F^+}^{\infty,p})\times\cG(\cO_{F^+,p})$ is a compact open subgroup,  the space of algebraic automorphic forms on $G$ of level $U$ and coefficients in $W$ is defined as the following $\cO_E$-module:
\begin{equation*}
S(U,W)\defeq \left\{f:\,G(F^{+})\backslash G(\bA^{\infty}_{F^{+}})\rightarrow W\,|\, f(gu)=u_p^{-1}f(g)\,\,\forall g\in G(\bA^{\infty}_{F^{+}}), u\in U\right\}
\end{equation*}
(with the obvious notation $u=u^pu_p$ for the elements in $U$).

We recall that the level $U$ is said to be \emph{sufficiently small} if for all $t \in G(\bA^{\infty}_{F^+})$, the finite group $t^{-1} G(F^+) t \cap U$
is of order prime to $p$. For a finite place $v$ of $F^+$ we say that $U$ is \emph{unramified} at $v$ if one has a decomposition $U=\cG(\cO_{F_v^+})U^{v}$ for some compact open subgroup $U^v\leq G(\bA^{\infty,v}_{F^+})$. 

Let $\cP_U$ denote the set consisting of finite places $w$ of $F$ such that $v\defeq w\vert_{F^+}$ is split in~$F$, $v\nmid p$ and $U$ is unramified at $v$.
If $\cP\subseteq \cP_U$ is a subset of finite complement that is closed under complex conjugation, we write $\bT^{\cP}=\cO_E[T^{(i)}_w : w\in\cP,\,i\in\{0,1,2,3\}]$ for the abstract Hecke algebra on $\cP$, where the Hecke operator $T_w^{(i)}$ acts on the space $S(U,W)$ as the usual double coset operator
$$
\iota_w^{-1}\left[ \GL_3(\cO_{F_w}) \left(\begin{matrix}
      \varpi_{w}\mathrm{Id}_i &  \cr  & \mathrm{Id}_{3-i} \end{matrix} \right)
\GL_3(\cO_{F_w}) \right],
$$
where $\varpi_w$ denotes a uniformizer of $F_w$.

\subsection{Serre weights}\label{sec:serre weights}
In this section we recall the notion of Serre weights, as well as define the set of Serre weights of a Galois representation $\rbar: G_F\rightarrow \GL_3(\F)$. 
\begin{df}
A \emph{Serre weight for} $\cG$ (or just Serre weight if $\cG$ is clear from the context) is an isomorphism class of a smooth, (absolutely) irreducible representation of $\cG(\cO_{F^+,p})$ over $\F$.
If $w|p$ is a place of $F$, a \emph{Serre weight at $w$} is an isomorphism class of a smooth, (absolutely) irreducible representation of $\GL_3(\cO_{F_w})$ over $\F$, or equivalently an isomorphism class of an (absolutely) irreducible representation of $\GL_3(k_w)$ over $\F$.
\end{df}

Let $S_p^+$ (resp.\ $S_p$) denote the set of places of $F^+$ (resp.\ $F$) that divide $p$.
Suppose that for each $w\in S_p$ we are given a $p$-restricted triple ${a}_w=(a_{w,2},a_{w,1},a_{w,0})\in \Z^3$ such that $a_{w,i}+a_{w^c,2-i}=0$ for all $0\leq i\leq 2$ and all $w\in S_p$.
Let $F_{{a}_w}$ denote the irreducible representation $F(a_{w,2},a_{w,1},a_{w,0})$ of $\GL_3(k_w)$ over $\F$ defined in Section \ref{sec:LAS}. 
Then $F_{{a}_v}\defeq F_{{a}_w}\circ \iota_w$ is an irreducible representation of $\cG(\cO_{F^+_v})$ that is independent of the choice of place $w$ dividing $v$, and we define the Serre weight 
$F_{{a}}\defeq \underset{v\vert p}{\bigotimes} F_{{a}_v}$.
By \cite{EGH}, Lemma 7.3.4 any Serre weight is of this form.

\begin{df}
\label{definition modularity}
Let $\rbar:G_F\rightarrow \GL_3(\F)$ be a continuous, absolutely irreducible Galois representation and let $V$ be a Serre weight for $\cG$. We say that $\overline{r}$ is \emph{automorphic of weight $V$} or that $V$ is a \emph{Serre weight of $\overline{r}$} if there exist a compact open subgroup $U$ of $G(\bA^{\infty,p}_F)\times \cG(\cO_{F^+,p})$ which is unramified at all places dividing $p$ and a cofinite subset $\cP$ of $\cP_U$ such that
 $\overline{r}$ is unramified at each place of $\cP$ and $S(U,V)_{\mathfrak{m}_{\rbar}}\neq0$,
where $\mathfrak{m}_{\rbar}$ is the maximal ideal of $\bT^{\cP}$ with residue field $\F$ defined by the formula
$$
\det\left(1-\overline{r}^{\vee}(\mathrm{Frob}_w)X\right)=\sum_{j=0}^3 (-1)^j(\mathbf{N}_{F/\Q}(w))^{\binom{j}{2}}(T_w^{(j)}\bmod \mathfrak{m}_{\rbar})X^j\quad \forall w\in \cP
$$
(and $\mathbf{N}_{F/\Q}(w)$ denotes the norm from $F$ to $\Q$ of the prime $w$).
We write $W(\rbar)$ for the set of all Serre weights of $\rbar$.
\end{df}

\emph{From now on we fix a place $w\in S_p$ and assume that $\rbar$ is automorphic of weight $F_{{a}}$, where ${a}_{w'}$ lies in the lowest alcove, i.e.\ $a_{w',2}-a_{w',0}<p-2$, for all $w'\in S_p\setminus\{w,w^c\}$.}

Let $v\defeq w\vert_{F^+}$ and $V'\defeq \bigotimes_{v'\in S_p^+\setminus\{v\}}F_{{a}_{v'}}$. 
Now the construction in \cite{EGH}, \S7.1.4 (replacing the coefficients $\Zpbar$ by $\cO_E$) gives a finite
free $\cO_E$-module $W_{a_{v'}}$ with a linear action of $\cG(\cO_{F^+_{v'}})$ for any $v' \in S_p^+\setminus\{v\}$, and we let
$\wt V' \defeq \bigotimes_{v'\in S_p^+\setminus\{v\}} W_{a_{v'}}$. Note that since $a_{w'}$ 
lies in the lowest alcove for all $w'$ we have $\wt V'\otimes_{\cO_E}\F\cong V'$.

\begin{df}
We define
\begin{equation*}
W_w(\rbar)\defeq \left\{F_{{a}_w}: \text{${a}_w\in \Z^3$ is $p$-restricted and $(F_{{a}_w}\circ\iota_w)\otimes_{\F} V'\in W(\rbar)$}\right\}.
\end{equation*}
Note that by construction we have $W_w(\rbar)\neq \varnothing$.
\end{df}

\emph{From now on we suppose that $U$ is a compact open subgroup of $G(\bA^{\infty,p}_F)\times \cG(\cO_{F^+,p})$ that is sufficiently small, unramified at all places dividing $p$, and a subset $\cP \subseteq \cP_U$ as above such that $S(U,(V\circ\iota_w)\otimes_{\F} V')_{\m_{\rbar}}\neq 0$ for some Serre weight $V$ at $w$.}
(Such a group $U$ exists, since we are free to shrink $U$ further to ensure it is sufficiently small, cf.\ \cite{EGH}, Remark 7.3.6.) 

Having fixed $V'$, for any representation $V$ of $\GL_3(\cO_{F_w}) = \GL_3(\zp)$ over $\cO_E$, let
\begin{equation*}
  S(V)\defeq S(U,(V\circ \iota_w)\otimes_{\cO_E} \wt V').
\end{equation*}
In particular, $W_w(\rbar) \supseteq \{\text{Serre weights $V$ at $w$}: S(V)_{\m_{\rbar}}\neq 0\}$ (and equality holds if we shrink $U$ and $\cP$ sufficiently).

\subsection{Potentially semistable lifts}
\label{sec:potent-sst-lifts}

In this section we discuss the relationship between Serre weights and potentially semistable lifts.
We start by recalling some facts about Deligne--Lusztig representations, referring to \cite{florian-duke}, Section 4 for details.

For any $n\geq 1$ let $k_{w,n}$ denote the extension of $k_w$ of degree $n$. 
We write $T$ for a maximal torus of ${\GL_3}_{/k_w}$. We have an identification
\begin{equation}
\label{torus}
T(k_w) \stackrel{\sim}{\longrightarrow} \prod_{j=1}^rk_{w,n_j}^{\times},
\end{equation}
where $1\leq  n_j\leq 3$ and $\sum_{j=1}^r n_j=3$ (cf.\ \cite{florian-duke}, Lemma 4.7), which is well defined up to $\prod_{j=1}^r\Gal(k_{w,n_j}/k_w)$.

For any homomorphism $\theta: T(k_w)\rightarrow \qb_p^{\times}$ we have a Deligne--Lusztig representation $R^{\theta}_T$ of $\GL_3(k_w)$ over $\qpb$. Via the identification~(\ref{torus}) we have $\theta=\otimes_{j=1}^r\theta_j$, where $\theta_j:k_{w,n_j}^{\times}\rightarrow\qb_p^{\times}$. We say that $\theta$ is \emph{primitive} if 
for each $j$ the $\Gal(k_{w,n_j}/k_w)$-conjugates are pairwise distinct.
By letting $\Theta(\theta_j)$ be the cuspidal representation of $\GL_{n_j}(k_w)$ associated to the primitive character $\theta_j$ via \cite{florian-duke}, Lemma 4.7, we have
$$
R^{\theta}_{T}\cong (-1)^{3-r}\Ind_{P_{\underline{n}}(k_w)}^{\GL_3(k_w)}{\left(\otimes_{j}\Theta(\theta_j)\right)},
$$
where $P_{\underline{n}}$ is the standard parabolic subgroup containing the Levi $\prod_{j}\GL_{n_j}$.

Let $F_{w,n}\defeq W(k_{w,n})[\frac 1p]$ be the unramified extension of $F_w$ of degree $n$. We consider $\theta_j$ as a character on $\cO_{F_{w,n_{j}}}^{\times}$ by inflation and let $\Art_{F_{w,n}}$ be the isomorphism $F_{w,n}^\times\isoto W_{F_{w,n}}\ab$ of local class field theory, normalized so that geometric Frobenius elements correspond to uniformizers.
We define the inertial type $\rec(\theta) : I_{F_w} \to \qpb\s$ as follows:
\begin{enumerate}
	\item $\rec(\theta)\defeq\bigoplus_{j=1}^{3}\,\theta_j\circ\Art^{-1}_{F_{w}}$ if $\theta_j:k_w^\times\rightarrow \qb_p\s$ for all $j=1,2,3$;
	\item $\rec(\theta)\defeq\theta_1\circ\Art^{-1}_{F_{w}}\oplus \underset{\sigma\in \Gal(k_{w,2}/k_w)}{\bigoplus}\sigma\left(\theta_2\circ\Art^{-1}_{F_{w,2}}\right)$ if 
$\theta_1:k_w^\times\rightarrow \qb_p\s$ and $\theta_2: k_{w,2}^\times\rightarrow \qpb^{\times}$ are primitive characters;
	\item $\rec(\theta)\defeq\underset{\sigma\in \Gal(k_{w,3}/k_w)}{\bigoplus}\sigma\left(\theta_1\circ\Art^{-1}_{F_{w,3}}\right)$ if $\theta_1:k_{w,3}^\times\rightarrow \qb_p\s$ is a primitive character.
\end{enumerate}

\begin{lm}
\label{main global lifting lemma}
Let $V_w$ be a Serre weight at $w$ for the Galois representation $\rbar:G_F\rightarrow \GL_3(\F)$ and assume that $V_w\otimes_{\F}\Fpbar$ is a Jordan--H\"older constituent of the mod-$p$ reduction of a Deligne--Lusztig representation $R_T^{\theta}$ of $\GL_3(k_w)$, where $T$ is a maximal torus of ${\GL_3}_{/k_w}$ and $\theta: T(k_w)\rightarrow \qb_p^{\times}$ is a primitive character.
Then $\overline{r}\vert_{G_{F_w}}$ has a potentially semistable lift with Hodge--Tate weights $\{-2,-1,0\}$ and inertial type $\mathrm{rec}(\theta)$.
\end{lm}
\begin{proof}
This is proved in \cite{MP}, Theorem 5.5 (cf.\ also \cite{EGH}, Proposition 2.4.1 and Proposition 7.4.4 in cases (i), (iii)), except when
$n_i = 1$ for all $i$ and some of the characters $\theta_i$ are equal. The exceptional case follows in the same way, by using the lemma below,
which is an extension of \cite{EGH}, Proposition 2.4.1(ii). (We recall the idea: by the Deligne--Serre lifting lemma we obtain
a global lift $r$ of the dual $\o r^\vee$ which is attached to an automorphic representation in trivial weight and dual type $(R_T^{\theta})^\vee$.
Classical local-global compatibility together with the lemma below shows that $r^\vee|_{G_{F_w}}$ provides the required lift.)
\end{proof}

Suppose $n \ge 1$ and let $B_n$ denote the Borel subgroup of $\GL_n$ of upper-triangular matrices. Let $\rec_{F_w}$ denote the local Langlands correspondence 
for $\GL_n(F_w)$ over $\qpb$, as in \cite{EGH}, \S1.3.

\begin{lm}\label{lm:inertial-llc-principal-series}
  Suppose that $\pi$ is an irreducible admissible representation of $\GL_n(F_w)$ over $\qpb$.
  Assume that $\pi|_{\GL_n(\O_{F_w})}$ contains a Jordan--H\"older factor of the  $\GL_n(\O_{F_w})$-inflation of a principal series representation of the form
  \begin{equation}\label{eq:2}
    \Ind_{B_n(k_w)}^{\GL_n(k_w)} (\theta_1 \otimes \cdots \otimes \theta_n)
  \end{equation}
  for some characters  $\theta_i : k_w\s \to \qpb\s$. 
Then $\rec_{F_w}(\pi)|_{I_{F_w}} \cong \bigoplus_{i=1}^n \theta_i \circ \Art_{F_w}^{-1}$.
\end{lm}

\begin{proof}
  We indicate how the proof of \cite{EGH}, Proposition 2.4.1(ii) needs to be adjusted. 
  Recall the pair $(I,\rho)$ described in that proof, which has the property that $\Ind_I^{\GL_n(\O_{F_w})} \rho$ is isomorphic to the 
  inflation of the representation~\eqref{eq:2} to $\GL_n(\O_{F_w})$. Thus by assumption $\Hom_I(\rho, \pi) \ne 0$.
  By \cite{bib:roche}, Theorem 7.7 we deduce that $\pi$ is a subquotient of $\Ind_{B_n(F_w)}^{\GL_n(F_w)} \wt \theta$ for some character
  $\wt\theta = \wt\theta_1 \otimes \cdots \otimes \wt\theta_n$ with $\wt\theta_i|_{\O_{F_w}\s} = \theta_i$ for all $i$.
  Now from \cite{henniart-unique}, \S1.9 it follows that if $\pi'$ is another irreducible admissible representation of $\GL_n(F_w)$, then $\pi$, $\pi'$
  have the same supercuspidal support iff $\rec_{F_w}(\pi)|_{W_{F_w}} \cong \rec_{F_w}(\pi')|_{W_{F_w}}$. The claim then follows from
  \cite{harris-taylor}, \S VII.2 (pp.\ 251f).
\end{proof}

\subsection{Serre weights of $\rbar$ in the maximally non-split case}
\label{sec:Weight elimination}
We keep the notation and assumptions of Section \ref{sec:serre weights}. 
In particular, $p$ splits completely in $F$, there is a distinguished place $w\in S_p$, and $\rbar: G_F\rightarrow \GL_3(\F)$ is an automorphic, absolutely irreducible Galois representation. 
The aim of this section is to almost completely determine the Serre weights of $\rbar$ at $w$ when $\rbar\vert_{G_{F_w}}$ is maximally non-split and generic. The techniques of this section follow the strategy of \cite{EGH}, \S 5, and have been further pursued in \cite{MP} to cover other cases when $\rbar$ is ordinary at $w$.

\begin{thm}
\label{theorem weight elimination}
Assume that $\rbar|_{G_{F_w}}$ is maximally non-split and generic, of the form
\begin{equation*}
\rbar|_{G_{F_w}}\sim \maq{\omega^{a+1}\un{\mu_2}}{\ast}{\ast}{}{\omega^{b+1}\un{\mu_1}}{\ast}{}{}{\omega^{c+1}\un{\mu_0}}.
\end{equation*}
\begin{enumerate}
	\item If $\FL(\rbar\vert_{G_{F_w}})\notin\{0,\infty\}$ 
we have
\begin{align*}
&\{F(a-1,b,c+1)\}\subseteq W_{w}(\rbar)\subseteq\\
&\qquad\subseteq \{F(a-1,b,c+1),\,F(c+p-1,b,a-p+1)\}.
\end{align*}
\item If $\FL(\overline{r}\vert_{G_{F_w}})=\infty$ we have
\begin{align*}
&\{F(a-1,b,c+1), F(a,c,b-p+1)\}\subseteq W_{w}(\rbar)\subseteq\\
&\qquad\subseteq \{F(a-1,b,c+1),\,F(c+p-1,b,a-p+1), F(a,c,b-p+1)\}.
\end{align*}
\item  If $\FL(\overline{r}\vert_{G_{F_w}})=0$ we have
\begin{align*}
&\{F(a-1,b,c+1), F(b+p-1,a,c)\}\subseteq W_{w}(\rbar)\subseteq\\
&\qquad\subseteq \{F(a-1,b,c+1),\,F(c+p-1,b,a-p+1), F(b+p-1,a,c)\}.
\end{align*}
\end{enumerate}
\end{thm}

\begin{rk}
  In fact, the proof shows that the containments in (i)--(iii) hold if we replace $W_w(\rbar)$ by 
  $\{\text{Serre weights $V$ at $w$}: S(V)_{\m_{\rbar}}\neq 0\}$ (for any group $U$ and subset $\cP \subseteq \cP_U$ considered in \S\ref{sec:serre weights}).
\end{rk}

In the remainder of this section we will prove this theorem, except for showing the existence of the shadow weights $F(a,c,b-p+1)$, $F(b+p-1,a,c)$ in parts (ii) and (iii). We will complete the proof in Proposition \ref{existence of shadows}. 
We start with the following preliminary lemma.
\begin{lm}
\label{lemma weight elimination niveau 1}
Let $V\in W_{w}(\rbar)$ be a Serre weight of $\rbar$ at $w$. 
Then $V$ is isomorphic to one of the weights in the following list:
\begin{gather*}
F(b+p-1,a-1,c+1),\,F(a-1,c+1,b-p+1),\, F(c+p,b,a-p),\\
F(b+p-1,a,c),\,F(a,c,b-p+1),\\
F(c+p-1,b,a-p+1),\,F(a-1,b,c+1).
\end{gather*} 
\end{lm}
\begin{proof}
Let us write $V\cong F(x,y,z)$ for a $p$-restricted weight $(x,y,z)$.
By Frobenius reciprocity the weight $F(x,y,z)$ is a constituent of the principal series $\Ind_{\overline{B}(k_w)}^{\overline{G}(k_w)}\big(\omega^x\otimes\omega^y\otimes\omega^z\big)$.
Hence, by Lemma \ref{main global lifting lemma}, $\rbar|_{G_{F_w}}$ admits a potentially semistable lift of Hodge--Tate weights $\{-2,-1,0\}$ and inertial type $\theta\defeq\teich{\omega}^x\oplus\teich{\omega}^y\oplus\teich{\omega}^z$.
Proposition \ref{proposition weight elimination niveau 1} implies
$\{x,y,z\}\equiv \{a,b,c\}$ or $\{x,y,z\}\equiv \{a-1,b,c+1\}$ as subsets of $\Z/(p-1)\Z$.

If $x-z \ge p-2$, we conclude by genericity that $V$ belongs to the following list:
\begin{gather*}
F(b+p-1,a-1,c+1),\,F(a-1,c+1,b-p+1),\,F(c+p,b,a-p),\\
F(b+p-1,a,c),\,F(a,c,b-p+1),\,F(c+p-1,b,a-p+1).
\end{gather*}

If $x-z < p-2$, then noting that $x+1$, $y$, $z-1$ are distinct modulo $p-1$ by above, the proof
of \cite{EGH}, Lemma 6.1.1, shows that $V$ is also a constituent of the principal series 
$\Ind_{\overline{B}(k_w)}^{\overline{G}(k_w)}\big(\omega^{z-1}\otimes\omega^y\otimes\omega^{x+1}\big)$, hence
$\{x+1,y,z-1\}\equiv \{a,b,c\}$ or $\{x+1,y,z-1\}\equiv \{a-1,b,c+1\}$ as subsets of $\Z/(p-1)\Z$.
We deduce by genericity that $V\cong F(a-1,b,c+1)$.
\end{proof}

\begin{proof}[Proof of Theorem \ref{theorem weight elimination} \(modulo the existence of shadows\)]
We first establish the upper bound on $W_w(\rbar)$ (this is sometimes called ``weight elimination'').
Let $V\in W_{w}(\rbar)$. Then $V$ is one of the weights appearing in the list of Lemma \ref{lemma weight elimination niveau 1}. 
We now prove that the weights $F(a-1,c+1,b-p+1)$, $F(c+p,b,a-p)$ and $F(b+p-1,a-1,c+1)$ cannot be Serre weights of $\rbar$ at $w$.
Indeed, by the proof of \cite{florian-duke}, Proposition 7.4, it is easily checked that $F(b+p-1,a-1,c+1)\otimes_{\F}\Fpbar$ is an element of $\mathrm{JH}\big(\overline{R}_{T}^{\theta}\big)$, where the maximal torus $T$ verifies $T(k_w)\cong k_{w,2}\s\times k_{w}\s$ and 
$\theta\cong \teich{\omega}_2^{a-1+pb}\otimes\teich{\omega}^{c+1}$. 
If $F(b+p-1,a-1,c+1)$ is a Serre weight of $\rbar$ then by Lemma \ref{main global lifting lemma}, $\rbar|_{G_{F_w}}$ has a potentially semistable lift of Hodge--Tate weights $\{-2,-1,0\}$ and inertial type $\teich{\omega}_2^{a-1+pb}\oplus\teich{\omega}_2^{b+p(a-1)}\oplus\teich{\omega}^{c+1}$, which is not possible by Proposition \ref{proposition weight elimination niveau 2}.

In a similar fashion, the weights $F(a-1,c+1,b-p+1)\otimes_{\F}\Fpbar$, $F(c+p,b,a-p)\otimes_{\F}\Fpbar$ appear in $\mathrm{JH}\big(\overline{R}_{T}^{\theta}\big)$, where now $\theta\cong \teich{\omega}_2^{b+p(c+1)}\otimes\teich{\omega}^{a-1}$, providing a contradiction with Proposition \ref{proposition weight elimination niveau 2}.

Assume now that $\FL(\rbar|_{G_{F_w}})\neq \infty$; we claim that the weight $F(a,c,b-p+1)$ cannot be a Serre weight of $\rbar$ at $w$. 
Again, using \cite{florian-duke}, proof of Proposition 7.4, we see that $F(a,c,b-p+1)\otimes_{\F}\Fpbar$
is an element of $\mathrm{JH}\big(\overline{R}_{T}^{\theta}\big)$, where now $\theta\cong \teich{\omega}_2^{a+1+p(b-1)}\otimes\teich{\omega}^{c}$. By Lemma \ref{main global lifting lemma}, $\rbar|_{G_{F_w}}$ has a potentially semistable lift of Hodge--Tate weights $\{-2,-1,0\}$ and inertial type $\teich{\omega}_2^{a+1+p(b-1)}\oplus\teich{\omega}_2^{b-1+p(a+1)}\oplus\teich{\omega}^{c}$, which implies $\FL(\rbar|_{G_{F_w}})=\infty$ by Proposition \ref{weight elimination FL}.
By duality, $F(b+p-1,a,c)$ cannot be a Serre weight of $\rbar$ at $w$ if $\FL(\rbar|_{G_{F_w}})\neq 0$.

By weight cycling we will now deduce that $F(a-1,b,c+1)\in W_w(\rbar)$. 
(We refer to the introduction of \cite{EGH} for the term ``weight cycling'' as well as its history.)
Recall the (commuting) Hecke operators $\o{T}_1$, $\o{T}_2$ at $w$ defined in \cite{EGH}, \S4.2. They act on the finite-dimensional $\F$-vector space $S(V)$ for any Serre weight $V$ at $w$.
If $\o{T}_i$ has a non-zero eigenvalue on $S(F(c+p-1,b,a-p+1))_{\m_{\rbar}}$ for some $i$, we deduce exactly as in the proof of \cite{EGH}, Corollary 4.5.4 that $\rbar|_{G_{F_w}}$ has a crystalline lift $\rho$ over $E$, where $\rho$ admits a subrepresentation of Hodge--Tate weight $-(c+p+1)$ or a quotient of Hodge--Tate weight $-(a-p+1)$. This implies that $\rbar|_{I_{F_w}}$ admits $\omega^{c+2}$ as subrepresentation or $\omega^{a}$ as quotient, contradicting our assumptions on $\rbar|_{G_{F_w}}$.
Hence $\o{T}_i$ act nilpotently on $S(V)_{\m_{\rbar}}$ for $V=F(c+p-1,b,a-p+1)$. 
By a similar argument the same holds when $V=F(b+p-1,a,c)$ or $F(a,c,b-p+1)$.

Suppose that $F(a-1,b,c+1)\notin W_w(\rbar)$. 
If $V \defeq F(c+p-1,b,a-p+1)\in W_w(\rbar)$, then $S(V)_{\m_{\rbar}} \ne 0$ (for suitable $U$, $\cP$). By the previous paragraph we may apply 
weight cycling (\cite{EGH}, Proposition 6.1.3(ii)) to the smooth $\GL_3(F_w)$-representation $\pi \defeq S(U^v,V')_{\m_{\rbar}}$ (defined just after
(\ref{eq:7})) and our weight $V$.
Using the nilpotency of $\o{T}_1$ (resp.\ $\o{T}_2$) we see that $S(V'')_{\m_{\rbar}} \ne 0$, hence $V'' \in W_w(\rbar)$, for at least one weight $V''$ listed there.
However, the genericity of $\rbar|_{G_{F_w}}$ together with our upper bound for $W_w(\rbar)$ lets us rule out four out of five weights in each case and deduce that
$F(a,c,b-p+1) \in W_w(\rbar)$ (resp.\ $F(b+p-1,a,c) \in W_w(\rbar)$). But this contradicts the upper bound on $W_w(\rbar)$ we obtained above.
Similarly, if $F(b+p-1,a,c)\in  W_w(\rbar)$, then weight cycling with $\o{T}_2$ and genericity imply that $W_w(\rbar)$ contains $F(a,c,b-p+1)$, contradicting our upper bound. 
We get a similar contradiction if $F(a,c,b-p+1)\in W_w(\rbar)$.
We conclude that $F(a-1,b,c+1)\in W_w(\rbar)$.

This ends the proof of Theorem \ref{theorem weight elimination} modulo the existence of shadow weights.
\end{proof}

\begin{corollary}\label{corollary weight elimination}
In the setting of Theorem \ref{theorem weight elimination}, $\rbar$ has an automorphic lift $r:G_F\to \GL_3(\cO_E)$ \(after possibly enlarging $E$\) such that $r|_{G_{F_w}}$ is crystalline and ordinary of Hodge--Tate weights $\{-a-1,-b-1,-c-1\}$.
\end{corollary}
\begin{proof}
Keeping our notation of the proof of Theorem \ref{theorem weight elimination}, consider the Hecke operators $\o{T}_1$, $\o{T}_2$ acting on $S(F(a-1,b,c+1))_{\m_{\rbar}}$, which is non-zero by the proof above.
If $\o{T}_1$ fails to be invertible on this vector space, then one of the weights $F(c+p,a-1,b)$, $F(a-1,c+1,b-p+1)$ is contained in $W_w(\rbar)$ by \cite{EGH}, Proposition 6.1.3(i).
This contradicts the upper bound of Theorem \ref{theorem weight elimination} by genericity, hence $\o{T}_1$ is invertible on $S(F(a-1,b,c+1))_{\m_{\rbar}}$. 
The same is true for $\o{T}_2$.
Enlarging $E$ if necessary we can pick simultaneous eigenvalues $\alpha_i\in \F\s$ of $\o{T}_i$ ($i=1,2$) on this space.
Consider the double coset operators $T_w^{(i)}$ ($i=1,2$) defined as in Section \ref{automorphic unitary groups} but for our fixed place $w| p$ and let $\widetilde{T}_w^{(1)}\defeq p^{-c-1}T_w^{(1)}$, $\widetilde{T}_w^{(2)}\defeq p^{-b-c-1}T_w^{(2)}$.
Let $\T\defeq \T^{\cP}[\widetilde{T}_w^{(1)}, \widetilde{T}_w^{(2)}]$ and consider the homomorphism $\theta: \T\rightarrow \F$ whose kernel contains $\m_{\rbar}$ and which sends $\widetilde{T}_w^{(i)}$ to $\alpha_i$ ($i=1,2$).
Redefining $a_w=(a-1,b,c+1)$, the construction in \cite{EGH}, \S4.1.1 (replacing the coefficients $\Zpbar$ by $\cO_E$) gives a finite free $\cO_E$-module $W_{a_w}$ with a linear action of $\GL_3(\cO_{F_w})$. 
As $a_{w}$ lies in the lowest alcove we have $S(W_{a_w})\otimes_{\cO_E}\F\cong S(F(a-1,b,c+1))$.
The commutative $\cO_E$-algebra $\T$ acts naturally on both sides, with $\widetilde{T}_w^{(i)}$ acting as $\o{T}_i$ on the right-hand side (see \cite{EGH}, \S 4.4).
Using the Deligne--Serre lifting lemma and enlarging $E$ if necessary we can lift the Hecke eigenvalues $\theta$ to $S(W_{a_w})$ and obtain an automorphic lift $r: G_F\to \GL_3(\cO_E)$ of $\rbar$ such that $r|_{G_{F_w}}$ is crystalline of Hodge--Tate weights $\{-a-1,-b-1,-c-1\}$.
It is moreover ordinary, i.e.\ its Hodge and Newton polygons coincide, using that $\alpha_i\neq 0$ for $i=1,2$
(cf.\ the proof of \cite{EGH}, Corollary 4.5.4).
\end{proof}

\subsection{Local-global compatibility}
\label{sec:local/global compatibility}
We keep the notation and assumptions of Section \ref{sec:serre weights}. 
Suppose that $\rbar$ is as in Theorem \ref{theorem weight elimination}.
We recall the operators $S,\,S'\in \F[\GL_3(\Fp)]$ obtained from (\ref{operators char p}) by setting $(a_2,a_1,a_0)= (-c,-b,-a)$:
\begin{equation}
\begin{aligned}
S&\defeq \underset{x,y,z\in\fp}{\sum}x^{p-(a-c)}z^{p-(a-b)}\text{\tiny$\maq{1}{x}{y}{}{1}{z}{}{}{1}$}
\dot{w}_0,
\\
S'&\defeq \underset{x,y,z\in\fp}{\sum}x^{p-(b-c)}z^{p-(a-c)}\text{\tiny$\maq{1}{x}{y}{}{1}{z}{}{}{1}$}
\dot{w}_0.
\end{aligned}\label{eq:7}
\end{equation}
As the compact open subgroup $U$ we fixed is unramified at $v$, we can write $U=\cG(\cO_{F^+_v})U^v$ for some compact open subgroup $U^v\leq G(\bA^{\infty,v}_{F^+})$.
If $W$ is an $\cO_E$-module endowed with an action of $\prod_{v' \in S_p^+\setminus\{v\}} \cG(\cO_{F^+_{v'}})$ we define $S(U^v, W) \defeq \ilim_{U_v} S(U^v \cdot
U_v, W)$, where the limit runs over all compact open subgroups $U_v \le \cG(\cO_{F_v^+})$. 
Then $S(U^v,W)$ has a smooth action of $G(F_v^+)$ and hence via $\iota_w$ of $\GL_3(F_w)=\GL_3(\Qp)$.
We also define $\wt S(U^v, W) \defeq \plim_s S(U^v, W/\varpi_E^s)$, which has a linear action of $\GL_3(\Qp)$. Note that
$\wt S(U^v, W)$ is isomorphic to the $p$-adic completion of $S(U^v, W)$, provided $W$ is finitely generated.

We recall that in Section \ref{sec:LAS} we defined an Iwahori subgroup $I$ of $K = \GL_3(\Zp)$ and the notation $V^{I,(a_2,a_1,a_0)}$ for any representation $V$ of $K$ over $\cO_E$ and any triple $(a_2,a_1,a_0)\in \Z^3$.
We note that from the definitions it follows that $\wt S(U^v,\wt V')^{I,(-b,-c,-a)} \cong S(U^v,\wt V')^{I,(-b,-c,-a)}$,
a finite free $\cO_E$-module. Recall furthermore that $\Pi = \text{\tiny$\maq{}{1}{}{}{}{1}{p}{}{}$}$.

\begin{thm}
\label{main local global}
We make the following assumptions:
\begin{enumerate}
\item $\rbar|_{G_{F_w}}$ is maximally non-split and generic, of the form
  \begin{equation*}
    \rbar|_{G_{F_w}}\sim \maq{\omega^{a+1}\un{\mu_2}}{\ast}{\ast}{}{\omega^{b+1}\un{\mu_1}}{\ast}{}{}{\omega^{c+1}\un{\mu_0}}.
  \end{equation*}
\item $\FL(\rbar|_{G_{F_w}})\not\in\{0,\infty\}$.
\item The $\cO_{E}$-dual of $S(U^v,\wt V')_{\m_{\rbar}}^{I,(-b,-c,-a)}$ is free over $\T$, where 
  $\T$ denotes the $\cO_{E}$-subalgebra of $\End\big(S(U^v,\wt V')_{\m_{\rbar}}^{I,(-b,-c,-a)}\big)$ generated by $\bT^{\cP}$, $U_1$, and $U_2$.
\end{enumerate}
Then we have the equality
\begin{equation}
{S}'\circ \Pi =(-1)^{a-b}\cdot\frac{b-c}{a-b}\cdot\FL(\rbar|_{G_{F_w}})\cdot {S}\label{eq:main local global}
\end{equation}
of maps
\[ S(U^v,V')[\m_{\rbar}]^{I,(-b,-c,-a)}[U_1,U_2] \to S(U^v,V')[\m_{\rbar}]^{I,(-c-1,-b,-a+1)}.  \]
Moreover, these maps are injective with non-zero domain.
In particular, $\FL(\rbar|_{G_{F_w}})$ is determined by the smooth $\GL_3(\Qp)$-representation $S(U^v,V')[\m_{\rbar}]$.
\end{thm}
\begin{rk}
The proof shows that for any non-zero $v \in S(U^v,V')[\m_{\rbar}]^{I,(-b,-c,-a)}[U_1,U_2]$, the
$K$-subrepresentation generated by $v$ and $\Pi v$ is of the form 
$$
F(-c-1,-b,-a+1)\textbf{---}\big(F(-b+p-1,-c,-a)\oplus F(-c,-a,-b-p+1)\big)
$$
with socle $F(-c-1,-b,-a+1)$. Together with Proposition~\ref{main mod p GL3} this explains why there exists
a constant $x \in \F\s$ such that $S'\circ \Pi = x S$ on $S(U^v,V')[\m_{\rbar}]^{I,(-b,-c,-a)}[U_1,U_2]$.
\end{rk}
\begin{proof}
Let $\eta:I\to\oe\s$ denote the character $\teich{\omega}^b\otimes\teich{\omega}^c\otimes\teich{\omega}^a$ and set
$$
M\defeq S(U^v,\wt V')^{I,(-b,-c,-a)}_{\m _{\rbar}}\cong 
S\big(U^v\cdot I,\eta\otimes\wt V '\big)_{\m_{\rbar}}.
$$
We remark that $M\neq 0$, as $S(F(a-1,b,c+1))_{\m_{\rbar}}\ne 0$ and $F(a-1,b,c+1)$ is a Jordan--H\"older factor of $\o{\Ind_I^K \eta}$, and also that $M$ is a finite free $\cO_E$-module. For any $\cO_E$-algebra $A$ we let $M_A \defeq M\otimes_{\cO_E} A$ and $\T_A \defeq \T\otimes_{\cO_E} A$.

Picking a $\qp$-linear embedding $E \into \qpb$, as well as an isomorphism $\imath : \qpb \congto \C$, we see that
\begin{equation}
\label{decomposition}
M_{\qpb} \cong\underset{\pi}{\bigoplus}m(\pi)\cdot \pi_v^{I,(-b,-c,-a)}\otimes\big(\pi^{\infty,v}\big)^{U^v},
\end{equation}
where the sum runs over irreducible representations $\pi \cong \pi_\infty \otimes \pi_v \otimes \pi^{\infty,v}$ of $G(\A_{F^+})$ such that
$\pi \otimes_\imath \C$ is a cuspidal automorphic representation of multiplicity $m(\pi) \in \Z_{>0}$ with $\pi_\infty \otimes_\imath \C$ 
determined by the algebraic representation $(\wt V ')^{\vee}$ and with Galois representation $r_{\pi}$ lifting $\rbar^{\vee}$
(cf.\ \cite{EGH}, Lemma 7.1.6).
For any $\pi$ contributing to \eqref{decomposition} we have 
\begin{itemize}
\item[(a)] $\pi_v\cong \Ind_{B(\Qp)}^{G(\Qp)}\big(\psi_b\norm{\cdot}^2\otimes\psi_c\norm{\cdot}\otimes\psi_a\big)$ for some smooth characters $\psi_i:\Qp\s\to\Qpbar\s$ (depending on $\pi$) with $\psi_i|_{\Zp\s}=\teich{\omega}^{-i}$ for $i\in\{a,b,c\}$;
\item[(b)] $r_{\pi}^\vee|_{G_{F_w}}$ is potentially crystalline with Hodge--Tate weights $\{-2,-1,0\}$ and $\WD(r_{\pi}^\vee|_{G_{F_w}})^{\mathrm{F-ss}}\cong \psi_b^{-1}\oplus\psi_c^{-1}\oplus\psi_a^{-1}$.
\end{itemize}
Here, part (a) follows from \cite{EGH}, Propositions 2.4.1 and 7.4.4 and (b) follows from classical local-global compatibility.
In particular from (b) we deduce that the $\varphi$-eigenvalue on 
$\Dst^{\Qp,2}(r_{\pi}^\vee|_{G_{F_w}})^{I_{F_w}=\teich{\omega}^b}$ equals $p^2\psi_b(p)^{-1}$, so from Theorem \ref{main theorem Galois} and assumptions (i), (ii) it follows that
\begin{equation}
\label{local global 2}
\ord(\psi_b(p))=1\quad \text{and}\quad \FL(\rbar|_{G_{F_w}})=\left[\frac{\psi_b(p)}{p}\right].
\end{equation}
We also deduce from Corollary~\ref{corollary Shape filtration str div modules} that $\ord(p^2\psi_i(p)^{-1}) \in (0,2)$ for $i \in \{a,b,c\}$,
with sum equal to 3.
In particular, the eigenvalue of $U_1$ (resp.\ $U_2$) on $\pi_v^{I,(-b,-c,-a)}$, which equals $p^2 \psi_b(p)^{-1}$
(resp.\ $p^3 \psi_b(p)^{-1} \psi_c(p)^{-1}$) by Lemma~\ref{lemma Up-operators},
has positive valuation. (Note for later reference that this is true even if $\FL(\rbar|_{G_{F_w}})\in\{0,\infty\}$.)

Note that the image $\T_a$ of $\T^\cP$ in $\End(M)$ is local with maximal ideal the image of $\m_{\rbar}$.
As $\T_a \subseteq \T$ is a finite ring extension, $\m_{\rbar} \subseteq \rad(\T)$. Since all eigenvalues of $U_i$ on $M_{\qpb}$ have positive
valuation, we deduce that $U_i \in \rad(\T)$ ($i = 1$, $2$). 
It then follows that $\T$ is local with maximal ideal $\m$ generated by $\m_{\rbar}$, $U_1$, and $U_2$.
As $\T_{\qpb}$ acts faithfully and semisimply on $M_{\qpb}$ (by \eqref{decomposition}), we see that
$\T_{\qpb}$ is semisimple, hence reduced. In particular, $\T$ is reduced.

As $U$ is sufficiently small, our assumptions as well as the statement we want to prove is insensitive to a finite base change $E \to E'$.
Hence, by passing to a finite extension of $E$, we may assume that $\T_E \cong E^r$ for some $r > 0$.

We have $M_E=\bigoplus_{\p} M_E[\p_E]$, where the sum runs over the minimal primes of $\T$ and $\p_E\defeq\p \T_E$.  Note
that, by the above, $\T_E/\p_E = E$ for any such $\p$. Then $M_E[\p_E]\otimes_E \qpb$ is a direct summand of
\eqref{decomposition}, where $\pi$ runs over a subset of the automorphic representations in \eqref{decomposition},
and we claim that each $\pi_v$ occurring in this direct summand is the same. 
To see this, note that $r_\pi$ is determined by $\p$, by using Cebotarev density and classical local-global compatibility at
the places in ${\cP}$. (Note that classical local-global compatibility, which is known only up to Frobenius-semisimplification, 
determines $\tr(r_\pi|_{G_{F_{w'}}})$ for $w' \in \cP$ and
that $\tr(r_\pi)$ determines $r_\pi$.) Then we deduce the claim from classical local-global compatibility at $w$.

By Proposition \ref{main char 0 GL3} we have 
\begin{equation}
\label{key equation char 0}
\widehat{S}'\circ \Pi =\frac{\psi_b(p)}{p}\kappa\, \widehat{S}
\end{equation}
on $M_E[\p_E]$, where $\kappa\in \zp\s$ is such that $\kappa\equiv(-1)^{a-b}\cdot\frac{b-c}{a-b}\,\mod{p}$ and where $\widehat{S},\,\widehat{S}'\in \oe[\GL_3(\Fp)]$ are obtained from (\ref{operators char 0}) by setting $(a_2,a_1,a_0)= (-c,-b,-a)$.

Let $M^{\mathrm d} \defeq \Hom_{\cO_E}(M,\cO_E)$, which is finite free over $\T$ by assumption (iii). 
Fix any minimal prime $\p$ of $\T$. As $\T_E/\p_E = E$ we have $\T/\p=\cO_E$, so in particular $M^{\mathrm d}/\p$ is a finite free $\cO_E$-module.
The identity \eqref{key equation char 0} holds also on $M_E^{\mathrm d}/\p_E\cong \Hom_{E}(M_E[\p_E],E)$ and on its submodule $M^{\mathrm d}/\p$.
By (\ref{key equation char 0}) and \eqref{local global 2}, and since $\p+\varpi_E\T=\m$, we deduce the identity \eqref{eq:main local global}
on $(M^{\mathrm d}/\p) \otimes_{\cO_E} \F \cong M^{\mathrm d}/\m$, hence also on $M_{\F}[\m] \cong \Hom_{\F}(M^{\mathrm d}/\m,\F)$. Now observe that
$M_{\F}\cong S(U^v,V')^{I,(-b,-c,-a)}_{\m_{\rbar}}$  and hence $M_{\F}[\m] \cong S(U^v,V')[\m_{\rbar}]^{I,(-b,-c,-a)}[U_1,U_2]$.

We claim that $U_i \Pi (M_{\F}[\m]) = 0$ for $i = 1$, $2$. Let $N \defeq S(U^v,\wt V')^{I,(-c,-a,-b)}_{\m _{\rbar}}$ and let $\T'
\subseteq \End(N)$ be the $\cO_E$-subalgebra generated by $\T^{\cP}$, $U_1$, and $U_2$. By above, $\Pi (M[\p]) \subseteq N[\p']$
for the minimal prime $\p'$ of $\T'$ lying over the same prime of $\T^{\cP}$ as $\p$ such that moreover the image of $U_i$ in
$\T'/\p' = \cO_E$ equals the eigenvalue of $U_i$ on $\pi_v^{I,(-c,-a,-b)}$. (In fact, $\Pi (M[\p]) = N[\p']$, by reversing the
argument.) From the above, the natural injection $M[\p] \otimes_{\cO_E} \F \into M_{\F}[\m]$ is an
isomorphism, so $\Pi (M_{\F}[\m])$ is killed by the unique maximal ideal $\m'$ of $\T'$, which contains $U_1$, $U_2$ by the same
argument as for $\T$. This proves the claim. (In fact, if $\Hom_{\cO_E}(N,\cO_E)$ is moreover free as $\T'$-module, which is true
under the hypotheses of Theorem~\ref{freeT}, then we even get that $\Pi$ induces an isomorphism $M_{\F}[\m] \congto N_{\F}[\m']$.)

Now by Theorem \ref{theorem weight elimination}(i) and \cite{EGH}, Lemma 7.4.3 any irreducible $K$-subrepresentation of $S(U^v,V')_{\m_{\rbar}}$ is isomorphic
to $F(-c-1,-b,-a+1)$ or $F(-a+p-1,-b,-c-p+1)$. Thus by Corollary \ref{cor killed by U1 U2}(i),
if $v\in M_{\F}[\m]$ is non-zero, then $\langle K v\rangle_\F$ is the unique quotient of 
$\Ind_I^K\big(\omega^{-b}\otimes\omega^{-c}\otimes\omega^{-a}\big)$ with socle $F(-c-1,-b,-a+1)$. By \cite{le},
Proposition 2.2.2 we see that $\langle K v\rangle_\F$ is the uniserial length 2 representation of shape $F(-c-1,-b,-a+1)\textbf{---}F(-b+p-1,-c,-a)$.
A similar argument shows that $\langle K \Pi v\rangle_\F$ is the uniserial length 2 representation of shape 
$F(-c-1,-b,-a+1)\textbf{---}F(-c,-a,-b-p+1)$. In particular, $S v\ne 0$ by Proposition \ref{main mod p GL3}.
\end{proof}

\begin{rk}\label{rk:main-local-global1}
  Theorem \ref{main local global} holds equally well if we replace $\T$ in the statement and proof by the subalgebra
  $\T_2 \subseteq \T$ generated by $\T^{\cP}$ and $U_2$. Assumption (iii) will be verified in Section~\ref{sec:freeness-over-hecke}
  both for $\T_2$ and for $\T$, under suitable hypotheses. The same comments apply for the subalgebra $\T_a\subseteq \T$, provided
  $W_w(\rbar) = \{F(a-1,b,c+1)\}$ (see Remark~\ref{rk:no-shadow}).
\end{rk}

\begin{rk}\label{rk:main-local-global2}
  Theorem \ref{main local global} holds if we replace assumption (iii) by either of the following two statements:
  \begin{equation}
    \label{eq:4}\tag{iii$'$}
    \dim_{\F} S(U^v,V')[\m_{\rbar}]^{I,(-b,-c,-a)}[U_1,U_2] = 1;
  \end{equation}
  \begin{equation}
    \label{eq:5}\tag{iii$''$}
    \dim_{\F} S(U^v,V')[\m_{\rbar}]^{I,(-b,-c,-a)} = 1.
  \end{equation}
  It suffices to note that \eqref{eq:5} $\Rightarrow$ \eqref{eq:4} $\Rightarrow$ (iii). The first implication is obvious. For the second, we have
  $\dim_{\F} M^{\mathrm d}/\m = \dim_{\F} M_{\F}[\m] = 1$, in the notation of the proof. As in the proof, we may assume that $\T/\m = \F$.
  By Nakayama's lemma, we have a surjective
  $\T$-linear map $\T \onto M^{\mathrm d}$, which has to be an isomorphism, as $M^{\mathrm d}$ is a faithful $\T$-module.
  (Note that assumption \eqref{eq:5} even implies that $\T_a = \T$.)
\end{rk}

We now establish the existence of shadow weights in Theorem \ref{theorem weight elimination}, thereby completing the proof of that theorem.

\begin{prop}
\label{existence of shadows}
Suppose $\rbar$ satisfies assumption \(i\) of Theorem \ref{main local global}.
\begin{enumerate}
\item If $\FL(\overline{r}\vert_{G_{F_w}})=\infty$ then
$F(a,c,b-p+1)\in W_{w}(\rbar)$.
\item  If $\FL(\overline{r}\vert_{G_{F_w}})=0$ then
$F(b+p-1,a,c)\in W_{w}(\rbar)$.
\end{enumerate}
\end{prop}

\begin{rk}
Suppose that $\rbar$ is as in Proposition~\ref{existence of shadows}(i).
Then the proof below shows that for some non-zero $v \in S(U^v,V')[\m_{\rbar}]^{I,(-b,-c,-a)}[U_1,U_2]$, the
$K$-subrepresentation generated by $v$ and $\Pi v$ is of the form 
$$F(-b+p-1,-c,-a)\oplus \big(F(-c-1,-b,-a+1)\textbf{---}F(-c,-a,-b-p+1)\big).$$
In particular, $Sv = 0$, $S' \Pi v \ne 0$.
Similar remarks apply in case $\FL(\overline{r}\vert_{G_{F_w}})=0$.
\end{rk}

\begin{proof}
(i) We follow the proof of Theorem \ref{main local global}, noting that with our current assumption only three things change:
(a) we get $\ord(\psi_b(p))<1$ instead of \eqref{local global 2},
(b) the natural injection $M[\p] \otimes_{\cO_E} \F \into M_{\F}[\m]$ need not be an isomorphism,
and (c) we only know that $\langle K v\rangle_\F$ (for non-zero $v\in M_{\F}[\m]$) is a quotient
of the uniserial representation of shape $F(-c-1,-b,-a+1)\textbf{---}F(-b+p-1,-c,-a)$.
From (a) we deduce that $\widehat{S}=\frac{p}{\psi_b(p)}\kappa^{-1}\,\widehat{S}'\circ \Pi$ on $M_E[\p_E]$, so $Sv=0$ for some non-zero $v \in M_\F[\m]$
(but maybe not all, see (b)). 
Then (c) together with Proposition~\ref{main mod p GL3} implies that $\langle K v\rangle_\F \cong F(-b+p-1,-c,-a)$ for such $v$.
In particular, $F(a,c,b-p+1)\in W_{w}(\rbar)$. Part (ii) is analogous.
\end{proof}

\section{Freeness over the Hecke algebra}
\label{sec:freeness-over-hecke}

In this section, we prove Theorem \ref{freeT}, which states that the dual 
\[\Hom_{\cO_E}(S(U^v,\widetilde{V}')_{\m_{\rbar}}^{I,(-b,-c,-a)}, \cO_E)\]
of the space of automorphic forms is free over a Hecke algebra under certain conditions on $U^v$ and $\rbar$
(and $\widetilde{V}'$ and $\m_{\rbar}$ are as defined in Section \ref{sec:serre weights}).

\subsection{The setup}
\label{sec:setup}

As before, let $F/\Q$ be a CM field in which $p$ splits completely, and let $F^+$ be its maximal totally real subfield. Assume moreover:
\begin{enumerate}
\item \label{field unram} $F/F^+$ is unramified at all finite places.
\end{enumerate}
Fix a place $w|p$ of $F$, and let $v \defeq w|_{F^+}$.
Let $\rbar:G_F \rightarrow \GL_3(\F)$ be a Galois representation with $\rbar|_{G_{F_w}}$ maximally non-split and generic as in Theorem~\ref{main local global}(i),
satisfying the following additional properties:
\begin{enumerate}\setcounter{enumi}{1}
\item \label{rbar-unram} $\rbar$ is unramified at all finite places not dividing $p$;
\item \label{rbar-FL} $\rbar$ is Fontaine--Laffaille and regular at all places dividing $p$;
\item \label{adequate} $\rbar$ has image containing $\GL_3(k)$ for some $k \subseteq \F$ with $\#k>9$;
\item \label{zeta-p} $\overline{F}^{\ker \ad \rbar}$ does not contain $F(\zeta_p)$.
\end{enumerate}

Note that condition \ref{adequate}, which is stronger than the usual condition of adequacy (see Definition 2.3 of \cite{Thorne}), allows us to choose a finite place $v_1$ of $F^+$ which is prime to $p$ satisfying the following properties (see Section 2.3 of \cite{CEGGPS}):
\begin{itemize}
 \item $v_1$ splits in $F$ as $v_1 = w_1 w_1^c$;
 \item $v_1$ does not split completely in $F(\zeta_p)$;
 \item $\rbar(\Frob_{w_1})$ has distinct $\F$-rational eigenvalues, no two of which have ratio $(\mathbf{N}v_1)^{\pm 1}$.
\end{itemize}

We choose a unitary group $G_{/F^+}$ and a model $\cG_{/\cO_{F^+}}$ as in \S\ref{automorphic unitary groups}.
We note that $G$ is automatically quasi-split, hence unramified by \ref{field unram}, at all finite places, as we
are in odd rank. (See the proof of \cite{BCh-book}, Lemma 6.2.4.)
Let $U^v = \underset{v'\neq v}{\prod} U_{v'} \leq G(\bA^{\infty,v}_{F^+})$ be a compact open subgroup satisfying the following properties:
\begin{enumerate}\setcounter{enumi}{5}
\item $U_{v'} = \cG(\cO_{F^+_{v'}})$ for all places $v'$ which split in $F$ other than $v_1$ and those dividing $p$;
  \item  $U_{v_1}$ is the preimage of the upper-triangular matrices under the map
\[\cG(\cO_{F^+_{v_1}})\to \cG(k_{v_1}) \underset{\iota_{w_1}}{\stackrel{\sim}{\longrightarrow}} \GL_3(k_{w_1});\]
 \item $U_{v'}$ is a hyperspecial maximal compact open subgroup of $G(F^+_{v'})$
    if $v'$ is inert in $F$.
\end{enumerate}
The choice of $U_{v_1}$ implies that $U^v U_v$ is sufficiently small in the sense of Section \ref{automorphic unitary groups} for any compact
open subgroup $U_v$ of $G(F^+_v)$.

As in Section~\ref{sec:serre weights} suppose that for each $w'\in S_p\setminus\{w, w^c\}$ we are given a $p$-restricted triple ${a}_{w'}=(a_{w',2},a_{w',1},a_{w',0})\in \Z^3$ such that $a_{w',i}+a_{{w'}^c,2-i}=0$ for all $0\leq i\leq 2$ and all $w'\in S_p\setminus\{w, w^c\}$. 
Suppose moreover that $a_{w',2}-a_{w',0}<p-3$, for all $w'\in S_p\setminus\{w,w^c\}$, a slightly stronger condition than before.
Recall that we defined $V' = \bigotimes_{v'\in S_p^+\setminus\{v\}}F_{{a}_{v'}}$ and $\wt V' = \bigotimes_{v'\in S_p^+\setminus\{v\}} W_{a_{v'}}$.

Let $\cP$ denote the set of finite places $w'$ of $F$ that split over~$F^+$ and do not divide $p$ or $v_1$,
and define the maximal ideal $\m_{\rbar}$ of $\T^{\cP}$ as in Section~\ref{sec:serre weights}.
We make the additional automorphy assumption:
\begin{enumerate}\setcounter{enumi}{8}
\item \label{automorphy} $S(U^v,{V}')_{\m_{\rbar}}$ is non-zero.
\end{enumerate}
By Theorem \ref{theorem weight elimination} this implies $S(U^v \cG(\cO_{F^+_{v}}),{V}'\otimes F(a-1,b,c+1))_{\m_{\rbar}} \ne 0$ and hence
also assumptions \ref{rbar-unram}--\ref{rbar-FL} above.
Just as in the proof of Theorem~\ref{main local global} it also implies that $S(U^v,\widetilde{V}')_{\m_{\rbar}}^{I,(-b,-c,-a)} \ne 0$.

Let $\bT_a$ (resp.\ $\bT$) denote the $\cO_E$-subalgebra of \[\End(S(U^v,\widetilde{V}')_{\m_{\rbar}}^{I,(-b,-c,-a)}) \cong \End((S(U^v,\widetilde{V}')_{\m_{\rbar}}^{I,(-b,-c,-a)})^\mathrm{d})\] generated by $\bT^{\cP}$ (resp.\ $\bT^{\cP}$, $U_1$, and $U_2$).
Here the subscript ``$a$'' stands for the ``anemic" Hecke algebra, and the superscript ``$\mathrm d$'' denotes the Schikhof dual (see Section 1.8 of \cite{CEGGPS}).
Note that $\bT_a$ and $\bT$ are local $\cO_E$-algebras (see the proof of Theorem \ref{main local global}).

\begin{thm} \label{freeT}
Let $\rbar$ be as in Theorem \ref{main local global}\(i\) with $\FL(\rbar\vert_{G_{F_w}}) \ne \infty$.
Assume \ref{field unram}--\ref{automorphy} in the setup above.
Then \[(S(U^v,\widetilde{V}')_{\m_{\rbar}}^{I,(-b,-c,-a)})^\mathrm{d}\]
is free over $\bT$.
\end{thm}

Because of the following lemma, we may assume in the proof of Theorem~\ref{freeT}
that $E$ is as large as we like. (In fact, for the remainder of this section it will suffice to assume that $E$ is large
enough such that $\T_a[1/p] \cong E^r$ for some $r$ and that $\T/\m_{\T} = \F$.)

\begin{lm}
  Suppose that $M$ is a finite free $\cO_E$-module and that $A$ is a local $\cO_E$-algebra acting faithfully on $M$ with $A/\m_A = \F$. 
  If $E'/E$ is a finite extension and $M \otimes_{\cO_E} \cO_{E'}$ is free over $A \otimes_{\cO_E} \cO_{E'}$, then $M$ is free over $A$.
\end{lm}

\begin{proof}
  As $A/\m_A = \F$ we deduce that $A \otimes_{\cO_E} \cO_{E'}$ is local with residue field $\cO_{E'}/\m_{E'}$. Pick any surjection $f : A^d \onto M$
  with $d$ minimal. It is easy to see that $f \otimes_{\cO_E} \cO_{E'}$ is an isomorphism, hence so is $f$.
\end{proof}

\subsection{The Taylor--Wiles method}

Let $S$ be the set of places $S_p^+ \cup \{v_1\}$.
For each place $v'$ in $S$, fix a place $\widetilde{v}'$ of $F$ lying over $v'$ and let $\widetilde{S}$ be the set of these places $\widetilde{v}'$.
We will assume that $\wt v = w$. For places $w'$ in $F$, let $R_{w'}^\square$ be the universal $\mathcal{O}_E$-lifting ring of $\rbar|^\vee_{G_{F_{w'}}}$.
For $w'\in S_p\setminus\{w, w^c\}$ let $\psi_{w'} = a_{w'} \in \Z^3$ and
let $R^{\square,\psi_{w'}}_{w'}$ be the framed crystalline deformation ring for $\rbar|^\vee_{G_{F_{w'}}}$ of Hodge--Tate weights $\psi_{w'}+(2,1,0)$.
Let $\delta_{F/F^+}$ denote the quadratic character of $F/F^+$.
Consider the deformation problem
\[\cS \defeq \(F/F^+,S,\widetilde{S},\cO_E,\rbar^\vee, \varepsilon^{-2} \delta_{F/F^+}, \{ R^\square_{\widetilde{v}_1} \} \cup \{ R^{\square, \psi_{\widetilde{v}'}}_{\widetilde{v}'}\}_{v'|p, v'\neq v} \cup \{ R^\square_{\widetilde{v}} \} \)\]
in the terminology of \cite[\S 2.3]{CHT}.
There is a universal deformation ring $R_\cS^{\textrm{univ}}$ and a universal $S$-framed deformation ring $R_\cS^{\square_S}$ in the sense of \cite[\S 2.2]{CHT}.
(We work with deformations of $\rbar^\vee$ to be consistent with \cite{CEGGPS}. Note that $\rbar^\vee(\Frob_w)$ is used in our definition of $\m_{\rbar}$
in \S \ref{sec:serre weights}.)

Let 
\[R^{\textrm{loc}} \defeq \widehat{\otimes}_{v'|p, v'\neq v} R^{\square,\psi_{\widetilde{v}'}}_{\widetilde{v}'} \widehat{\otimes} R^\square_{\widetilde{v}} \widehat{\otimes} R^\square_{\widetilde{v}_1},\]
where all completed tensor products are taken over $\cO_E$.
Choose an integer $q\geq 3[F^+:\Q]$ as in Section 2.6 of \cite{CEGGPS}.
We introduce the local ring \[R_\infty \defeq R^{\textrm{loc}}[[x_1,\ldots ,x_{q-3[F^+:\Q]}]],\]
over which \cite{CEGGPS} constructs a patched module of automorphic forms.
Let $\tau \defeq \Ind_I^K (\teich{\omega}^{-b}\otimes \teich{\omega}^{-c} \otimes \teich{\omega}^{-a}) \otimes_{\Z_p} \cO_E$ be the natural $\cO_E$-lattice in a principal series type over $E$.
In the following we will identify $G(F^+_v)$ with $\GL_3(\qp)$ and $\cG(\cO_{F^+_{v}})$ with $K=\GL_3(\zp)$ via $\iota_w$.
We collect some results from \cite{CEGGPS}.

\begin{thm} \label{patchedmod}
There exist a map
\[S_\infty \defeq \cO_E[[z_1,\ldots , z_{9\# S},y_1,\ldots, y_q]] \rightarrow R_\infty\]
and an $R_\infty[[K]]$-module $M_\infty$ together with a compatible $\GL_3(\qp)$-action, satisfying the following properties.
For a finitely generated $\cO_E$-module $W$ with continuous $K$-action, let $M_\infty(W)$ denote 
$\Hom\cont_{K}(W,M_\infty^\vee)^\vee$,
where $\cdot^\vee$ denotes the Pontryagin dual.
\begin{enumerate}
\item \label{proj} $M_\infty$ is a finitely generated projective $S_\infty[[K]]$-module. In particular,
  if $W$ is $p$-torsion free, then $M_\infty(W)$ is a finite free $S_\infty$-module.
\item \label{cm} Suppose that $W[1/p]$ is a locally algebraic type \(as defined in \cite{CEGGPS}, \S4\).
Assume that $M_\infty(W)$ is non-zero and let $\Spec R_\infty(W)$ be its support in $\Spec R_\infty$.
Then $M_\infty(W)$ is a maximal Cohen--Macaulay $R_\infty(W)$-module and $M_\infty(W)[1/p]$ is a projective $R_\infty(W)[1/p]$-module.
\item \label{special} Let $\a \defeq (z_1,\ldots, z_{9\# S}, y_1, \ldots, y_q)$ be the augmentation ideal of $S_\infty$.
There is a natural $\GL_3(\qp)$-equivariant identification 
\[(M_\infty/\a)^\mathrm{d} \cong \wt S(U^v,\widetilde{V}')_{\m_{\rbar}},\]
which induces a $\GL_3(\qp)$-equivariant identification 
\[(M_\infty(\tau)/\a)^\mathrm{d} \cong S(U^v,\widetilde{V}')_{\m_{\rbar}}^{I,(-b,-c,-a)}.\]
Furthermore, there is a surjection $R_\infty \onto \bT_a$ so that the latter identification is $R_\infty$-equivariant.
\item \label{patchedfiber} Let $\wp_r$ denote the prime ideal of $R_\infty$ corresponding to an $\cO_E$-point of $\Spec \bT_a$ and by abuse, the corresponding prime ideal of $\bT^{\cP}$. By abuse, let $\m_{\rbar}$ also denote the maximal ideal of $R_\infty$.
Then
\[(M_\infty(\tau)/\wp_r)^\mathrm{d} \cong S(U^v,\widetilde{V}')[\wp_r]^{I,(-b,-c,-a)}\]
and
\[(M_\infty(\tau)/\m_{\rbar})^\vee \cong S(U^v,V')[\m_{\rbar}]^{I,(-b,-c,-a)}.\]
\end{enumerate}
\end{thm}

Note that if $W$ is $p$-torsion free, we have $M_\infty(W) \cong \Hom\cont_{K}(W,M_\infty^\mathrm{d})^\mathrm{d}$ by Remark 4.15 of \cite{CEGGPS}.

\begin{proof}The construction of $M_\infty$ is in Section 2.8 of \cite{CEGGPS}, except that we allow the Hodge--Tate weights $\psi_{\widetilde{v}'}$ to depend on $v'$, and we do not include Hecke operators at $v_1$ in our Hecke algebras.
Without these Hecke operators, the rank of $M_\infty(W)[1/p]$ as a projective $R_\infty(W)[1/p]$-module in part \ref{cm} is greater than one. Otherwise, the necessary modifications are minor and straightforward.
The map
\[S_\infty = \cO_E[[z_1,\ldots , z_{9\# S},y_1,\ldots, y_q]] \rightarrow R_\infty\]
and the module $M_\infty$ are defined in Section 2.8 of \cite{CEGGPS}.
Part \ref{proj} follows from Proposition 2.10 and the proof of Lemma 4.18(1) of \cite{CEGGPS}.
Part \ref{cm} follows from Lemma 4.18(1) of \cite{CEGGPS}.

The first identification in \ref{special} follows from Corollary 2.11 of \cite{CEGGPS}.
For the second identification, we have
\begin{align*}
(M_\infty(\tau)/\a)^\mathrm{d} \cong \Hom\cont_{K}(\tau, M_\infty^\mathrm{d})[\a]
&\cong \Hom\cont_{K}(\tau, (M_\infty/\a)^\mathrm{d}) \\
&\cong ((M_\infty/\a)^\mathrm{d})^{I,(-b,-c,-a)},
\end{align*}
where the first isomorphism follows from Remark 4.15 of \cite{CEGGPS}
and the final isomorphism follows from Frobenius reciprocity.

We now define the map $R_\infty \onto \bT_a$.
The surjection $R_\infty \onto \bT_a$ comes from the construction in Section 2.8 of \cite{CEGGPS}, as we now explain.
We will freely use the notation of Section 2.8 of \cite{CEGGPS} with the caveat that we exclude the Hecke operators at $v_1$ from our Hecke algebras (which allows us to deduce that the maps from $R_\infty$ to our Hecke algebras are surjections).
The patching argument produces compatible tuples $(\phi, M^{\square}, \psi, \alpha)$ of all levels $N \ge 1$. 
Since the intersection $\cap_N \mathfrak{d}_N$ is $0$, by completeness, the maps $\phi : R_\infty \onto R_{\cS}^{\mathrm{univ}}/\mathfrak d_N$ induce a surjection $R_\infty \onto R_{\cS}^{\mathrm{univ}}$.
Furthermore, in the paragraph before Corollary 2.11 of \cite{CEGGPS}, a map $R_{\cS}^{\mathrm{univ}} \rightarrow \bT_{\xi,\tau}^{S_p}(U^\p,\cO)_\m$ is described as the inverse limit of maps $R_{\cS}^{\mathrm{univ}} \rightarrow \bT_{\xi,\tau}^{S_p}(U_{m},\cO)_\m$ (or equivalently of maps $R_{\cS}^{\mathrm{univ}} \rightarrow \bT_{\xi,\tau}^{S_p}(U_{2N},\cO/\varpi^N)_\m$) at finite level.
Using the explicit generators given by the surjection $\bT^{S_p, \mathrm{univ}} \onto \bT_{\xi,\tau}^{S_p}(U^\p,\cO)_\m$ (again for us with operators at $v_1$ excluded) and Proposition 3.4.4(2) of \cite{CHT}, each of these maps at finite level is surjective.
By completeness, the map $R_{\cS}^{\mathrm{univ}} \rightarrow \bT_{\xi,\tau}^{S_p}(U^\p,\cO)_\m$ is surjective.
We conclude that the composition $R_\infty \rightarrow R_{\cS}^{\mathrm{univ}} \rightarrow \bT_{\xi,\tau}^{S_p}(U^\p,\cO)_\m$ is surjective.
Composing with the natural map $\bT_{\xi,\tau}^{S_p}(U^\p,\cO)_\m \onto \bT_a$ gives a surjective homomorphism $R_\infty \onto \bT_a$.

The proof of Corollary 2.11 of \cite{CEGGPS} shows that the first identification of \ref{special} is $R_\infty$-equivariant, with $R_\infty$
acting via the map $R_{\cS}^{\mathrm{univ}} \rightarrow \bT_{\xi,\tau}^{S_p}(U^\p,\cO)_\m$ on the right-hand side, 
since at each finite level $N$, the map $\psi : M^{\square}/\a \congto \tilde S_{\xi,\tau}(U_{2N},\cO/\varpi^N)_\m$ is $R_\infty$-equivariant,
with $R_\infty$ acting via the map $R_{\cS}^{\mathrm{univ}} \rightarrow \bT_{\xi,\tau}^{S_p}(U_{2N},\cO/\varpi^N)_\m$ on the right-hand side.
This shows that the second identification is also $R_\infty$-equivariant, with $R_\infty$ acting on the right-hand side via the map $R_\infty \onto \bT_a$ constructed above.
Part \ref{patchedfiber} follows from part \ref{special}, noting that $\a \subseteq \ker(R_\infty \to \bT_a)$ by part \ref{special}.
\end{proof}

Let $R$ be the $R_\infty$-subalgebra of $\End_{R_\infty}(M_\infty(\tau))$ generated by $U_1$ and $U_2$.

\begin{lm}\label{lm:R-local}
  The ring $R$ is local with maximal ideal $(\m_{\rbar},U_1,U_2)$.
\end{lm}

Here, by abuse, $\m_{\rbar}$ denotes the image of the maximal ideal of $R_\infty$ in $R$.

\begin{proof}
  As $M_\infty(\tau)$ is a finitely-generated $R_\infty(\tau)$-module, we deduce that the ring extension $R_\infty(\tau) \subseteq R$ is
  finite, hence $\m_{\rbar} \subseteq \rad(R)$. In particular, $\a \subseteq \rad(R)$, hence also $\sqrt{\a R} \subseteq \rad(R)$.
  By Theorem \ref{patchedmod}\ref{special}, the action of $R_\infty/\a$ on $M_\infty(\tau)/\a$ factors through $\bT_a$, hence
  the action of $R/\a$ on $M_\infty(\tau)/\a$ factors through $\bT$, giving rise to a surjective homomorphism $R/\a \onto \bT$.
  As the action of $R$ on $M_\infty(\tau)$ is faithful, the $R/\a$-module $M_\infty(\tau)/\a$ has full support by \cite[Lemma 2.2]{tay}
  (i.e., is nearly faithful in the terminology of \cite{tay}). Since $\bT$ is reduced, we get an induced isomorphism $(R/\a)_{\red} \congto \T$.
    It follows that $R$ is local, and $U_i \in \m_R$, as $U_i \in \m_\T$ by the proof of Theorem~\ref{main local global} ($i = 1$, $2$). 
  It is then easy to see that $\m_R = (\m_{\rbar},U_1,U_2)$.
\end{proof}

We now show that Theorem \ref{freeT} follows from the following theorem.

\begin{thm} \label{freeR}
Let $\rbar$ be as in Theorem \ref{main local global}\(i\) with $\FL(\rbar\vert_{G_{F_w}}) \ne \infty$. The module $M_\infty(\tau)$ is finite free over $R$.
\end{thm}

\begin{proof}[Proof of Theorem \ref{freeT}]
As $M_\infty(\tau)$ is finite free over $R$ by Theorem \ref{freeR}, $M_\infty(\tau)/\a$ is finite free over $R/\a$.
Hence the natural surjection $R/\a \onto \bT$ considered in the proof of Lemma~\ref{lm:R-local} is an isomorphism. 
We conclude by Theorem \ref{patchedmod}\ref{special}.
\end{proof}

\subsection{The Diamond--Fujiwara trick}

In this section, we prove Theorem \ref{freeR} using the method of Diamond and Fujiwara.
Let $B \subseteq {\GL_3}_{/\Z_p}$ be the algebraic subgroup of upper-triangular matrices.
For a dominant character $\lambda$ of $B$, let $W_{\cO_E}(\lambda)$ be the algebraic induced module $\( \Ind_{B}^{\GL_3} w_0 \lambda\)^{\mathrm{alg}} \otimes_{\Z_p} \cO_E$,
considered as representation of $K = \GL_3(\zp)$. (This is the characteristic $0$ version of the module defined in Section \ref{sec:LAS}.)
Let $W \defeq W_{\cO_E}(-c-1,-b,-a+1)$. 

\begin{prop} \label{flfree}The module $M_\infty(W)$ is non-zero and finite free over its support in $R_\infty$.
\begin{proof}
Note that $F(-c-1,-b,-a+1) \cong W \otimes_{\cO_E} \F$ by Proposition 3.18(ii) of \cite{florian-duke}.
Hence $M_\infty(W) \ne 0$ by Theorem \ref{theorem weight elimination}, Theorem \ref{patchedmod}\ref{special}, as well as Proposition 5.1.1 in \cite{le}.
Let $\psi_{\widetilde{v}} \defeq (a-1, b, c+1)$ and 
\[R_\infty(W)' \defeq \widehat{\otimes}_{v'|p} R^{\square,\psi_{\widetilde{v}'}}_{\widetilde{v}'} \widehat{\otimes} R^\square_{\widetilde{v}_1}[[x_1,\ldots, x_{q-3[F^+:\Q]}]],\]
where $R^{\square,\psi_{\widetilde{v}}}_{\widetilde{v}}$ is the framed crystalline deformation ring for $\rbar|^\vee_{G_{F_{\wt v}}}$ of Hodge--Tate weights $\psi_{\wt v}+(2,1,0)$.
Then $M_\infty(W)$ is supported on $\Spec R_\infty(W)'$ by Lemma 4.17(1) of \cite{CEGGPS}.
(Note the dual in the definition of ``Hodge--Tate weights prescribed by $\sigma_{\mathrm{alg}}$'' in Section 4 of \cite{CEGGPS}.)
By the proof of Lemma 4.18 of \cite{CEGGPS}, $M_\infty(W)$ is a maximal Cohen--Macaulay $R_\infty(W)'$-module.
By the choice of $v_1$, $R^\square_{\widetilde{v}_1}$ is formally smooth over $\cO_E$, see Lemma 2.5 of \cite{CEGGPS}.
By Lemma 2.4.1 of \cite{CHT}, $R^{\square,\psi_{\widetilde{v}'}}_{\widetilde{v}'}$ is formally smooth over $\cO_E$ for all $v'|p$. 
As $R_\infty(W)'$ is a regular local ring, $M_\infty(W)$ is a projective and hence free $R_\infty(W)'$-module by the Auslander--Buchsbaum--Serre theorem and the Auslander--Buchsbaum formula.
\end{proof}
\end{prop}

Let $n$ be the (free) rank of $M_\infty(W)$ over its support.

\begin{lm} \label{generate}
Let $\rbar$ be as in Theorem \ref{main local global}\(i\) with $\FL(\rbar\vert_{G_{F_w}}) \ne \infty$. 
The $R$-module $M_\infty(\tau)$ is generated by $n$ elements.
\begin{proof}
Let $\taubar \defeq \tau \otimes_{\cO_E} \F$.
By Nakayama's lemma, it suffices to show that $M_\infty(\taubar)/(\m_{\rbar},U_1,U_2)$ has dimension $n$.
Consider
\begin{align*}
(M_\infty(\taubar)/(\m_{\rbar},U_1,U_2))^\vee &\cong \Hom_{K} (\taubar,M_\infty^\vee)[\m_{\rbar},U_1,U_2] \\
& \cong \Hom_{K} (\taubar,M_\infty^\vee[\m_{\rbar}])[U_1,U_2] \\
& \cong \Hom_{K} (\taubar/M,M_\infty^\vee)[\m_{\rbar}]\\
& \cong (M_\infty(\taubar/M)/\m_{\rbar})^\vee,
\end{align*}
where the third isomorphism holds by Lemma \ref{lemma killed by U1 U2}, and where $M$ is the minimal subrepresentation of $\taubar$ containing $F(-c,-a,-b-p+1)$
and $F(-a+p-1,-b,-c-p+1)$ as Jordan--H{\"o}lder factors.
By the isomorphism $M_\infty(\taubar)/(\m_{\rbar},U_1,U_2) \cong M_\infty(\taubar/M)/\m_{\rbar}$, it suffices to show that $M_\infty(\taubar/M)/\m_{\rbar}$ has dimension $n$.
By Proposition 2.2.2 of \cite{le} and Theorem \ref{theorem weight elimination}, $M$ is of length 3 and $(\taubar/M)^\vee$ contains no element of $W_w(\rbar)$ other than 
$F(a-1,b,c+1)$. Hence by Proposition 5.1.1 of \cite{le}, the natural inclusion $M_\infty(F(-c-1,-b,-a+1)) \into M_\infty(\taubar/M)$ is an isomorphism.
We conclude that the maps 
\[M_\infty(W)/\m_{\rbar} \rightarrow M_\infty(F(-c-1,-b,-a+1))/\m_{\rbar} \rightarrow M_\infty(\taubar/M)/\m_{\rbar}\]
are isomorphisms of $n$-dimensional $\F$-vector spaces, by definition of $n$.
\end{proof}
\end{lm}

\begin{proof}[Proof of Theorem \ref{freeR}]
By Lemma \ref{generate}, we can and do fix a surjection $R^n \onto M_\infty(\tau)$ with kernel $P$.
This induces an exact sequence 
\begin{equation} \label{exactproj}
0 \rightarrow P[1/p] \rightarrow (R[1/p])^n \rightarrow M_\infty(\tau)[1/p] \rightarrow 0.
\end{equation}
We now freely use the notation of Section 4 of \cite{CEGGPS}.
First, we note that $R_\infty(\tau)[1/p] \rightarrow R[1/p]$ is an isomorphism by Theorem 4.1 and Lemma 4.17(2) of \cite{CEGGPS}.
Note that while the construction of $M_\infty$ is slightly different in this paper,
the proof of Lemma 4.17(2) still applies verbatim.

By Proposition 2.4.7 of \cite{EGH} (and its easy converse) $\tau \otimes_{\cO_E} \qpb$ is an irreducible smooth representation of $K$ that is associated to the inertial type 
$\teich{\omega}^{-b}\oplus\teich{\omega}^{-c}\oplus\teich{\omega}^{-a}$ by the inertial local Langlands correspondence (\cite{CEGGPS}, Theorem 3.7).
Hence by Theorem \ref{patchedmod}\ref{cm}, $M_\infty(\tau)[1/p]$ is a projective $R_\infty(\tau)[1/p]$-module.
We claim that it is of constant rank $n$.
The rank is clearly no larger than $n$ by (\ref{exactproj}).
Thus it suffices to show that each irreducible component of $\Spec (R_\infty(\tau)[1/p])$, or equivalently of $\Spec R_\infty(\tau)$, contains an $E$-point over which the fiber of $M_\infty(\tau)[1/p]$ has $E$-dimension at least $n$.

First, we claim that each irreducible component $Z$ of $\Spec R_\infty(\tau)$ has an $E$-point in the closed subscheme $\Spec \bT_a$.
Indeed as $M_\infty(\tau)$ is finite free over $S_\infty$ by Theorem \ref{patchedmod}\ref{proj}, $R_\infty(\tau)$ is a subring of a matrix ring over $S_\infty$ and hence a finite torsion-free $S_\infty$-module.
From this, we see that $R_\infty(\tau)$ is Cohen--Macaulay and in particular equidimensional.
As $Z \rightarrow \Spec S_\infty$ is a finite map between irreducible schemes of the same dimension, it is surjective.
In particular, the base change $Z\times_{\Spec S_\infty} \Spec (S_\infty/\mathfrak{a}) \subseteq \Spec R_\infty(\tau)/\mathfrak{a}$ has a non-empty generic fibre.
Just as in the proof of Lemma~\ref{lm:R-local}, we have an isomorphism $(R_\infty(\tau)/\mathfrak{a})_{\mathrm{red}} \congto \bT_a$.
Since $E$ is sufficiently large, $\Spec (\bT_a[1/p])$ is a union of copies of $\Spec E$.
We conclude that the underlying reduced subscheme of the generic fiber of $Z\times_{\Spec S_\infty} \Spec (S_\infty/\mathfrak{a})$ is a non-empty union of copies of $\Spec E$ in $\Spec \bT_a$.

Second, we claim that for each $\cO_E$-point of $\Spec \bT_a$ with corresponding prime ideal $\wp_r$, $M_\infty(\tau)/\wp_r$ has generic rank at least $n$.
Let $M$ be the minimal subrepresentation of $\taubar$ containing $F(-c-1,-b,-a+1)$ as a Jordan--H{\"o}lder factor.
By Proposition 2.2.2 of \cite{le} and Theorem \ref{theorem weight elimination}, $(\taubar/M)^\vee$ contains no element of $W_w(\rbar)$. By Proposition 5.1.1 of \cite{le}, the natural map $M_\infty(M) \rightarrow M_\infty(\taubar)$ is an isomorphism.
The natural surjection $M_\infty(M) \onto M_\infty(F(-c-1,-b,-a+1))$ gives a surjection $M_\infty(\taubar) \onto M_\infty(F(-c-1,-b,-a+1))$, and therefore a surjection $M_\infty(\taubar)/\m_{\rbar} \onto M_\infty(F(-c-1,-b,-a+1))/\m_{\rbar}$.
Since $M_\infty(F(-c-1,-b,-a+1))/\m_{\rbar}$ has dimension $n$ by definition of $n$, $M_\infty(\taubar)/\m_{\rbar}$ has dimension at least $n$.
For any $\cO_E$-point of $\Spec \bT_a$ with corresponding prime ideal $\wp_r$, $M_\infty(\tau)/\wp_r$ is $p$-torsion free by Theorem \ref{patchedmod}\ref{patchedfiber} and so has generic rank at least $n$.

As $M_\infty(\tau)[1/p]$ is projective of rank $n$ over $R[1/p]$ ($\cong R_\infty(\tau)[1/p]$), we conclude by (\ref{exactproj}) that $P[1/p]$ is projective of rank $0$ and hence $0$.
As $R$ is $p$-torsion free, $P$ is ~$0$.
\end{proof}

Let $R_2$ be the $R_\infty$-subalgebra of $\End_{R_\infty}(M_\infty(\tau))$ generated by $U_2$, and let $\T_2$ be the
$\cO_E$-subalgebra of $\End(S(U^v,\widetilde{V}')_{\m_{\rbar}}^{I,(-b,-c,-a)})$ generated by $\bT^{\cP}$ and $U_2$.

\begin{thm}
  Let $\rbar$ be as in Theorem \ref{main local global}\(i\)--\(ii\). We have $R_2 = R$ and $\T_2 = \T$, i.e.\ $U_1 \in R_2$ and $U_1 \in \T_2$.
  Moreover, $S(U^v,V')[\m_{\rbar}]^{I,(-b,-c,-a)}[U_2] = S(U^v,V')[\m_{\rbar}]^{I,(-b,-c,-a)}[U_1,U_2]$.
\end{thm}

\begin{proof}
  The argument for Theorem~\ref{freeR} goes through to show that the module $M_\infty(\tau)$ is also finite free over $R_2$,
  and the free rank is the same as for $R$. (The main thing that changes in the proof is that the analogue of $M$ in the
  proof of Lemma~\ref{generate} now has length 2, and we need that $\FL(\rbar\vert_{G_{F_w}}) \not \in \{0,\infty\}$
  in order to deduce that $(\taubar/M)^\vee$ contains no element of $W_w(\rbar)$.) Now $R_2 \subseteq R$ is a finite ring
  extension, hence $\m_R \cap R_2 = \m_{R_2}$. As $E$ is sufficiently large, the residue fields of $R_2$ and $R$ are both
  equal to $\F$. It follows that the surjection $M_\infty(\tau)/\m_{R_2} \onto M_\infty(\tau)/\m_R$ is an isomorphism.
  Therefore any $R_2$-basis of $M_\infty(\tau)$ is also an $R$-basis, so $R_2 = R$. The other claims follow.
  \end{proof}

\begin{rk}\label{rk:no-shadow}
  If $W_w(\rbar) = \{F(a-1,b,c+1)\}$, then we even deduce (by the same argument) that $R_\infty(\tau) = R$, $\T_a = \T$, and
  that $U_1$, $U_2$ annihilate $S(U^v,V')[\m_{\rbar}]^{I,(-b,-c,-a)}$.
\end{rk}

\begin{rk} \label{rk:dual-lattice}
  The analogues of all the above results in this section hold for the ``dual'' lattice $\tau' \defeq \Ind_I^K
  (\teich{\omega}^{-c}\otimes \teich{\omega}^{-a} \otimes \teich{\omega}^{-b}) \otimes_{\Z_p} \cO_E$, provided
  one interchanges $U_1$ with $U_2$ and the condition $\FL(\rbar\vert_{G_{F_w}}) \ne \infty$
  with $\FL(\rbar\vert_{G_{F_w}}) \ne 0$.
\end{rk}

We now show that Theorem~\ref{main local global} applies to Galois representations of any possible invariant outside
$0$, $\infty$.

\begin{thm}\label{thm:all-invariants-arise}
  Suppose that $\rhobar : G_{\qp} \to \GL_3(\F)$ is upper-triangular, maximally non-split, and generic. Then, after possibly replacing $\F$ by a
  finite extension, there exist a CM field $F$, a Galois representation $\rbar : G_F \to \GL_3(\F)$, a place $w|p$ of $F$,
  groups $G_{/F^+}$ and $\cG_{/\cO_{F^+}}$, and a compact open subgroup $U^v$ \(where $v = w|_{F^+}$\) satisfying all hypotheses
  of the setup in Section \ref{sec:setup} such that $\rbar|_{G_{F_w}} \cong \rhobar$.
 
  In particular, if $\FL(\rhobar) \not \in \{0,\infty\}$, Theorem~\ref{main local global} applies to $\rbar$.
\end{thm}

\begin{proof}
  We suppose that $\rhobar$ is as in \eqref{Galoisrep} with $(a_2,a_1,a_0) = (a,b,c)$.
  We note that $\rhobar$ satisfies the two hypotheses of \cite{CEGGPS}, \S2.1: $p \nmid 6$ by genericity, and $\rhobar$ has a
  potentially diagonalizable lift of regular weight: for example, by
  Corollary~\ref{corollary weight elimination} or by \cite{gee-geraghty}, Lemma 3.1.5 there exists a
  crystalline lift with distinct Hodge--Tate weights in the Fontaine--Laffaille range. Therefore the procedure of \cite{CEGGPS}, \S 2.1,
  building on \cite{EG}, \S A, yields a ``suitable globalization'' $\rbar : G_F \to \GL_3(\F)$ (after possibly replacing $\F$)
  for a CM field $F$ with $F/F^+$ unramified at all finite places. For each place
  $v' \in S_p^+$ there exists a place $\wt v'$ of $F$ lying over $v'$ with $F_{\wt v'} = \qp$ and $\rbar|_{G_{F_{\wt v'}}} \cong \rhobar$, 
  which is Fontaine--Laffaille and regular. We see that hypotheses \ref{rbar-unram}--\ref{zeta-p} in
  \S\ref{sec:setup} are satisfied. 

  Choose a unitary group $G$ with integral model $\cG$ as in \S\ref{sec:setup}. Fix any $v \in S_p^+$, let $w \defeq \wt v$,
  and define $\wt V' = \bigotimes_{v'\in S_p^+\setminus\{v\}}W_{{a}_{v'}}$ with $a_{\wt v'} = (a-1,b,c+1)$ for any $v' \in S_p^+ \setminus \{v\}$. 
  There exists a compact open subgroup $U = \prod_{v'} U_{v'} \le G(\A_{F^+}^\infty)$ with $U_{v'} =
  \cG(\cO_{F^+_{v'}})$ (resp.\ $U_{v'}$ hyperspecial) for all places $v'$ which split (resp.\ are inert) in $F$ such that $S(U^v,
  \wt V')_{\m_{\rbar}} \ne 0$, where $\cP$ and $\m_{\rbar}$ are defined as in \S\ref{sec:setup}. 
  (Note that a priori the proof of \cite{EG}, Corollary A.7 does not give us any information about the type $(\lambda,\tau)$
  in \cite{EG}, Definition 5.3.1. However, we can get the desired type by choosing appropriate lifts locally at places dividing $p$ in applying \cite{BLGGT}, Theorem 4.3.1 
  in the proof of \cite{EG}, Lemma A.5.)
  By choosing a place $v_1$ and redefining $U_{v_1}$ as in \S\ref{sec:setup}, we see that hypothesis \ref{automorphy} in \S\ref{sec:setup} is satisfied.
\end{proof}

\begin{appendices}
\section{Some integral $p$-adic Hodge theory}
\label{sec:appendix}
The first aim of this appendix is to collect some results in integral $p$-adic Hodge theory. Several of these are
slight generalizations of known results in the literature, incorporating coefficients and descent data.
The second aim of the appendix is to prove some lemmas and propositions of Sections~\ref{sec:p-adic-hodge} and \ref{sec:breuil-submod}. 

\subsection{On certain categories of torsion modules}
\label{appendix: Modules}

The goal of this section is to recall the definition of various categories of $p$-torsion modules with additional structures
(Breuil modules, Kisin modules, $\phi$-modules) and the relations between them.

Recall that $K_0= W(k)[\frac{1}{p}]$ and that $K= K_0(\varpi)$, where $\varpi=\sqrt[e]{-p}$ and $e \ge 1$ divides $p^{[k:\Fp]}-1$. Also recall our definitions of
$\barS_k = k[u]/u^{ep}$ and $S_{W(k)}$ in Section~\ref{sec:p-adic-hodge}. We will always suppose that $0 \le r \le p-2$. For the purpose
of establishing Proposition~\ref{proposition comparison} we will consider general $e$, but
sometimes we will specialize to $e = 1$. We hope that this will not cause any confusion.

If $\cA$ is any additive category we let $\F\text{-}\cA$ denote the $\F$-linear additive category whose objects consist
of pairs $(A,i)$, where $A \in \cA$ and $\F \to \End(A)$ is a ring homomorphism (the morphisms are the obvious ones). We
will often use below, without further comment, that a functor $\cA \to \cB$ of additive categories gives rise to a functor 
$\F\text{-}\cA \to \F\text{-}\cB$.

We consider the general setting of \S \ref{subsubBrModdd}. In particular, we write $\BrMod$ to denote the category of Breuil modules over $\barS_k$ and $\BrModN$ for the category of Breuil modules without monodromy (also called \emph{quasi-Breuil modules}), defined in the evident way. Note that these categories are respectively denoted by $\widetilde{\underline{\cM}}^r$, $\widetilde{\underline{\cM}}^r_0$ in \cite{Caruso-Crelle}, \cite{breuil-normes}.

In \cite{breuil-inventiones}, \S 2.2.1 (cf.\ also \cite{Caruso-Crelle}, \S 2.1) the authors introduce certain categories of modules, noted as $\underline{\cM}^r$, $\underline{\cM}_0^r$. We actually only consider their full subcategories formed by $p$-torsion objects, though we keep the same notation. The objects of these subcategories are finite free
modules over $S_{W(k)}/pS_{W(k)}$. We have the following commutative diagram:
\begin{align*}
\xymatrix{
\underline{\cM}^r\ar[d]\ar^-{\cT}_-{\sim}[r]&\BrMod\ar[d]\\
\underline{\cM}^r_0\ar_-{\cT_0}^-{\sim}[r]&\BrModN
}
\end{align*}
where the functors $\cT$, $\cT_0$ inducing equivalences of categories are defined in \cite{breuilAENS98}, \S2.2.2 when $e = 1$ and 
\cite{Caruso-Crelle}, \S2.3 in general, and the vertical arrows are forgetful functors.
(Strictly speaking, this diagram and the ones that follow only commute up to equivalence.)
 
Choose a  sequence $(p_n)_n\in\overline{\Q}_p^{\N}$ verifying $p_{n+1}^p=(-1)^{p-1} p_{n}$  for all $n$ and $p_0\defeq -p$, so $N_{K_0(p_{n+1})/K_0(p_n)}(p_{n+1}) = p_n$.
Let $k((\underline{p}))$ be the field of norms associated to the extension $(K_{0})_{\infty}/K_0$, where $(K_{0})_{\infty}\defeq \cup_{n\in\N}K_0(p_n)$,
and let $k[[\underline{p}]]$ be its ring of integers. In particular, $\underline{p}$ is identified with the sequence $(p_n)_n$.
We write  $\mathfrak{Mod}^r_{k[[\underline{p}]]}$ for the category of $(\phi, k[[\underline{p}]])$-modules of height $\le r$ defined in \cite{breuil-normes}, \S 2.3.

We have an equivalence of categories:
\begin{equation*}
M_{\mathfrak{S}}:\,\,\mathfrak{Mod}^r_{k[[\underline{p}]]}{\stackrel{\sim}{\longrightarrow}}\BrModN
\end{equation*}
(cf.\ \cite{breuil-normes}, \S4, where the functor is noted by $\tilde\Theta_r$).

We finally define the category $\mathfrak{Mod}_{k((\underline{p}))}$ of \'etale $(\phi,k((\underline{p})))$-modules as the category of finite-dimensional $k((\underline{p}))$-vector spaces $\mathfrak{D}$  endowed with a semilinear map $\phi:\mathfrak{D}\rightarrow \mathfrak{D}$ (with respect to the Frobenius on $k((\underline{p}))$) and inducing an isomorphism $\phi^{\ast}\mathfrak{D}\rightarrow \mathfrak{D}$ (with obvious morphisms between objects).

By work of Fontaine \cite{fontaine-fest}, we have an exact anti-equivalence
\begin{align*}
\F\mbox{-}\mathfrak{Mod}_{k((\underline{p}))}& \stackrel{\sim}{\longrightarrow}\mathrm{Rep}_{\F}(G_{(K_0)_{\infty}})\\
\mathfrak{D}&\longmapsto \mathrm{Hom}(\mathfrak{D},k((\underline{p}))^{\text{sep}}).
\end{align*}

When $e = 1$ we have the fully faithful embedding 
$\mathcal{F}_r:\mathcal{FL}^{[0,r]}{\longrightarrow}\BrModN$ 
of \cite{breuil-normes}, \S5 (via the equivalence $\cT_0$ above). By composition we obtain a fully faithful functor 
\begin{equation*}
\xymatrix{
\F\text{-}\mathcal{FL}^{[0,r]}
\ar@/^2pc/^-{\mathcal{F}}[rrr]\ar^-{\mathcal{F}_r}[r]&\FBrModN
&\F\text{-}\mathfrak{Mod}^r_{k[[\und{p}]]}\ar_-{M_{\mathfrak{S}}}^-{\sim}[l]\ar[r]
&\F\text{-}\mathfrak{Mod}_{k((\und{p}))},
}
\end{equation*}
where the last functor in the sequence is localization at $\underline p$ (which is fully faithful by \cite{breuil-normes}, Proposition 2.3.7).

We are now going to define the analogous categories and functors with descent data.
Recall that $\varpi= \sqrt[e]{-p}\in K$. There is a unique sequence $(\varpi_n)_n\in\overline{\Q}_p^{\N}$  such that $\varpi_n^e=p_n$, $\varpi_{n+1}^p = (-1)^{p-1} \varpi_n$  for all $n$, and $\varpi_0\defeq \varpi$. In particular, $N_{K(\varpi_{n+1})/K(\varpi_n)}(\varpi_{n+1}) = \varpi_n$.
By letting $K_{\infty}\defeq \cup_{n\in\N}K(\varpi_n)$ we have a canonical isomorphism $\Gal(K_{\infty}/(K_{0})_{\infty})\isoto \Gal(K/K_0)$ which lets us identify $\omega_\varpi$ with a character of $\Gal(K_{\infty}/(K_{0})_{\infty})$.

As $\und\varpi^e = \und p$ and $e$ divides $|k\s|$, the field of norms $k((\underline{\varpi}))$ associated to $K_\infty/K$ is a cyclic extension of $k((\und p))$ of degree $e$.
Note that $\Gal(K_{\infty}/(K_0)_{\infty}) \cong \Gal(k((\und\varpi))/k((\und p)))$ acts on $k((\underline{\varpi}))$; concretely this action is determined by 
$
g\cdot \underline{\varpi}= \omega_{\varpi}(g)\underline{\varpi}
$
for $g\in \Gal(K_{\infty}/(K_0)_{\infty})$.

We write $\mathfrak{Mod}^r_{k[[\und{\varpi}]]}$ for the category of $(\phi, k[[\underline{\varpi}]])$-modules of height $\le r$, described in \cite{caruso-fourier}, \S2.1 (where it is denoted by $\mathrm{Mod}^{\phi}_{/\tilde{\mathfrak{S}}}$).
By \cite{caruso-fourier}, Th\'eor\`eme 2.1.2 (building on Breuil \cite{breuil-normes}, Th\'eor\`eme 4.1.1)
we have an exact equivalence of categories
\begin{equation*}
M_{\mathfrak{S}}:\,\,\mathfrak{Mod}^r_{k[[\und{\varpi}]]}{\stackrel{\sim}{\longrightarrow}}\BrModN.
\end{equation*}

We define in the evident way the category $\mathfrak{Mod}^r_{k[[\und{\varpi}]],\,\mathrm{dd}}$ consisting of $(\phi, k[[\underline{\varpi}]])$-modules endowed with a descent data from $K$ to $K_0$, i.e.\ a semilinear action of $\Gal(K_{\infty}/(K_0)_{\infty})$ that commutes with $\phi$ (the morphisms now being required to be compatible with the descent datum). 
\begin{lm}
\label{lemma equivalence dd}
There exists an exact covariant functor
$$
M_{\mathfrak{S}}:\, \F\text{-}\mathfrak{Mod}^r_{k[[\und{\varpi}]],\, \mathrm{dd}}\rightarrow \F\text{-}\BrModddN
$$
which establishes an equivalence of categories.
\end{lm}

Here, the category $\F\text{-}\BrModddN$ has the same definition as $\F\text{-}\BrModdd$, except that we drop monodromy.

\begin{proof}
We may assume that $\F = \Fp$. By above we only need to check compatibility with the descent data. For this we follow the strategy outlined in \cite{EGH}, proof of Proposition 3.2.6.
More precisely, given  an object $\mathfrak{M}\in \mathfrak{Mod}^r_{k[[\und{\varpi}]],\,\mathrm{dd}}$ and an element $g\in\Gal(K/K_0)$ we define its twist $\mathfrak{M}^{(g)}$ as $k[[\und{\varpi}]] \otimes_{g,k[[\und{\varpi}]]} \mathfrak{M}$ with $\phi$ acting diagonally. 
It is easily checked that $\mathfrak{M}^{(g)}\in \mathfrak{Mod}^r_{k[[\und{\varpi}]]}$ and that $\widehat{g}:\mathfrak{M} \rightarrow \mathfrak{M}$ induces
a morphism $\widehat{g}:\mathfrak{M}^{(g)}\rightarrow \mathfrak{M}$ in $\mathfrak{Mod}^r_{k[[\und{\varpi}]]}$; moreover,  
$\widehat{g h}=\widehat{g}\circ\widehat{h}^{(g)}$ for any $g,h\in\Gal(K/K_0)$.
It is formal to verify that the datum of a family of morphisms $\widehat g: \mathfrak{M}^{(g)}\rightarrow \mathfrak{M}$ verifying the cocycle relation is equivalent to a descent datum on $\mathfrak{M}$ from $K$ to $K_0$. 

Recall now that the functor ${M}_{\mathfrak{S}}: \mathfrak{Mod}^r_{k[[\und{\varpi}]]}\rightarrow \BrModN$ is defined as the base change induced by the morphism $k[[\und{\varpi}]]\onto\barS_k\stackrel{\varphi}{\rightarrow}\barS_k$, where the first map sends $\sum \lambda_i \und\varpi^i$ to $\sum \lambda_i u^i$. Therefore  
$M_{\mathfrak{S}}(\mathfrak{M}^{(g)})\cong (M_{\mathfrak{S}}(\mathfrak{M}))^{(g)}$, where the twist in the right-hand side is the one defined in \cite{EGH}, proof of Proposition 3.2.6.

By functoriality, the quasi-Breuil module $\cM\defeq {M}_{\mathfrak{S}}(\mathfrak{M})$ is thus endowed with a family of morphism
$\cM^{(g)}\rightarrow\cM$ verifying the above cocycle relation, i.e.\ $\cM$ is endowed with descent data from $K$ to $K_0$. 
We conclude that the composite functor
$\mathfrak{Mod}^r_{k[[\und{\varpi}]],\,\mathrm{dd}}\rightarrow \mathfrak{Mod}^r_{k[[\und{\varpi}]]}\stackrel{\sim}{\longrightarrow} \BrModN$
factors though $\BrModddN\rightarrow \BrModN$.
Its exactness and the induced equivalence follow as a formal consequence of the analogous properties of $M_{\mathfrak{S}}$, cf.\ \cite{caruso-fourier}, Th\'eor\`eme 2.1.2.
\end{proof}

We can finally introduce the category $\mathfrak{Mod}_{k((\und{\varpi})),\,\mathrm{dd}}$ of \'etale $(\phi,k((\underline{\varpi})))$-modules with descent data: an object $\mathfrak{D}$ is defined in the analogous way as for the category  $\mathfrak{Mod}_{k((\und{p}))}$, but we moreover require that $\mathfrak{D}$ is endowed with a semilinear action of $\Gal(K/K_0)$ that commutes with $\phi$.

We can now define a functor by composition: 
\begin{equation*}
\xymatrix{
\FBrModdd\ar@/^2pc/^-{M_{k((\underline{\varpi}))}}[rrr]\ar[r]&\FBrModddN
&\F\text{-}\mathfrak{Mod}_{k[[\und{\varpi}]],\,\mathrm{dd}}^r\ar_-{M_{\mathfrak{S}}}^-{\sim}[l]\ar[r]
& \F\text{-}\mathfrak{Mod}_{k((\und{\varpi})),\,\mathrm{dd}}.
}
\end{equation*}

\subsection{Some basic lemmas on Breuil modules}
\label{sec:some-basic-lemmas}

\begin{lm}
\label{Lemma1Note Florian}
Let $\cN$ be an $\barS$-submodule of a Breuil module $\cM\in\FBrMod$.
Then $\cN$ is an $\barS_k$-direct summand of $\cM$ if and only if $\cN$ is free as $\barS$-module.
\end{lm}
\begin{proof}
By Baer's criterion the noetherian ring $\barS_k$ is self-injective; in particular an $\barS_k$-module is injective if and only if it is projective. 
As $\barS_k$ is local, $\cN$ is an $\barS_k$-direct summand of $\cM$ if and only if $\cN$ is free as $\barS_k$-module.
This is equivalent to $\cN$ being free as $\barS$-module by \cite{EGH}, Remark 3.2.1.
\end{proof}

\begin{lm}
  \label{lm:brmod-fil}
  If $\cM\in\BrMod$, then $\Fil^r\cM / u^e \Fil^r\cM$ is a free $k[u]/u^e$-module of rank $\rank_{\barS_k}\cM$ and
  the map 
  \begin{equation}\label{eq:1}
    \barS_k \otimes_{\varphi, k[u]/u^e} \Fil^r\cM / u^e \Fil^r\cM \to \cM
  \end{equation}
  given by $1 \otimes \varphi_r$ is an isomorphism.
\end{lm}

\begin{proof}
  (Cf.\ \cite{breuilAENS98}, Lemme 2.2.1.2 when $e = 1$.)
  As $\cM$ is finite free over $\barS_k$ and $u^{er} \cM \subseteq \Fil^r\cM$ the first claim follows.
  Clearly the map~\eqref{eq:1} is well defined and surjective. It is an isomorphism for dimension reasons.
\end{proof}

\begin{corollary}
  \label{cor:brmod-fil}
  If $\cM\in \underline{\cM}^r$, let 
  \[ \Fil^{r+1} \cM \defeq (\Fil^1 S_{W(k)})(\Fil^r \cM) + (\Fil^{r+1} S_{W(k)}) \cM. \]
  Then the natural map 
  \[ \Fil^r \cM/ \Fil^{r+1} \cM \to \Fil^r \cT(\cM)/u^e \Fil^r \cT(\cM) \]
  is an isomorphism of $S_{W(k)}/(p, \Fil^1 S_{W(k)}) \cong k[u]/u^e$-modules, and 
  the map 
  \begin{equation}\label{eq:3}
    S_{W(k)}/p \otimes_{\varphi, k[u]/u^e} \Fil^r\cM / \Fil^{r+1}\cM \to \cM
  \end{equation}
  given by $1 \otimes \varphi_r$ is an isomorphism.
\end{corollary}

\begin{proof}
  The first claim follows, as $\ker(\cM \to \cT(\cM)) = (\Fil^p S_{W(k)})\cM \subseteq \Fil^{r+1} \cM$.
  Thus by Lemma~\ref{lm:brmod-fil} both sides of (\ref{eq:3}) are free $S_{W(k)}/p$-modules of the same rank.
  As $S_{W(k)}$ is local, it suffices to check that (\ref{eq:3}) is an isomorphism after the base change $S_{W(k)}/p \onto \barS_k$,
  and hence we are done by Lemma~\ref{lm:brmod-fil}.
\end{proof}

\subsection{Functors towards Galois representations}
\label{Appendix: Galois}

The aim of this section is to briefly recall the constructions of various functors towards Galois representations defined on the categories of \S \ref{appendix: Modules}. 
The functors towards Galois representations are defined via certain period rings introduced in \cite{breuil-inventiones}, \S 2.2.2 (cf.\ also \cite{Caruso-Crelle}, \S 2.2 and 2.3).  
We have a natural surjective morphism $S_{W(k)}\rightarrow \barS_k$, as well as $S_{W(k)}$-algebras $\widehat{A}_{\mathrm{st}}$, $A_{\mathrm{cris}}$,  which are endowed with an action of $G_K$. We have
the following commutative diagram:
\begin{equation}
\label{period rings}
\begin{gathered}
  \xymatrix{
    \widehat{A}_{\mathrm{st}}/p\ar@{->>}[r]\ar@{->>}[d]&\widehat{A}_{\mathrm{st}}\otimes_{S_{W(k)}}\barS_k\defeq \widehat{A}\ar@{->>}[d]\\
    A_{\mathrm{cris}}/p\ar@{->>}[r]&{A}_{\mathrm{cris}}\otimes_{S_{W(k)}}\barS_k\defeq \widehat{A}_0 }
\end{gathered}
\end{equation}
which becomes equivariant by endowing the rings\footnote{We remark that the rings $\widehat{A}$, $\widehat{A}_0$, $A_{\mathrm{cris}}/p$ admit alternative descriptions. We refer the reader to \cite{Caruso-Crelle}, \S 2.2 and 2.3, where it is shown (Lemme 2.3.2) that $A_{\mathrm{cris}}/p\cong R^{DP}$, $R^{DP}$ being a certain period ring introduced in \cite{breuil-normes}, \S 3.2, and where the ring $\widehat{A}$ is explicitly described.}
$\widehat{A}$, $\widehat{A}_0$ with the actions of $G_K$, $G_{K_\infty}$ inherited from $\widehat{A}_{\mathrm{st}}/p$, $A_{\mathrm{cris}}/p$, and $\barS_k$.
The choice of $\und\varpi$ lets us extend these actions to $G_{K_0}$ (respectively $G_{(K_0)_\infty}$ in the case of $\wh A_0$), as explained in \cite{EGH}, \S  3.1.
The horizontal maps in the diagram are then $G_{K_0}$-equivariant, while the vertical maps are $G_{(K_0)_\infty}$-equivariant.

Let us consider the contravariant functors
\begin{align*}
\underline{\cM}^r&\rightarrow \mathrm{Rep}_{\fp}(G_K)&
\underline{\cM}^r_0&\rightarrow \mathrm{Rep}_{\fp}(G_{K_{\infty}})\\
\cM&\mapsto\mathrm{Hom}(\cM,\widehat{A}_{\mathrm{st}}/p)&
\cM_0&\mapsto\mathrm{Hom}(\cM_0,{A}_{\mathrm{cris}}/p)
\end{align*}
(where homomorphisms, here and in the following, respect all structures)
and their analogues on the categories $\BrMod$, $\BrModN$, defined by replacing the rings $\widehat{A}_{\mathrm{st}}/p$, ${A}_{\mathrm{cris}}/p$ with $\widehat{A}$, $\widehat{A}_0$ respectively.
(The functor $\mathrm{Hom}(-,\widehat{A})$ on the category $\BrMod$ is often denoted by $\Tst^*$. This is compatible with the notation $\Tst^*$ in \S\ref{sec:p-adic-hodge} for the functor on the category $\FBrModdd$.)

These functors have been extensively studied in \cite{breuil-inventiones}, \S 2.3.1, \cite{breuil-normes}, \S 3.2 and \S 4, \cite{Caruso-Crelle}, \S 2.2 and \S 2.3 and the morphisms between the various period rings (\ref{period rings}) induce natural transformation between the functors toward Galois representations. 
The following result provides us with a precise description of the relations between the categories and functors introduced above.

\begin{prop}
\label{proposition comparison classic}
The natural maps \(\ref{period rings}\) induce a commutative diagram
\begin{equation*}
\xymatrix@=5pc{
\F\text{-}\underline{\cM}^r\ar_-{\cT}^-\sim[r]\ar[d]\ar@/^2pc/^-{\text{\tiny$\mathrm{Hom}(-,\widehat{A}_{\mathrm{st}}/p)$}}[rr]
&\F\text{-}\BrMod\ar[d]\ar_-{\text{\tiny$\mathrm{Hom}(-,\widehat{A})$}}[r]
&\mathrm{Rep}_{\F}(G_K)\ar^{\res}[d]\\
\F\text{-}\underline{\cM}^r_0\ar^-{\cT_0}_-\sim[r]\ar@/_2pc/_-{\text{\tiny$\mathrm{Hom}(-,{A}_{\mathrm{cris}}/p)$}}[rr]
&\F\text{-}\BrModN\ar^-{\text{\tiny$\mathrm{Hom}(-,\widehat{A}_0)$}}[r]
&\mathrm{Rep}_{\F}(G_{K_{\infty}})
}
\end{equation*}
where the unlabeled morphisms are forgetful functors.
\end{prop}
\begin{proof}
We may assume that $\F=\Fp$.
The left square is commutative by above, and the external square is commutative by \cite{breuil-inventiones}, Lemme 2.3.1.1.

We now indicate why the upper and lower triangles are commutative, completing the proof. For the upper triangle this follows
from \cite{Caruso-Crelle}, Proposition 2.3.4. However, the second paragraph of the published proof needs to be fixed as
follows. (We thank Xavier Caruso for the argument.)  Given $\o\psi : \cT(\cM) \to \wh A$, we get an induced map $\o\psi :
\Fil^r \cT(\cM) \to \Fil^r \wh A$. Also, as $\varphi_r$ kills the kernel of the map $\widehat{A}_{\mathrm{st}}/p \to \wh A$,
we get an induced map $\o\varphi_r : \Fil^r \wh A \to \widehat{A}_{\mathrm{st}}/p$. The composite $\o\varphi_r \circ \o\psi$ 
induces a map $\Fil^r \cT(\cM)/u^e \Fil^r \cT(\cM) \to \widehat{A}_{\mathrm{st}}/p$ that is linear with respect to the
homomorphism $\varphi: k[u]/u^e \to S_{W(k)}/p$. 
Extending scalars, and using the isomorphisms of Corollary~\ref{cor:brmod-fil}, we get an $S_{W(k)}$-linear map
$\psi : \cM \to \widehat{A}_{\mathrm{st}}/p$. It is now straightforward to check that $\psi$ lifts $\o\psi$ and that
$\psi$ is compatible with all structures. (We note that the results of \cite{Caruso-Crelle} that we cite do not depend
on the assumption $er < p-1$ made there.)

The lower triangle is commutative by an analogous argument.
\end{proof}

Recall the rings $k[[\underline{p}]]$, $k[[\underline{\varpi}]]$ introduced in \S \ref{appendix: Modules}.
The theory of the field of norms \cite{wintenberger} lets us endow the ring 
$k[[\underline{\varpi}]]^s\cong k[[\underline{p}]]^{s}$ with an action of $G_{(K_0)_\infty}$ in such a way that the natural morphism 
\begin{equation}
\label{period ring1}
k[[\underline{p}]]^{s}\to A_{\mathrm{cris}}/p \end{equation} 
becomes $G_{(K_0)_{\infty}}$-equivariant (cf.\ \cite{breuil-normes}, \S 4.2, where $k[[\underline{p}]]^{s}$ is noted by $\cX^{s}$; note that we need to take the composite of $\iota: \cX^{s}\to R$ in \emph{loc.\ cit.} with the natural map $R \to R^{DP} \cong A_{\mathrm{cris}}/p$).
We have an exact, contravariant functor $\mathrm{Hom}(-,k[[\underline{p}]]^s)$ toward Galois representations (cf.\ \cite{breuil-normes}, \S 2.3, \cite{caruso-liu}, \S 2.1).
\begin{prop}
\label{proposition comparison classic1}
The morphism \(\ref{period ring1}\) and the equivalence $M_{\mathfrak{S}}$ induce a commutative diagram
\begin{equation*}
\xymatrix@=5pc{
\FBrModN\ar_-{\text{\tiny$\mathrm{Hom}(-,\widehat A_0)$}}[r]
&\mathrm{Rep}_{\F}(G_{K_{\infty}})\\
\F\text{-}\mathfrak{Mod}^r_{k[[\und{\varpi}]]}\ar^-{M_{\mathfrak{S}}}_-\wr[u]\ar_-{\ \ \text{\tiny$\mathrm{Hom}(-,k[[\underline{p}]]^s)$}}[ur]
}
\end{equation*}
\end{prop}
\begin{proof}
By Proposition~\ref{proposition comparison classic} this follows from \cite{breuil-normes}, Proposition 4.2.1 (which is stated for $e=1$ but whose proof generalizes line by line to the case $e>1$).
\end{proof}

We consider finally the category of $(\phi,k((\underline{\varpi})))$-modules. By a classical result of Fontaine \cite{fontaine-fest}, we have an anti-equivalence 
\begin{align*}
\mathfrak{Mod}_{k((\underline{\varpi}))}&\stackrel{\sim}{\longrightarrow}
\mathrm{Rep}_{\fp}(G_{K_{\infty}})\\
\mathfrak{D}&\longmapsto\mathrm{Hom}(\mathfrak{D},k((\underline{p}))^{s})
\end{align*}
and it is an exercise to prove that the composite functor
$$
\mathfrak{Mod}_{k((\underline{\varpi})),\,\mathrm{dd}}\rightarrow\mathfrak{Mod}_{k((\underline{\varpi}))}\rightarrow
\mathrm{Rep}_{\fp}(G_{K_{\infty}})
$$
factors naturally through the restriction $\mathrm{Rep}_{\fp}(G_{(K_0)_{\infty}})\rightarrow \mathrm{Rep}_{\fp}(G_{K_{\infty}})$.
Indeed, for any $g\in G_{(K_0)_{\infty}}$, $f\in \mathrm{Hom}(\mathfrak{D},k((\underline{p}))^s)$ we define the element
\begin{equation}
\label{descent galois action}
g\cdot f\defeq g\circ f\circ \widehat{\overline{g}}^{-1},
\end{equation}
where $\widehat{\overline{g}}$ denotes the endomorphism of $\mathfrak{D}$ associated to $\overline{g}\in\Gal(K_{\infty}/(K_0)_{\infty})$ via the descent data. It is easy to see that the assignment (\ref{descent galois action}) is well-defined and endows the former $\mathrm{Hom}$-space with a continuous $G_{(K_0)_{\infty}}$-action, thus providing us with the claimed factorization of functors.

\begin{lm}
\label{lemma base change}
The base change along $k((\underline{p}))\into k((\underline{\varpi}))$ induces a commutative diagram of equivalences of categories
$$
\xymatrix{
&\F\text{-}\mathfrak{Mod}_{k((\underline{\varpi})),\,\mathrm{dd}}\ar^-\cong[dl]\\
\mathrm{Rep}_{\F}(G_{(K_0)_{\infty}})&\\
&\F\text{-}\mathfrak{Mod}_{k((\underline{p}))}\ar_-\cong[ul]\ar_-{-\otimes_{k((\underline{p}))}k((\underline{\varpi}))}^-\cong[uu]
}
$$
\end{lm}
\begin{proof}
We assume that $\F=\Fp$.
For commutativity one just observes that for any $(\phi,k((\underline{p})))$-module $\mathfrak{D}_0$,
the natural isomorphism
$$
\mathrm{Hom}_{(\phi,k((\underline{p})))}(\mathfrak{D}_0,k((\underline{p}))^s)\stackrel{\sim}{\longrightarrow}
\mathrm{Hom}_{(\phi,k((\underline{\varpi})))}
(\mathfrak{D}_0\otimes_{k((\underline{p}))}k((\underline{\varpi})),k((\underline{p}))^s)
$$
is $G_{(K_0)_{\infty}}$-equivariant once we endow the right-hand side with the $G_{(K_0)_{\infty}}$-action induced by the descent data on $\mathfrak{D}_0\otimes_{k((\underline{p}))}k((\underline{\varpi})).$ 
The vertical arrow is an equivalence by Galois descent for vector spaces.
\end{proof}
As in the case of $(\phi,k((\underline{\varpi})))$-modules, the contravariant functor 
$\mathfrak{Mod}^r_{k[[\underline{\varpi}]],\,\mathrm{dd}}\rightarrow \mathfrak{Mod}^r_{k[[\underline{\varpi}]]}\rightarrow \mathrm{Rep}_{\fp}(G_{K_{\infty}})$ factors through the restriction $\mathrm{Rep}_{\fp}(G_{(K_0)_{\infty}})\rightarrow \mathrm{Rep}_{\fp}(G_{K_{\infty}})$ and we have the following lemma.

\begin{lm}
\label{lemma galois dd}
The functor $\BrModdd\rightarrow \BrModddN\stackrel{\sim}{\longrightarrow}\mathfrak{Mod}^r_{k[[\underline{\varpi}]],\,\mathrm{dd}}$ fits into a commutative diagram:
$$
\xymatrix{
\F\text{-}\BrModdd
\ar_-{\Tst^*=\mathrm{Hom}(-,\widehat{A})}[d]\ar[r]
&\F\text{-}\mathfrak{Mod}_{k[[\underline{\varpi}]],\,\mathrm{dd}}^r\ar^-{\mathrm{Hom}(-,k[[\underline{p}]]^s)}[d]\\
\mathrm{Rep}_{\F}(G_{K_0})\ar^-{\res}[r]&\mathrm{Rep}_{\F}(G_{(K_0)_{\infty}}).
}
$$
\end{lm}
\begin{proof}
As usual we prove the result when $\F=\Fp$.
Combining Propositions \ref{proposition comparison classic} and \ref{proposition comparison classic1} we have a commutative diagram (the arrows being the evident ones):
\begin{equation*}
\xymatrix{
\BrModdd\ar[r]\ar[d]&\BrModddN\ar[d]&
\mathfrak{Mod}_{k[[\underline{\varpi}]],\,\mathrm{dd}}^{r}\ar[d]\ar_{\sim}^{M_{\mathfrak{S}}}[l]\\
\BrMod\ar[r]\ar[d]&\BrModN\ar[d]
&\mathfrak{Mod}^{r}_{k[[\underline{\varpi}]]}\ar[d]\ar^{\sim}_{M_{\mathfrak{S}}}[l]\\
\mathrm{Rep}_{\fp}(G_{K})\ar^-{\res}[r]&\mathrm{Rep}_{\fp}(G_{K_{\infty}})\ar@{=}[r]&\mathrm{Rep}_{\fp}(G_{K_{\infty}})
}
\end{equation*}
where the natural functorial isomorphism on the right side
\begin{equation}
\label{functorial iso galois}
\mathrm{Hom}(M_{\mathfrak{S}}(-),\widehat{A}_0)\stackrel{\sim}{\longrightarrow}
\mathrm{Hom}(-,k[[\underline{p}]]^s)
\end{equation}
 is induced by the morphisms between the period rings described by (\ref{period rings}), (\ref{period ring1}).
By using twists as in the proof of Lemma~\ref{lemma equivalence dd} and the fact that
the morphisms between period rings are $G_{(K_0)_{\infty}}$-equivariant, it follows that (\ref{functorial iso galois})
is compatible with the $G_{(K_0)_{\infty}}$-action induced by the descent data assignment (\ref{descent galois action}).
An analogous argument applies to the left side of the diagram.
\end{proof}
In a similar fashion, we obtain the following.
\begin{lm}
\label{lemma comparison final}
The localization functor 
$\F\text{-}\mathfrak{Mod}^r_{k[[\underline{\varpi}]],\,\mathrm{dd}}\rightarrow 
\F\text{-}\mathfrak{Mod}_{k((\underline{\varpi})),\,\mathrm{dd}}$ induces a commutative diagram
$$
\xymatrix{
\F\text{-}\mathfrak{Mod}^r_{k[[\underline{\varpi}]],\,\mathrm{dd}}\ar[r]\ar[dr]&
\F\text{-}\mathfrak{Mod}_{k((\underline{\varpi})),\,\mathrm{dd}}\ar_-\cong[d]\\
&\mathrm{Rep}_{\F}(G_{(K_0)_{\infty}}).
}
$$
\end{lm}
\begin{proof}
Without descent data this is \cite{breuil-normes}, Lemme 2.3.3 (the latter being stated for $e=1$ but its proof generalizes). The statement with descent data follows now as in the proof of Lemma \ref{lemma galois dd} (noting that $k[[\und p]]^s\rightarrow k((\underline{p}))^s$ is obviously $G_{(K_0)_{\infty}}$-equivariant).
\end{proof}

\subsection{On dualities}
\label{subsection:dualities}
We recall that we have a notion of duality on the category $\FBrModdd$; more precisely, there is an involutive contravariant functor
\begin{align*}
\FBrModdd&\rightarrow \FBrModdd\\
\cM&\mapsto\cM^{\ast}
\end{align*}
such that 
$
\Tst^{r}(\cM)\cong \Tst^*(\cM^{\ast}).
$
For details, see \cite{EGH}, Definition 3.2.8 (building on work of Caruso \cite{xavier-thesis}, Chapitre V).

\subsection{Proofs of some results in Section \texorpdfstring{\ref{sec:p-adic-hodge}}{2.2}}
\label{appendix: results}

\begin{proof}[Proof of Proposition \ref{proposition comparison}]
As usual we may assume that $\F=\Fp$.
By Lemmas \ref{lemma base change}, \ref{lemma galois dd}, and \ref{lemma comparison final} the commutativity of the upper square and of the right triangle in the statement of Proposition \ref{proposition comparison} is clear. Concerning the lower square, we consider the fully faithful functor $\cF_{p-2} : \mathcal{FL}^{[0,p-2]} \to \BrMod[p-2]$ defined in \cite{breuilAENS98}, \S2.4,
which extends the functor of the same name considered in \S\ref{appendix: Modules}. On the other hand, we have the functor $M_{k((\und p))} : \BrMod[p-2] \to \mathfrak{Mod}_{k((\underline{p}))}$ obtained by specializing the functor $M_{k((\und\varpi))}$ of \S\ref{appendix: Modules} to the case when $e = 1$. The diagram
$$
\xymatrix{
& \mathrm{Rep}_{\fp}(G_{K_0}) \ar^-{\res}[r]& \mathrm{Rep}_{\fp}(G_{(K_0)_{\infty}}) \\
\mathcal{FL}^{[0,p-2]}\ar^-{\cF_{p-2}}[r]\ar^{\Tcris^*}[ur]&\BrMod[p-2] \ar^-{M_{k((\und p))}}[r]\ar[u] &\mathfrak{Mod}_{k((\underline{p}))} \ar^-{\cong}[u]
}
$$
commutes: for the triangle, see \cite{breuilAENS98}, Proposition 3.2.1.1, whereas the square is the specialization of the upper
square of the statement of Proposition \ref{proposition comparison} to the 
case when $e = 1$. The composite morphism $\mathcal{FL}^{[0,p-2]}\rightarrow \mathrm{Rep}_{\fp}(G_{(K_0)_{\infty}})$ 
is fully faithful, as the functor $\cF$ in \S\ref{appendix: Modules} has this property.
\end{proof}

\begin{proof}[Proof of Lemma~\ref{lemma generators dd}]
Lemma~\ref{lm:brmod-fil} shows that $\varphi_r$ induces a $\Gal(K/K_0)$-equivariant and $\varphi\otimes 1$-semilinear isomorphism $\Fil^r\cM/u\Fil^r\cM\stackrel{\sim}{\longrightarrow}\cM/u\cM$. We conclude by Nakayama's lemma and the semisimplicity of  $(k\otimes_{\Fp} \F)[\Gal(K/K_0)]$.
\end{proof}

\begin{proof}[Proof of Lemma~\ref{lemma lawdd 1}]
By passing to the new basis $\und e \cdot A$ we may assume, without loss of generality, that $A = 1$.
As $u^{er} \cM \subseteq \Fil^r \cM$ we deduce that $u^{er} = V V'$ in $\mathrm{M}_{n}(\barS)$ for some $V' \in \mathrm{M}_{n}(\barS)$.
Letting $\wh V \in \mathrm{M}_{n}(k\otimes_{\Fp}\F[[\underline{\varpi}]])$ denote any lift of $V$ such that 
$\wh V_{ij} \in (k\otimes_{\Fp}\F[[\underline{\varpi}]])_{\omega_{\varpi}^{p^{-1}a_j-a_i}}$ we see that
$\und\varpi^{er} = \wh V \wh W = \wh W \wh V$ for some $\wh W \in \mathrm{M}_{n}(k\otimes_{\Fp}\F[[\underline{\varpi}]])$.
Let $\und e^*$ denote the basis of $\cM^*$ that is dual to $\und e$, and let $\und f^* \defeq \und e^* \cdot {}^t W$, where
$W \in \mathrm{M}_{n}(\barS)$ denotes the reduction of $W$.
One easily computes that $f_i^* \in \Fil^r \cM$ and $\varphi_r(f_i^*) = e_i^*$ for all $i$.
As $\varphi_r$ induces an isomorphism $\Fil^r\cM^*/u\Fil^r\cM^*\stackrel{\sim}{\longrightarrow}\cM^*/u\cM^*$ (cf.\ the proof of Lemma~\ref{lemma generators dd}),
it follows that $\und f^*$ generates $\Fil^r \cM$ as $\barS$-module. Also, $\wh g e_i^* = (\omega_{\varpi}^{-a_i}(g)\otimes 1) e_i^*$.

Define $\fM \in \F\text{-}\mathfrak{Mod}^r_{k[[\underline{\varpi}]],\,\mathrm{dd}}$ of rank $n$ by giving it basis $\und{\mf e}$ and defining the maps $\phi$, $\wh g$ by 
$\mathrm{Mat}_{\underline{\mathfrak{e}}}(\phi) = {}^t \wh V$ and $\wh g \mf e_i = (\omega_{\varpi}^{-p^{-1} a_i}(g)\otimes 1) \mf e_i$.
(It is easy to check that $\fM$ is of height $\le r$, and the condition on $\wh V_{ij}$ above implies that $\wh g \circ \phi = \phi \circ \wh g$.)

We can compute $\cM' \defeq M_{\mathfrak{S}} (\fM)$ without coefficients and descent data by the recipe in \cite{caruso-liu}, \S2.1, so 
$\cM' \cong \barS_k \otimes_{\varphi, k[[\und\varpi]]} \fM$. It has $\barS$-basis $\und e' \defeq 1 \otimes_{\varphi} \und {\mf e}$.
By functoriality, $\cM'$ has an $\F$-action with $\lambda \in \F$ acting as $1 \otimes \lambda$. By the twisting argument of the proof of
Lemma \ref{lemma equivalence dd}, $\cM'$ has descent data with $g \in \Gal(K/K_0)$ acting as $g \otimes \wh g$.
Letting $\und f' \defeq \und e' \cdot {}^t W$, it is easy to check that $f_i' \in \Fil^r \cM'$ and $\varphi_r(f_i') = e_i'$ for all $i$.
In particular, as above, we deduce that $\und f'$ generates $\Fil^r \cM'$ as $\barS$-module. Also, 
$\wh g e_i' = (\omega_{\varpi}^{-a_i}(g)\otimes 1) e_i'$. Hence $\cM' \cong \cM^*$ in $\FBrModddN$.
\end{proof}

\begin{proof}[Proof of Lemma~\ref{lemma lawdd 2}]
In the following $e = 1$, so $\barS_k=k[u]/u^p$.
Note that $\cN \defeq\cF_{p-2}(M)\in \FBrModN[p-2]$ has $\barS$-basis $\und e$ (or more precisely $1 \otimes \und e \in \barS_k \otimes_k M$).
A computation shows that $\Fil^{p-2}\cN$ has generating set $\und f$ with $\mat_{\und{e},{\und{f}}}(\Fil^{p-2}\cN)=\mathrm{Diag}(u^{p-2-m_0},\dots,u^{p-2-m_{n-1}})$ 
and $\mat_{\und{e},\und{f}}(\varphi_{p-2})=F$.

We now apply Lemma \ref{lemma lawdd 1} with $e=1$ and $\cM=\cN^*$, noting that $M_{k((\und{p}))}(\cN)=\cF(M)$ in this case.
\end{proof}

\begin{proof}[Proof of Lemma~\ref{lemma base change matrix}]
This is elementary. Indeed, by (\ref{equation base change}), one has
$$
(V+u^{e(r+1)}M)B= A V'
$$
for some element $M\in \mathrm{M}_{n}(\barS)$. By considering descent data, we see that $\und{f}_1\defeq \und{e}\cdot (V+u^{e(r+1)}M)$ is a framed generating family of $\Fil^r\cM$ and one has
$$
\mathrm{Mat}_{\underline{e},\underline{f}_1}(\varphi_r)=\mathrm{Mat}_{\underline{e},\underline{f}}(\varphi_r)=A.
$$
Note that $\und{f}'=\und{e}'\cdot V'=\und{f}_1\cdot B$. As $B\in\GL_n(\barS)$ and by considering descent data, we see that $\und{f}'$ is a framed system of generators for $\Fil^r\cM$.
The last statement follows from  $\varphi_r(\und{f}')=\varphi_r(\und{f}_1)\cdot\varphi(B)=\und{e}\cdot A\varphi(B)=\und{e}'\cdot \varphi(B)$.
\end{proof}

\subsection{Proofs of some results in Section \texorpdfstring{~\ref{sec:breuil-submod}}{2.3}}
\label{appendix: results-fil}

Recall that in Definition \ref{definition:Breuil submodule} we defined the notion of a Breuil submodule.

\begin{proof}[Proof of Lemma \ref{Lemma2Note Florian}]
We note that $\cN$ is a finite free $\barS$-module by Lemma~\ref{Lemma1Note Florian}.
Let $\Fil^r\cN\defeq\Fil^r\cM\cap\cN$, so $u^{er}\cN \subseteq \Fil^r \cN$.
The map $\Fil^r\cN / u^e \Fil^r\cN \to \Fil^r\cM / u^e \Fil^r\cM$ is injective (as $\cN$ is an $\barS_k$-direct summand).
It is in fact a split injection of $k[u]/u^e$-modules by Lemma~\ref{lm:brmod-fil}, as $k[u]/u^e$ is self-injective.
Thus, using again the same lemma, we have the following commutative diagram:
\begin{equation*}
\xymatrix{
 \barS_k \otimes_{\varphi, k[u]/u^e} \Fil^r\cN / u^e \Fil^r\cN \ar@{^{(}->}[r]\ar_{1\otimes \varphi_r}[d]
& \barS_k \otimes_{\varphi, k[u]/u^e} \Fil^r\cM / u^e \Fil^r\cM \ar^{1\otimes \varphi_r}_\cong[d]\\
\cN\ar@{^{(}->}[r]&\cM.
}
\end{equation*}
It follows that the left vertical map is injective, hence surjective by dimension
reasons.
It is now obvious that $\cN$ inherits from $\cM$ the structure of a Breuil module with descent data.

By defining $\Fil^r(\cM/\cN)\defeq (\Fil^r\cM+\cN)/\cN \cong \Fil^r\cM/\Fil^r\cN$ we see that the Frobenius $\varphi_r$ on $\Fil^r\cM$ naturally induces a semilinear morphism $\varphi_r:\Fil^r(\cM/\cN)\rightarrow \cM/\cN$ and it is immediate that the triple $(\cM/\cN, \Fil^r(\cM/\cN),\varphi_r)$ inherits the structure of a Breuil module with descent data.

By construction, the complex $0\rightarrow \cN\rightarrow \cM\rightarrow \cM/\cN\rightarrow 0$ is an exact sequence in $\FBrModdd$. The last statement in the lemma is obvious.
\end{proof}

Recall that given an exact sequence
$$
0\rightarrow \cM_1\stackrel{f}{\rightarrow} \cM\stackrel{g}{\rightarrow} \cM_2\rightarrow0
$$
in $\FBrModdd$, the morphisms $f$ and $g$ are respectively said to be an \emph{admissible monomorphism} and an \emph{admissible epimorphism}.

\begin{lm}
\label{Lemma5Note Florian}
Let $f:\cM\rightarrow \cN$ be a morphism in $\FBrModdd$. Then:
\begin{enumerate}
	\item $f$ is an admissible epimorphism  if and only if $f$ induces $\barS$-linear surjections $\cM\onto \cN$ and $\Fil^r\cM\onto\Fil^r\cN$;
\item $f$ is an admissible monomorphism  if and only if $f$ induces an $\barS$-linear injection $\cM\into \cN$ and  $\Fil^r\cM=f^{-1}(\Fil^r\cN)$.
\end{enumerate}
\end{lm}
\begin{proof}
The ``only if'' part, both in (i) and (ii), is obvious.
Let us assume $f: \cM\rightarrow\cN$ induces $\barS$-linear surjections $\cM\onto \cN$ and $\Fil^r\cM\onto\Fil^r\cN$. Then $\ker(f)$ is clearly a Breuil submodule of $\cM$; by Lemma \ref{Lemma2Note Florian} it is an object in $\FBrModdd$ and the complex
$
0\rightarrow \ker(f)\rightarrow \cM\rightarrow \cN\rightarrow0
$ 
is an exact sequence in $\FBrModdd$.

Assume now that $f$ induces an $\barS$-linear injection $\cM\into \cN$ and  that $\Fil^r\cM=f^{-1}(\Fil^r\cN)$. Then, using Lemma \ref{Lemma1Note Florian} one sees that $f(\cM)$ is a Breuil submodule of $\cN$. As the map $f : \cM \to f(\cM)$ is an isomorphism in $\FBrModdd$, Lemma \ref{Lemma2Note Florian} implies that $f$ is an admissible monomorphism, as required.
\end{proof}

Recall from \S \ref{subsection:dualities} that we have a duality on the category $\FBrModdd$.

\begin{lm}
\label{Lemma6Note Florian}
The functor $\cM\mapsto\cM^{\ast}$ preserves exact sequences in $\FBrModdd$.
\end{lm}
\begin{proof}
By definition of $\Fil^r(\cM^{\ast})$ and the self-injectivity of $\barS/u^{er}$ (cf.\ Lemma \ref{Lemma1Note Florian}) one has an exact sequence
of $\barS$-modules
\begin{equation*}
0\rightarrow \Fil^r(\cM^{\ast})\rightarrow \cM^{\ast}\rightarrow 
\mathrm{Hom}_{\barS/u^{er}}(\Fil^r\cM/u^{er}\cM,\barS/u^{er})\rightarrow 0,
\end{equation*}
which is functorial in $\cM$.
In order to prove the lemma it is enough to prove that the functor $\cM\mapsto \Fil^r(\cM^{\ast})$ is exact (in $\barS$-modules) which is in turn equivalent, thanks to the exact sequence above and the self-injectivity of $\barS/u^{er}$, to the exactness of the functor $\cM\mapsto \Fil^r\cM/u^{er}\cM$ (in $\barS$-modules).
But the exactness of $\cM\mapsto \Fil^r\cM/u^{er}\cM$ is an easy consequence of the snake lemma. \end{proof}

We are now finally in a position to prove Propositions \ref{Proposition4Note Florian} and \ref{Proposition7Note Florian}.

\begin{proof}[Proof of Proposition \ref{Proposition4Note Florian}] 
Clearly $\FBrModdd$ is an additive category and the class of exact sequences is closed under isomorphism. Moreover, given an exact sequence
$$
0\rightarrow \cM_1\stackrel{f}{\rightarrow}\cM\stackrel{g}{\rightarrow}\cM_2\rightarrow0
$$
in $\FBrModdd$ it is clear that $f$ is the categorical kernel of $g$. By duality, $g$ is the categorical cokernel of $f$.

We shall now check the axioms $\text{Ex}0,\,\text{Ex}1,\,\text{Ex}2,\,\text{Ex}2^{op}$ of \cite{keller}, Appendix A.
\begin{itemize}
	\item[$(\text{Ex}0)$:] It is obvious as $0\rightarrow 0$ is an admissible epimorphism.
	\item[$(\text{Ex}1)$:] By Lemma \ref{Lemma5Note Florian} a composition of admissible epimorphisms is an admissible epimorphism.
	\item[$(\text{Ex}2)$:] Let $\alpha:\cM\rightarrow \cM_2$ be an admissible epimorphism and $\beta:\cN\rightarrow\cM_2$ a morphism in $\FBrModdd$. By Lemma \ref{Lemma5Note Florian} the map $\cM\oplus\cN\rightarrow\cM_2$, $(x,y)\mapsto \alpha(x)-\beta(y)$ is an admissible epimorphism, so its kernel $\cP$ is a Breuil submodule of $\cM\oplus\cN$. The same lemma shows moreover that the natural map $\cP\rightarrow\cN$ is an admissible epimorphism. It is then clear that $\cP$ is a categorical pullback of $(\alpha,\beta)$. 	\item[$(\text{Ex}2^{op})$:] It follows from  $(\text{Ex}2)$ via the duality functor on $\FBrModdd$.
\end{itemize}

It remains to prove the exactness of the functor $\Tst^{r}$. Let 
\begin{equation}
\label{exact sequence Florian Note}
0\rightarrow \cM_1\stackrel{f}{\rightarrow}\cM\stackrel{g}{\rightarrow}\cM_2\rightarrow0
\end{equation}
be an exact sequence in $\FBrModdd$. Then, by forgetting monodromy, descent data and coefficients,  we deduce that the exact sequence (\ref{exact sequence Florian Note}) arises from a sequence
$0\rightarrow \mathfrak{M}_1\rightarrow\mathfrak{M}\rightarrow\mathfrak{M}_2\rightarrow0$
in $\mathfrak{Mod}^r_{k[[\underline{\varpi}]]}$
via the equivalence $M_{\mathfrak{S}}$. By \cite{caruso-liu}, Proposition 2.3.2 the latter sequence is exact.
By Proposition \ref{proposition comparison classic1} and the exactness of the functor $\mathrm{Hom}(-,k[[\underline{p}]]^s)$ we deduce that 
$
0\rightarrow \Tst^{r}(\cM_1)\rightarrow \Tst^{r}(\cM)\rightarrow\Tst^{r}(\cM_2)\rightarrow 0
$
is an exact sequence in $\mathrm{Rep}_{\fp}(G_{K_{\infty}})$.
\end{proof}

\begin{proof}[Proof of Proposition \ref{Proposition7Note Florian}]

Let $T'\subseteq\Tst^{r}(\cM) \cong \Tst^*(\cM^{\ast})$. By the proof of \cite{EGH}, Proposition 3.2.6 we get a unique (up to isomorphism) surjective homomorphism $f:\cM^{\ast}\rightarrow \cN'$ in $\FBrModdd$ such that $\mathrm{im}(\Tst^*(f))=T'$.
We claim that $f$ is an admissible epimorphism. We recall from \cite{EGH}, proof of Proposition 3.2.6 and \cite{caruso-fourier}, proof of Proposition 2.2.5 that $f$ is obtained from a surjective homomorphism $\mathfrak{f}:\mathfrak{M}^{\ast}\rightarrow \mathfrak{N}'$ in $\mathfrak{Mod}^r_{k[[\underline{\varpi}]]}$ via the equivalence $M_{\mathfrak{S}}$,
and that we have an exact sequence in $\mathfrak{Mod}^r_{k[[\underline{\varpi}]]}$ given by
$
0\rightarrow\ker(\mathfrak{f})\rightarrow \mathfrak{M}^{\ast}\rightarrow\mathfrak{N}'\rightarrow 0.
$
By the exactness of $M_{\mathfrak{S}}$ we deduce that $f$ is indeed an admissible epimorphism.
By Lemma~\ref{Lemma6Note Florian} the image of $(\cN')^{\ast}$ inside $\cM$ is the desired Breuil submodule in $\cM$ mapping to $T'$ via $\Theta$. It is unique by construction.

We still need to check that $\Theta$ is order preserving. Suppose that $\cM_2$, $\cM_1$ are Breuil submodules of $\cM$ such that $\cM_2\subseteq \cM_1$. Then by Lemma \ref{Lemma3Note Florian} we have $\Theta(\cM_2)\subseteq \Theta(\cM_1)$. Moreover we have a natural isomorphism $\Theta(\cM_1)/\Theta(\cM_2)\cong\Tst^{r}(\cM_1/\cM_2)$ by the exactness of $\Tst^{r}$.
Conversely, given $T_2\subseteq T_1\subseteq \Tst^{r}(\cM)$ then, by the above, we have $T_1=\mathrm{im}(\Tst^{r}(f_1))$ for some admissible monomorphism $f_1:\cM_1\into\cM$ and, similarly, $T_2=\mathrm{im}(\Tst^{r}(f_2))$ for some admissible monomorphism $f_2:\cM_2\into\cM_1$. Then $f_1\circ f_2$ is an admissible monomorphism and, by uniqueness, we obtain finally
\begin{equation*}
\Theta^{-1}(T_2)=\mathrm{im}(f_1\circ f_2)\subseteq \mathrm{im}(f_1)=\Theta^{-1}(T_1).\qedhere
\end{equation*}
\end{proof}

\end{appendices}

\bibliography{Biblio}
\bibliographystyle{amsalpha} 

\end{document}